\documentclass[final,onefignum,onetabnum,onealgnum]{siamart220329}
\usepackage{lipsum}
\usepackage{epstopdf}
\ifpdf
  \DeclareGraphicsExtensions{.eps,.pdf,.png,.jpg}
\else
  \DeclareGraphicsExtensions{.eps}
\fi

\newsiamremark{remark}{Remark}
\newsiamremark{hypothesis}{Hypothesis}
\crefname{hypothesis}{Hypothesis}{Hypotheses}
\newsiamthm{claim}{Claim}
\newsiamthm{assumption}{Assumption}

\headers{Sketched Newton-Raphson}{Rui Yuan, Alessandro Lazaric, and Robert M.~Gower}

\title{Sketched Newton-Raphson\thanks{
Submitted to the editors February 8, 2021; accepted for SIAM Journal on Optimization publication (in revised form) by Coralia Cartis May 04, 2022. \\ 
\indent \indent An earlier version of this work has appeared in the Workshop on "Beyond first-order methods in ML systems" at the 37th International Conference on Machine Learning, Online, 2020.
}}

\author{Rui Yuan\thanks{Meta AI; LTCI, T\'el\'ecom Paris; Institut Polytechnique de Paris (\email{ruiyuan@fb.com}).}
\and Alessandro Lazaric\thanks{Meta AI (\email{lazaric@fb.com}).}
\and Robert M.~Gower\thanks{CCM, Flatiron Institute; Part of this work was done when author was affiliated with T\'el\'ecom Paris and Meta AI
(\email{gowerrobert@gmail.com}).}}

\usepackage{amsopn}

\ifpdf
\hypersetup{
  pdftitle={Sketched Newton-Raphson},
  pdfauthor={Rui Yuan, Alessandro Lazaric, and Robert M.~Gower}
}
\fi

\usepackage{layout}
\usepackage{multirow}

\usepackage{amsfonts}
\usepackage{amsmath}
\usepackage{amsbsy}
\usepackage{xcolor}
\usepackage{color}
\usepackage{graphicx}

\usepackage{hyperref}
\usepackage{textcomp}
\usepackage{cases}  

\newcommand{\R}{\mathbb{R}} 
\newcommand{\N}{\mathbb{N}} 

\newcommand{\cB}{{\cal B}}

\newcommand{\cD}{{\cal D}}

\newcommand{\cL}{{\cal L}}

\newcommand{\cN}{{\cal N}}
\newcommand{\cO}{{\cal O}}

\newcommand{\cY}{{\cal Y}}

\newcommand{\mA}{{\bf A}}

\newcommand{\mD}{{\bf D}}

\newcommand{\mG}{{\bf G}}
\newcommand{\mH}{{\bf H}}
\newcommand{\mI}{{\bf I}}

\newcommand{\mM}{{\bf M}}

\newcommand{\mS}{{\bf S}}

\newcommand{\mW}{{\bf W}}

\usepackage[colorinlistoftodos,bordercolor=orange,backgroundcolor=orange!20,linecolor=orange]{todonotes}

\newcommand{\robline}[1]{#1}

\newcommand{\ruiline}[1]{{#1}}
\newcommand{\ruilinex}[1]{{#1}}
\newcommand{\xline}[1]{{#1}}

\newcommand{\eqdef}{\overset{\text{def}}{=}} 

\newcommand{\dotprod}[1]{\left< #1\right>} 
\newcommand{\norm}[1]{ \| #1 \|}      

\newcommand{\Prob}[1]{\mathbb{P}[#1]}

\providecommand{\Image}[1]{{\bf Im}\left( #1\right)}

\DeclareMathOperator{\argmin}{argmin}        

\newcommand{\Diag}[1]{\mathbf{Diag}\left( #1\right)}
\providecommand{\null}[1]{{\bf Null}\left( #1\right)}

\newcommand{\E}[1]{\mathbb{E}\left[#1\right] } 
\newcommand{\EE}[2]{\mathbb{E}_{#1}\left[#2\right] }

\newcommand{\colvec}[1]{\begin{bmatrix}#1\end{bmatrix}}

\usepackage{booktabs}       

\usepackage{xspace}
\usepackage{algorithm}
\usepackage[noend]{algpseudocode}

\newcommand{\SNR}{\texttt{SNR}\xspace}
\newcommand{\SGD}{\texttt{SGD}\xspace}
\newcommand{\NR}{\texttt{NR}\xspace}
\newcommand{\GN}{\texttt{GN}\xspace}
\newcommand{\SNM}{\texttt{SNM}\xspace}
\newcommand{\RSN}{\texttt{RSN}\xspace}

\begin{document}


\maketitle

\begin{abstract}
  We propose a new globally convergent stochastic second order method. Our starting point is the development of a new Sketched Newton-Raphson (\SNR) method for solving large scale nonlinear equations of the form $F(x)=0$ with $F:\R^p \rightarrow \R^m$.
  We then show how to design several stochastic second order optimization methods by re-writing the optimization problem of interest as a system of nonlinear equations and applying \SNR. For instance, by applying \SNR to find a stationary point of a generalized linear model (GLM), we derive completely new and scalable stochastic second order methods. We show that the resulting method is very competitive as compared to state-of-the-art variance reduced methods.
  Furthermore, using a variable splitting trick, we also show that the \emph{Stochastic Newton method} (\SNM) is a special case of \SNR, and use this connection to establish the first global convergence theory of \SNM.

  We establish the global convergence of \SNR by showing that it is
  a variant of the online stochastic gradient descent (\SGD) method, and then leveraging proof techniques of \SGD.
  As a special case, our theory also provides a new global convergence theory for the original Newton-Raphson method under strictly weaker assumptions as compared to the classic monotone convergence theory.
\end{abstract}

\begin{keywords}
  Nonlinear systems, stochastic methods, iterative methods, stochastic Newton method, randomized Kaczmarz, randomized Newton, randomized Gauss-Newton, randomized fixed point, randomized subspace-Newton.
\end{keywords}

\begin{MSCcodes}
  58C15, 90C06, 90C53, 
  62L20, 46N10, 46N40, 49M15, 68W20, 
  68W40, 65Y20
\end{MSCcodes}
 
\begin{ACM}
  G.1.6 Optimization
\end{ACM}

\tableofcontents


\section{Introduction}

One of the fundamental problems in numerical computing is to find roots of systems of nonlinear equations such as
\begin{eqnarray}\label{eq:main}
F(x)  =  0,
\end{eqnarray}
where $F: \R^p \rightarrow \R^m$. \ruiline{We assume throughout that $F: \R^p \rightarrow \R^m$ is continuously differentiable and that there exists a solution to~\eqref{eq:main}, that is
\begin{assumption} \label{ass:zero}
$\exists x^* \in \R^p$ such that $F(x^*) = 0$.
\end{assumption}}
 This includes a wide range of applications
from solving the phase retrieval problems~\cite{PhaseRetrieval2015}, 
systems of polynomial equations related to cryptographic primitives~\cite{SysPolynomial2019}, 
discretized integral and differential equations~\cite{Ortega:2000}, 
the optimal power flow problem~\cite{871729} and,
our main interest here, solving nonlinear minimization problems in machine learning. Most convex optimization problems such as those arising from training a Generalized Linear Model (GLM), can be re-written as a system of nonlinear equations~\eqref{eq:main} either by manipulating the stationarity conditions or as the Karush-Kuhn-Tucker equations\footnote{Under suitable constraint qualifications~\cite{Nocedal1999a}.}. 

When dealing with non-convex optimization problems, such as training a Deep Neural Network (DNN), finding the global minimum is often infeasible (or not needed \cite{poorlocalmin}). Instead, the objective is to find a good stationary point $x$ such that $\nabla f(x) =0,$ where $f$ is the total loss we want to minimize. 

In particular, the task of training an overparametrized DNN (as they often are) can be cast as solving a special nonlinear system. That is, when the DNN is sufficiently overparametrized, the DNN can interpolate the data. As a consequence,  if $f_i(x)$ is the loss function over the $i$th data point, then there is a solution to the system of nonlinear equations $\norm{ \nabla f_i(x)}^2 =0,  \forall i.$

The building block of many iterative methods for solving nonlinear equations is the Newton-Raphson (\NR) method given by
\begin{align}
x^{k+1} &= x^k - \gamma\left(DF(x^k)^\top\right)^\dagger F(x^k)  \label{eq:newton}
\end{align}
at $k$th iteration, where $DF(x) \eqdef [\nabla F_1(x) \ \cdots \ \nabla F_m(x)] \in \R^{p \times m}$ is \ruiline{the transpose of} the Jacobian matrix of $F$ at $x$,  $\left(DF(x^k)^\top\right)^\dagger$ is the Moore-Penrose pseudoinverse of $DF(x^k)^\top$  and $\gamma>0$ is the stepsize.

The \NR method is at the heart of many commercial solvers for nonlinear equations~\cite{Ortega:2000}. The success of \NR can be partially explained by its invariance to affine coordinate transformations, which in turn means that the user does not need to tune any parameters (standard \NR sets $\gamma=1$). The downside of \NR is that we need to solve a linear least squares problem 
given in~\eqref{eq:newton} which costs
$\cO(\min\{pm^2, mp^2\})$ when using a direct solver. When both $p$ and $m$ are large, this cost per iteration is prohibitive. Here we develop a randomized \NR method based on the sketch-and-project technique~\cite{Gower2015} which can be applied in large scale, as we show in our experiments.


\subsection{The sketched Newton-Raphson method}

Our method relies on using \emph{sketching matrices} to reduce the dimension of the Newton system.

\begin{definition}
The sketching matrix $\mS \in \R^{m \times \tau}$ is a random matrix sampled from a distribution $\cD$, where $\tau \in \N$ is the sketch size. 
\ruiline{We use $\mS_k\in \R^{m \times \tau}$ to denote a sketching matrix sampled from a distribution $\cD_{x^k}$ that can depend  on the iterate $x^k.$}
\end{definition}

By sampling a sketching matrix $\mS_k \sim \cD_{x^k}$ at $k$th iteration, we \emph{sketch} (row compress) \NR update and compute an
approximate \emph{Sketched Newton-Raphson} (\SNR) step, see~\eqref{eq:update} in Alg.~\ref{algo:SkeNeR}.
We use $\cD_x$ to denote a distribution that depends on $x$, and allow the distribution of the sketching matrix to change from one iteration to the next.
\begin{algorithm}
\caption{ \SNR: Sketched Newton-Raphson}\label{algo:SkeNeR}
\begin{algorithmic}[1]
\State \textbf{parameters:} $\cD =$ distribution of sketching matrix; stepsize parameter $\gamma > 0$
\State \textbf{initialization:} Choose $x^0 \in \R^p$
\For{$k = 0, 1, \cdots $}
  \State Sample a fresh sketching matrix:  $\mS_k \sim \cD_{x^k}$ 
  \begin{eqnarray}\label{eq:update}
   x^{k+1} &=& x^k - \gamma DF(x^k)\mS_k \left(\mS_k^\top  DF(x^k)^\top DF(x^k)\mS_k \right)^\dagger\mS_k^\top F(x^k)  \hfill
  \end{eqnarray}
\EndFor
\State \textbf{return:} last iterate $x^k$
\end{algorithmic}
\end{algorithm}

Because the sketching matrix $\mS_k$ has $\tau$ columns, the dominating costs of computing the \SNR step~\eqref{eq:update} are linear in $p$ and $m$. 
\ruiline{
In particular, $DF(x^k)\mS_k \in \R^{p \times \tau}$ can be computed by using $\tau$ directional derivatives of $F(x^k)$, one for each column of $\mS_k$. Using automatic differentiation~\cite{Bruce1992}, these directional derivatives cost $\tau$ evaluations of the function $F(x).$ 
Furthermore, it costs $\cO(p\tau^2)$ to form the linear system in \eqref{eq:update} of Alg.~\ref{algo:SkeNeR} by using the computed matrix $DF(x^k)\mS_k $ and $\cO(\tau^3)$ to solve it , respectively. Finally the matrix vector product $\mS_k^\top F(x^k)$ costs $\cO(m \tau)$. Thus, without making any further assumptions to the structure of $F$ or the sketching matrix, the total cost in terms of operations of the update~\eqref{eq:update} is given by
\begin{equation}
\mbox{Cost(update~\eqref{eq:update})} = \cO \left((\mbox{eval}(F)+m)  \times \tau   + p \tau^2+ \tau^3\right).
\end{equation}
 Thus Alg.~\ref{algo:SkeNeR} can be applied when both $p$ and $m$ are large and $\tau$ is relatively small. \\ }

The rest of the paper is organized as follows. In the next section, we provide some background and contrast it with our contributions. After introducing some notations in Section~\ref{sec:notations} and presenting alternative sketching techniques in Section~\ref{sec:sketch}, we show that~\eqref{eq:update} can be viewed as a sketch-and-project type method in Section~\ref{sec:sketchproj}. This is the viewpoint that first motivated the development of this method. After which, we provide another crucial equivalent viewpoint of~\eqref{eq:update} in Section~\ref{sec:reform}, where we show that \SNR can be seen as \emph{Stochastic Gradient Descent} (\SGD) applied to an equivalent reformulation of~\eqref{eq:main}. We then provide a global convergence theory by leveraging this insight in Section~\ref{sec:convergence}. As a special case, our theory also provides a new global convergence theory for the original \NR method \eqref{eq:newton} under strictly weaker assumptions as compared to the monotone convergence theory in Section~\ref{sec:globalNR}, albeit for different step sizes. For the other extreme where the sketching matrix samples a single row, we present the new nonlinear Kaczmarz method as a variant of \SNR and its global convergence theory in Section~\ref{sec:Kaczmarz}. We then show how to design several stochastic second order optimization methods by re-writing the optimization problem of interest as a system of nonlinear equations and applying \SNR. For instance, using a variable splitting trick, we show that the \emph{Stochastic Newton method} (\SNM)~\cite{rodomanov16superlinearly,Kovalev2019stochNewt} is a special case of \SNR, and use this connection to establish the first global convergence theory of \SNM in Section~\ref{sec:SNM}. In Section~\ref{sec:GLM}, by applying \SNR to find a stationary point of a GLM, we derive completely new and scalable stochastic second order methods. We show that the resulting method is very competitive as compared to state-of-the-art variance reduced methods.

\subsection{Background and contributions}


\paragraph{a) Stochastic second-order methods} There is now a concerted effort to develop efficient second-order methods for solving high dimensional and stochastic optimization problems in machine learning. Most recently developed Newton methods fall into one of two categories: \emph{subsampling} and \emph{dimension reduction}. The subsampling methods~\cite{Erdogdu2015nips,Roosta-Khorasani2016,KohlerL17,Bollapragada2018,pmlr-v80-zhou18d} and~\cite{JMLR:v18:16-491,Pilanci2015a}\footnote{Newton sketch~\cite{Pilanci2015a} and LiSSa~\cite{JMLR:v18:16-491} use subsampling to build an estimate of the Hessian but require a full gradient evaluation. As such, these methods are not efficient for very large $n$.} use mini-batches to compute an approximate Newton direction. Though these methods can handle a large number of \emph{data points} ($n$), they do not scale well in the number of \emph{features} ($d$). On the other hand, second-order methods based on dimension reduction techniques such as~\cite{RSN_nips} apply Newton's method over a subspace of the features, and as such, do not scale well in the number of data points. Sketching has also been used to develop second-order methods in the online learning setting~\cite{G_rb_zbalaban_2015,Agarwal:sketchNewt:2016,NIPS2017_7194} and quasi-Newton methods~\cite{GowerGold2016}.

\noindent \emph{Contributions.}  We propose a new family of stochastic second-order method called \SNR.
Each choice of the sketching distribution and nonlinear equations used to describe the stationarity conditions, leads to a particular algorithm. For instance, we show that a nonlinear variant of the Kaczmarz method is a special case of \SNR. We also show that the subsampling based \SNM~\cite{rodomanov16superlinearly,Kovalev2019stochNewt} is a special case of \SNR. 
By using a different norm in the sketch-and-project viewpoint, we show that the dimension reduced method \emph{Randomized Subspace Newton} (\RSN)~\cite{RSN_nips} is also a special case of \SNR. 
We provide a concise global convergence theory, that when specialized to \SNM gives its first global convergence result. Furthermore, the convergence theory of \SNR
allows for any sketch size, which translates to any mini-batch size for the nonlinear Kaczmarz and \SNM.
In contrast, excluding \SNM, the subsampled based Newton methods~\cite{Erdogdu2015nips,Roosta-Khorasani2016,KohlerL17,Bollapragada2018,pmlr-v80-zhou18d,JMLR:v18:16-491,Pilanci2015a} rely on high probability bounds that in turn require large mini-batch sizes~\footnote{The batch sizes in these methods scale proportional to a condition number~\cite{JMLR:v18:16-491} or $\epsilon^{-1}$ where $\epsilon$ is the desired tolerance.}.
We detail the nonlinear Kaczmarz method in Sec.~\ref{sec:Kaczmarz}, the connection with \SNM in Sec.~\ref{sec:SNM}
 and \RSN in App.~\ref{sec:RSN}.

\paragraph{b) New method for GLMs}  There exist several specialized methods for solving GLMs, including variance reduced gradient methods such as SAG/SAGA~\cite{SAG,SAGA_Nips} and SVRG~\cite{Johnson2013}, and methods based on dual coordinate ascent like SDCA~\cite{SDCA}, dual free SDCA (dfSDCA)~\cite{Shalev-Shwartz2015} and Quartz~\cite{Qu2015b}.

\noindent \emph{Contributions.}  We develop a specialized variant of \SNR for GLMs in Sec.~\ref{sec:GLM}. Our resulting method scales linearly in the number of dimensions $d$ and the number of data points $n$, has the same cost as \SGD per iteration in average.
We show in experiments that our method is very competitive as compared to state-of-the-art variance reduced methods for GLMs.



\paragraph{c) Viewpoints of (Sketched) Newton-Raphson} We show in Sec.~\ref{sec:SGD} that \SNR can be seen as \SGD applied to an equivalent reformulation of our original problem.  We will show that this reformulation is \emph{always} a smooth and interpolated function~\cite{MaBB18,VaswaniBS19}. These gratuitous properties allow us to establish a simple global convergence theory by only assuming that the reformulation is a \emph{star-convex} function: a class of nonconvex functions that include convexity as a special case~\cite{Nesterov2006cubic,LeeV16,zhou2018sgd,hinder2020nearoptimal}. The details of the \SGD interpretation can be found in Sec.~\ref{sec:SGD}. 
In addition, we also show in App.~\ref{sec:viewpoints} that \SNR can be seen as a type of stochastic Gauss-Newton method or as a type of stochastic fixed point method.



\paragraph{d) Classic convergence theory of Newton-Raphson}\ruiline{
The better known convergence theorems for \NR (the Newton-Kantorovich-Mysovskikh Theorems) only guarantee local or semi-local convergence~\cite{kantorovitch1939,Ortega:2000}. To  guarantee global convergence of \NR,
 we often need an additional globalization strategy, such as damping sequences or adaptive trust-region methods~\cite{trustregionbook,Lu2010,DeuflhardNewton:2011,Kelley2018numerical}, continuation schemes such as interior point methods~\cite{nesterov1994interior,Wright2006}, and more recently cubic regularization~\cite{Kovalev2019stochNewt,Nesterov2006cubic,Cartis2009}.
Globalization strategies are used in conjunction with other second-order methods, such as inexact Newton backtracking type methods~\cite{Bellavia2001aglobally,An2007aglobally}, Gauss-Newton or Levenberg-Marquardt type methods~\cite{Zhou2010global,Zhou2013on,Yuan2011recent} and quasi-Newton methods~\cite{Yuan2011recent}\footnote{A recent paper~\cite{Gao2019quasi-newton} shows that quasi-Newton converges globally for self-concordant functions without globalization strategy.}.
The only global convergence theory that does not rely on such a globalization strategy, requires strong assumptions on $F(x)$, such as in the monotone convergence theory (MCT)~\cite{DeuflhardNewton:2011}. }

\noindent \emph{Contributions.} \ruiline{
We show in Sec.~\ref{sec:monotoneconv} that our main theorem specialized to the standard \NR method guarantees a global convergence under \emph{strictly} less assumptions as compared to the MCT, albeit under a different stepsize. Indeed, MCT holds for step size equal to one ($\gamma =1$) and our theory holds for step sizes less than one ($\gamma <1)$.
}

Furthermore, we give  an explicit sublinear $O(1/k)$ convergence rate, as opposed to only an asymptotic convergence in MCT. This appears to not be known before since, as stated by \cite{DeuflhardNewton:2011} w.r.t. the \NR method ``\emph{Not even an a-priori estimation for the number of iterations needed to achieve a prescribed accuracy {\bf may} be possible}''. We show that it is possible by monitoring which iterate achieves the best loss (suboptimality).


\paragraph{e) Sketch-and-project} The sketch-and-project method was originally introduced for solving linear systems in~\cite{Gower2015,Gower2015c}, where it was also proven to converge linearly and globally. In \cite{Reformulation}, the authors then go on to show that the sketch-and-project method is in fact \SGD applied to a particular reformulation of the linear system.

\noindent \emph{Contributions.} 
It is this \SGD viewpoint in the linear setting~\cite{Reformulation} that we extend to the nonlinear setting.
Thus the \SNR algorithm and our theory are generalizations of the original sketch-and-project method for solving linear equations to solving nonlinear equations, thus greatly expanding the scope of applications of these techniques.


\subsection{Notations}\label{sec:notations}

In calculating an update of \SNR~\eqref{eq:update} and analyzing \SNR, the following random matrix is key
\begin{eqnarray}\label{eq:Hk}
\mH_{\mS}(x) &\eqdef& \mS \left(\mS^\top DF(x)^\top DF(x) \mS \right)^\dagger \mS^\top.
\end{eqnarray}
The sketching matrix $\mS$ in~\eqref{eq:Hk} is sampled from a distribution $ \cD_x$ and $\mH_{\mS}(x)\in\R^{m \times m}$ is a random matrix that depends on $x.$ We use $\mI_p \in \R^{p \times p}$ to denote the identity matrix of dimension $p$
and 
\xline{use $\norm{x}_{\mM} \eqdef \sqrt{x^\top \mM x}$ to denote the seminorm of $x \in \R^p$ induced by a symmetric positive semi-definite matrix $\mM\in\R^{p\times p}$. Notice that $\norm{x}_{\mM}$ is not necessarily a norm as $\mM$ is allowed to be non invertible. We handle this with care in our forthcoming analysis.}
We also define the following sets: $F(U) = \{F(x) \mid x \in U\}$ for a given set $U \subset \R^p$; 
\xline{$W^\perp = \{v  \mid \dotprod{u,v} = 0, \ \forall u \in W\}$ to denote the orthogonal complement of a subspace $W$;
}
${\bf Im}(\mM) = \{y \in \R^m \mid \exists x \in \R^p \ \mbox{s.t.\ } \mM x = y \}$ to denote the image space and ${\bf Ker}(\mM) = \{x \in \R^p \mid \mM x = 0\}$ to denote the null space of a matrix $\mM \in \R^{m \times p}$. If $\mM$ is a random matrix sampled from a certain distribution $\cD$, we use $\EE{\mS \sim \cD}{\mM} = \int_{\mM}  \mM d\mathbb{P}_{\cD}(\mM)$ to denote the expectation of the random matrix. We omit the notation of the distribution $\cD$, i.e.\ $\E{\mM}$, when the random source is clear. In particular, when $\mM$ is sampled from a discrete distribution with $r \in \N$ s.t.\
$\Prob{\mM = \mM_i} = p_i > 0, \ \mbox{for } i = 1, \cdots, r$ and $\sum_{i=1}^r p_i = 1$, then $\E{\mM} = \sum_{i=1}^r p_i \mM_i$.

\subsection{Sketching matrices}
\label{sec:sketch}

\ruiline{Here we provide examples of sketching matrices that can be used in conjunction with \SNR. We point the reader to~\cite{woodruff2014sketching}
for a detailed exposure and introduction. The most straightforward sketch is given by the Gaussian sketch where every coordinate $\mS_{ij}$ of the sketch $\mS \in \R^{m\times\tau}$ is sampled i.i.d according to a Gaussian distribution with $\mS_{ij} \sim \cN(0,\frac{1}{\tau})$ for $i = 1, \cdots, m$ and $j = 1, \cdots, \tau$. 
The sketch we mostly use here is the uniform subsampling sketch, whereby
	\begin{equation}
		\label{eq:subsample}
		\Prob{\mS = \mI_{C}} \; = \; \frac{1}{\binom{m}{\tau}} \enspace, \quad \mbox{for all set } C \subset \{1, \cdots, m\} \mbox{ s.t.\ } |C| = \tau \enspace,
	\end{equation}
where $\mI_{C} \in \R^{m \times \tau}$ denotes the concatenation of the columns of the identity matrix $\mI_m$ indexed in the set $C$.
More sophisticated sketches that are able to make use of fast Fourier type routines include the 
random orthogonal sketches (ROS)~\cite{Pilanci2015a,Ailon:2009}.
 We will not cover ROS sketches here since these sketches are fast when applied only once to a fixed matrix $\mM,$ as opposed to being re-sampled at every iteration .}

\section{The sketch-and-project viewpoint}
\label{sec:sketchproj}

The viewpoint that motivated the development of Alg.~\ref{algo:SkeNeR} was the following iterative \emph{sketch-and-project} method applied to the Newton system. \ruiline{For this  viewpoint, we  assume that }
\begin{assumption}\label{ass:existence} 
$F(x)  \in   \Image{DF(x)^\top}   \forall x \in \R^p.$
\end{assumption}
This assumption 
guarantees that there exists a solution to the Newton system in~\eqref{eq:newton}.
Indeed, we can now re-write the \NR method~\eqref{eq:newton} as a projection of the previous iterate $x^k$ onto the solution space of a Newton system 
\begin{align}
x^{k+1} = \argmin_{x\in\R^p} \norm{x - x^k}^2 \quad \mbox{ s. t.}\quad  DF(x^k)^\top(x-x^k) = -\gamma F(x^k). \label{eq:newtonproj}
\end{align}
Since this is costly to solve when $DF(x^k)$ has many rows and columns, we \emph{sketch} the Newton system. That is, we apply a random row compression to the Newton system using the sketching matrix $\mS_k^\top \in \R^{\tau \times m}$ and then project the previous iterates $x^k$ onto this \emph{sketched} system as follows
\begin{align} \label{eq:sketch}
x^{k+1} = \argmin_{x \in \R^p} \norm{x-x^k}^2 \quad \mbox{ s. t.}\quad \mS_k^\top  DF(x^k)^\top(x-x^k)  = -\gamma\mS_k^\top F(x^k). 
\end{align}
That is, $x^{k+1}$ is the projection of $x^k$ onto the solution space of the sketched Newton system. \ruiline{This viewpoint was our motivation for developing the \SNR method. Next we establish our  core theory. The theory does not rely on the assumption 
$F(x)  \in   \Image{DF(x)^\top}$, though this assumption will appear again in several specialized corollaries.
Without this assumption, we can still interpret the Newton step~\eqref{eq:newton} as the least squares solution of the linear system~\eqref{eq:newtonproj}, as we show next.}

%



\section{Reformulation as stochastic gradient descent}\label{sec:SGD}
\label{sec:reform}

Our insight into interpreting and analyzing the \SNR in Alg.~\ref{algo:SkeNeR} is through its connection to the \SGD. Next, we show how \SNR can be seen as \SGD applied to a sequence of equivalent reformulations of~\eqref{eq:main}. 
Each reformulation is given by a vector $y\in\R^p$ and the following minimization problem
%
\begin{equation}\label{eq:sgdreform}
\min_{x \in \R^p} \EE{\mS \sim \cD_y}{\tfrac{1}{2}\norm{F(x)}_{\mH_{\mS}(y)}^2 },
\end{equation}
where 
 $\mH_{\mS}(y)$ is defined in~\eqref{eq:Hk}. To abbreviate notations, let
\begin{equation}
f_{\mS,y}(x)  \eqdef  \frac{1}{2}\norm{F(x)}_{\mH_{\mS}(y)}^2  \quad \mbox{and } \quad 
f_y(x)  \eqdef  \E{f_{\mS,y}(x)}  =  \frac{1}{2}\norm{F(x)}_{\E{\mH_{\mS}(y)}}^2.\label{eq:fk}
\end{equation}

\ruiline{
Every solution $x^*\in\R^p$ to~\eqref{eq:main} is a solution to~\eqref{eq:sgdreform}, since $f_y(x)$ is non-negative for every $x \in \R^p$ and $f_y(x^*) =0 $ is thus a global minima.
With an extra assumption, we can show that every solution to~\eqref{eq:sgdreform} is also a solution to~\eqref{eq:main} in the following lemma.
\begin{lemma}\label{lem:equiv}
If Asm~\ref{ass:zero} holds 
and the following 
\emph{reformulation} 
assumption 
\begin{equation}\label{ass:ker}
F(\R^p)\cap{\bf Ker}\left(\EE{\mS \sim \cD_y}{\mH_{\mS}(y)}\right)=\{0\}, \quad\forall y \in \R^p
\end{equation}
holds, then $\argmin_{x \in \R^p} f_y(x)  =  \left\{x  \mid  F(x) = 0\right\}$ for every $y \in \R^p.$
\end{lemma}}

\begin{proof}
Let $y \in \R^p$. 
Previously, we show that $\left\{x  \mid  F(x) = 0\right\} \subset \argmin_{x \in \R^p} f_y(x)$.
Now let $x^* \in \argmin_{x \in \R^p} f_y(x)$. By Asm.~\ref{ass:zero}, we know that any global minimizer $x^*$ of $f_y(x)$ must be s.t.\ $f_y(x^*) = 0$.
This implies that $F(x^*) \in {\bf Ker}\left(\E{\mH_{\mS}(y)}\right)$ since $\E{\mH_{\mS}(y)}$ is symmetric. However $F(x^*) \in F(\R^p)$ and thus from~\eqref{ass:ker}, we have that $F(x^*) \in F(\R^p) \cap {\bf Ker}\left(\E{\mH_{\mS}(y)}\right) = \{0\}$, which implies $F(x^*) = 0$. Thus, we have $\argmin_{x \in \R^p} f_y(x) \subset \left\{x  \mid  F(x) = 0\right\}$ which concludes the proof. 
\end{proof}

Thus with the extra reformulation assumption in~\eqref{ass:ker}, we can now use any viable optimization method to solve~\eqref{eq:sgdreform} for any fixed $y\in\R^p$  and arrive at a solution to~\eqref{eq:main}. In Lemma~\ref{lem:inv}, 
we give sufficient conditions on the sketching matrix and on the function $F(x)$ that guarantee~\eqref{ass:ker} hold. We also show how~\eqref{ass:ker} holds for our forthcoming examples in Sec.~\ref{sec:suffcondassker} as a direct consequence of Lemma~\ref{lem:inv}. 
 However,~\eqref{ass:ker} imposes for all $y \in \R^p$ which can be sometimes restrictive. 
 In fact, we do \emph{not need} for~\eqref{eq:sgdreform} to be equivalent to solving~\eqref{eq:main} for \emph{every $y\in\R^p$.} Indeed, by carefully and iteratively updating $y$, we can solve~\eqref{eq:sgdreform} and obtain a solution to~\eqref{eq:main} \emph{without} relying on~\eqref{ass:ker}. The trick here is to use an \emph{online \SGD} method for solving~\eqref{eq:sgdreform}.

%
%

%

Since~\eqref{eq:sgdreform} is a stochastic optimization problem, \SGD is a natural choice for solving~\eqref{eq:sgdreform}.
Let $\nabla f_{\mS, y}(x)$ denote the gradient of the function $f_{\mS, y}(\cdot)$ which is 
\begin{eqnarray}\label{eq:SGDgrad}
\nabla f_{\mS, y}(x) &=& DF(x)\mH_{\mS}(y) F(x).
\end{eqnarray}

Since we are free to choose $y$, we allow $y$ to \emph{change} from one iteration to the next by setting $y=x^k$ at the start of the $k$th iteration. We can now take a \SGD step by sampling $\mS_k \sim \cD_{x^k}$ at $k$th iteration and updating
\begin{eqnarray}
x^{k+1} &=& x^k - \gamma \nabla f_{\mS_k,x^k}(x^k)   \label{eq:SGD}.
\end{eqnarray}
It is straightforward to verify that the \SGD update~\eqref{eq:SGD} is exactly the same as the \SNR update in~\eqref{eq:update}.

The objective function $f_{\mS,y}(x)$ has many properties that makes it very favourable for optimization including the interpolation condition and a gratuitous smoothness property. Indeed, for any $x^* \in \R^p$ s.t.\ $F(x^*) =0$, we have that the stochastic gradient is zero, i.e.\ $\nabla f_{\mS,y}(x^*) =0. $ This is known as the \emph{interpolation condition}. When it occurs together with strong convexity, it is possible to 
shows that \SGD converges linearly~\cite{VaswaniBS19,MaBB18}. We will also give a linear convergence result in Sec.~\ref{sec:convergence} by assuming that $f_{y}(x)$ is quasi-strongly convex. 
We detail the smoothness property next.

However, we need to be careful, since~\eqref{eq:SGD} is not a classic \SGD method. In fact, from the $k$th iteration to the $(k+1)$th iteration, we change our objective function from $f_{x^k}(x)$ to $f_{x^{k+1}}(x)$ and the distribution from $\cD_{x^k}$ to $\cD_{x^{k+1}}$. \ruiline{Thus it is an online \SGD.} We handle this with care in our forthcoming convergence proofs.


\section{Convergence theory}\label{sec:convergence}

Using the viewpoint of \SNR in Sec.~\ref{sec:SGD}, we adapt proof techniques of \SGD to establish the global convergence of \SNR.

\subsection{Smoothness property}

In our upcoming proof, we rely on the following type of smoothness property thanks to our \SGD reformulation~\eqref{eq:sgdreform}.
\begin{lemma}\label{lem:smooth}
For every $x\in\R^p$ and any realization $\mS \sim \cD_x$ associated with any distribution $\cD_x$,
\begin{align}\label{eq:1smooth}
\frac{1}{2}\|\nabla f_{\mS, x}(x)\|^2 =  f_{\mS, x}(x).
\end{align}
\end{lemma}

\begin{proof}
Turning to the definition of $f_{\mS,x}$ in~\eqref{eq:fk}, we have that
\begin{align*}
\norm{\nabla f_{\mS,x}(x)}^2 &\overset{\eqref{eq:SGDgrad}}{=} \norm{DF(x)\mH_{\mS}(x)F(x)}^2 
= F(x)^\top \mH_{\mS}(x)^\top DF(x)^\top DF(x)\mH_{\mS}(x) F(x) \\
&= F(x)^\top\mH_{\mS}(x)F(x) = 2f_{\mS,x}(x),
\end{align*}
where we used the property 
 $ \mM^\dagger \mM\mM^\dagger = \mM^\dagger$ with $\mM = \mS^\top DF(x)^\top DF(x) \mS$ to establish that
$\mH_{\mS}(x)^\top DF(x)^\top DF(x)\mH_{\mS}(x) \overset{\eqref{eq:Hk}}{=} \mH_{\mS}(x).$
\end{proof}

This is not a standard smoothness property. 
Indeed, since  $\nabla f_{\mS,x}(x^*) =0 $ and $ f_x(x^*) =0$,
we have that~\eqref{eq:1smooth} implies that
\begin{eqnarray*}
\norm{\nabla f_{\mS,x}(x) -\nabla f_{\mS,x}(x^*) }^2 \leq 2 (f_{\mS, x}(x) -f_{\mS,x}(x^*)),
\end{eqnarray*}
which is usually a consequence of assuming that $f_{\mS, x}(x) $ is convex and $1$--smooth (see Theorem 2.1.5 and Equation 2.1.7 in~\cite{NesterovBook}).
Yet in our case, eq.~\eqref{eq:1smooth} is a direct consequence of the definition of $f_{\mS,x}$ as opposed to being an extra assumption. This gratuitous property will be key in establishing a global convergence result.


\subsection{Convergence for star-convex}
We use the shorthand $f_{k}(x) \eqdef f_{x^k}(x)$, $f_{\mS_k,k} \eqdef f_{\mS_k,x^k}$ and $\EE{k}{ \cdot } \eqdef \E{ \cdot  \mid x^k}$. 
Here we establish the global convergence of \SNR by supposing that $f_k$ is \emph{star-convex} which is a large class of nonconvex functions that includes convexity as a special case~\cite{Nesterov2006cubic,LeeV16,zhou2018sgd,hinder2020nearoptimal}.

\begin{assumption}[Star-Convexity]
\label{ass:convex}
Let $x^*$ satisfy Asm.~\ref{ass:zero}, i.e.\ let $x^*$ be a solution to~\eqref{eq:main}. For every $x^k$ given by Alg.~\ref{algo:SkeNeR} with $k \in \N$, we have that
\begin{eqnarray}\label{eq:cvx}
f_k(x^*) &\geq& f_k(x^k) + \dotprod{\nabla f_k(x^k), x^* - x^k}.
\end{eqnarray}
\end{assumption}

We now state our main theorem.
\begin{theorem} \label{theo:convex}
Let $x^*$ satisfy Asm.~\ref{ass:convex}. If $0 < \gamma < 1$,
then 
\begin{align}
\E{\min_{t=0, \ldots, k-1} f_t(x^t)}
\leq \frac{1}{k}\sum_{t=0}^{k-1}\E{f_t(x^t)} 
\leq \frac{1}{k}\frac{\norm{x^0-x^*}^2}{2\gamma\left(1 - \gamma\right)}. \label{eq:convminEf}
\end{align}
Written in terms of $F$ and for $\gamma =1/2$ the above gives
\[\E{\min_{t=0, \ldots, k-1}\norm{F(x^t)}_{\E{\mH_{\mS}(x^t)}}^2} \quad \leq \quad \frac{4\norm{x^0-x^*}^2}{k}.\]
Besides, if the stochastic function $f_{\mS,x}(x)$ is star-convex along the iterates $x^k$ 
, i.e.\
\begin{eqnarray}\label{eq:cvx-S-star}
 f_{\mS_k,x^k}(x^*) &\geq& f_{\mS_k,x^k}(x^k) + \dotprod{\nabla f_{\mS_k,x^k}(x^k), x^* - x^k}
\end{eqnarray}
for all $\mS_k \sim\cD_{x^k}$, then the iterates $x^k$ of \SNR~\eqref{eq:update} are bounded with
\begin{eqnarray}\label{eq:bounded-det}
\norm{x^k - x^*} &\leq& \norm{x^0 - x^*}.
\end{eqnarray}
\end{theorem}
\begin{proof}
Let $t \in  \{0, \ldots,  k - 1 \}$ and $\delta_t \eqdef x^t-x^* $. 
We have that
\begin{eqnarray}
\EE{t}{\norm{\delta_{t+1}}^2} &\overset{\eqref{eq:SGD}}{ =}& \EE{t}{\norm{x^t - \gamma \nabla f_{\mS_t,t}(x^t) - x^*}^2} \nonumber\\
&=& \norm{\delta_{t}}^2 - 2\gamma \dotprod{\delta_{t}, \nabla f_t(x^t)} + \gamma^2\EE{t}{\norm{\nabla f_{\mS_t,t}(x^t)}^2} \nonumber\\
&\overset{\eqref{eq:cvx}}{\leq}& \norm{\delta_{t}}^2 - 2\gamma (f_t(x^t) - f_t(x^*))+ \gamma^2\EE{t}{\norm{\nabla f_{\mS_t,t}(x^t)}^2} \nonumber \\
&\overset{ \eqref{eq:1smooth}}{ = }& \norm{\delta_{t}}^2 - 2\gamma\left(1 - \gamma\right)(f_t(x^t) -f_t(x^*))  \nonumber\\
&\overset{f_t(x^*) = 0}{=}& \norm{\delta_{t}}^2 - 2\gamma\left(1 - \gamma\right)f_t(x^t). \label{EEt}
\end{eqnarray}
Taking total expectation  for all $t \in \{0, \ldots,  k - 1 \}$, we have that
\begin{eqnarray}
\E{\norm{\delta_{t+1}}^2} &\leq& \E{\norm{\delta_{t}}^2} - 2\gamma\left(1 - \gamma\right)\E{f_t(x^t)}. \label{eq:Etotal}
\end{eqnarray}
Summing both sides of~\eqref{eq:Etotal} from $0$ to $k-1$ gives
\begin{eqnarray*}
\E{\norm{x^k-x^*}^2} + 2\gamma\left(1 - \gamma\right)\sum_{t=0}^{k-1}\E{f_t(x^t)} &\leq& \norm{x^0-x^*}^2.
\end{eqnarray*}
Dividing through by $2\gamma\left(1 - \gamma\right) >0$ and by $k$, we have that
\[\E{\min_{t=0, \ldots, k-1}f_t(x^t)}  \leq  \min_{t=0, \ldots, k-1} \E{f_t(x^t)}  \leq  \frac{1}{k}\sum_{t=0}^{k-1}\E{f_t(x^t)}  \leq  \frac{1}{k}\frac{\norm{x^0-x^*}^2}{2\gamma\left(1 - \gamma\right)},\]
where in the most left inequality we used Jensen's inequality.

Finally, if~\eqref{eq:cvx-S-star} holds, then we can repeat the steps leading up to~\eqref{EEt} without the conditional expectation, so that
\begin{eqnarray*}
\norm{\delta_{t+1}}^2 &\overset{\eqref{eq:SGD}+\eqref{eq:cvx-S-star}+\eqref{eq:1smooth}}{\leq}&\norm{\delta_{t}}^2 - 2\gamma\left(1 - \gamma\right)f_{\mS_t,t}(x^t).
\end{eqnarray*}
Since $f_{S_t,t}(x^t) \geq 0$, we have $\norm{\delta_{t+1}}^2 \leq \norm{\delta_{t}}^2$, i.e.\ \eqref{eq:bounded-det} holds.
\end{proof}


Thm.~\ref{theo:convex} is an unusual result for \SGD methods. Currently, to get a $\cO(1/k)$ convergence rate for \SGD, one has to assume smoothness and strong convexity~\cite{Gower19} or convexity, smoothness and interpolation~\cite{VaswaniBS19}.
Here we get a $\cO(1/k)$ rate by \emph{only} assuming star-convexity. This is because we have smoothness and interpolation properties as a by-product due to our  reformulation~\eqref{eq:sgdreform}. 
\ruiline{However, the star-convexity assumption of $f_k(\cdot)$ for all $k \in \N$ is hard to interpret in terms of assumptions on $F$ in general. But, we are able to interpret it in many important extremes.  That is, for the full \NR method, we show that it suffices for the Newton direction to be $2$--co-coercive (see~\eqref{eq:nxcvx} in Sec.~\ref{sec:globalNR}). For the other extreme where the sketching matrix samples a single row, then the star-convexity assumption is even easier to check, and is guaranteed to hold so long as $F_i(x)^2$ is convex for all $i=1, \cdots, m$ (see Sec.~\ref{sec:Kaczmarz}). 
}

\ruiline{
Next, we will show 
the convergence of $F(x^k)$ instead of $f_k(x^k)$ via Thm.~\ref{theo:convex}.}




\subsubsection{Sublinear convergence of the Euclidean norm \texorpdfstring{$\norm{F}$}{|F|}}
\label{sec:EuclideanF}
If $\E{\mH_{\mS}(x)}$ is invertible for all $x \in \R^p$, we can use Thm.~\ref{theo:convex} with the bound~\eqref{eq:bounded-det} to guarantee that $\norm{F}$ converges sublinearly. 
Indeed, \ruiline{when $\E{\mH_{\mS}(x)}$ is invertible, $\E{\mH_{\mS}(x)}$ is symmetric positive definite.}
Thus there exists $\lambda>0$ that bounds the smallest eigenvalue away from zero in any closed bounded set 
\ruiline{(e.g.\ $\{x \in \R^p \mid \norm{x-x^*} \leq \norm{x^0 - x^*} \}$\footnote{\ruiline{We can re-write the set as the closure of the ball $\{x \in \R^p \mid x \in \overline{\cB(x^*, \norm{x^0 - x^*})}\}$. }})}:
\begin{equation} \label{eq:smalleigbnd} 
\min_{x \in \{x \mid \norm{x-x^*} \leq \norm{x^0 - x^*}\}} \lambda_{\min} \left(\E{\mH_{\mS}(x)}\right)  =  \lambda  >  0,
\end{equation}
where  $\lambda_{\min}(\cdot)$ is the smallest eigenvalue operator.
Consequently, under the assumption of Thm.~\ref{theo:convex}
with the condition~\eqref{eq:cvx-S-star}, from~\eqref{eq:bounded-det} and~\eqref{eq:smalleigbnd},
we have
\begin{equation}
\lambda  \E{\min_{t=0, \ldots, k-1}\norm{F(x^t)}^2}   \leq  \E{\min_{t=0, \ldots, k-1}\norm{F(x^t)}^2_{\E{\mH_{\mS}(x^t)}}}  \overset{\eqref{eq:convminEf}}{\leq}  \frac{1}{k}\frac{\norm{x^0-x^*}^2}{\gamma\left(1 - \gamma\right)}.  \label{eq:FEHsk}
\end{equation}

%
It turns out that using the smallest eigenvalue of $\E{\mH_{\mS}(x)}$ in the above bound  is overly pessimistic.
\ruiline{
To improve it, first note that we do \emph{not need} that $\E{\mH_{\mS}(x)}$ is invertible. Instead, we only need that
 $F(x) \in {\bf Im}(DF(x)^\top) \subset {\bf Im}(\E{\mH_{\mS}(x)})$, as we show in Cor.~\ref{lem:EuclidieanF}.} But first, we need the following lemma.
\begin{lemma}[Lemma 10 in~\cite{RSN_nips}] \label{lem:ker}
For any matrix $\mW$ and symmetric positive semi-definite matrix $\mG$ s.t.\
${\bf Ker}(\mG)  \subset  {\bf Ker}(\mW),$
we have
${\bf Ker}(\mW^\top)  =  {\bf Ker}(\mW\mG\mW^\top).$
\end{lemma}


\ruilinex{Note $L \eqdef \sup_{x \in \{x \mid \norm{x-x^*} \leq \norm{x^0 - x^*}\}} \|DF(x)\| > 0$.
Such $L$ exists because $x$ is in a closed bounded convex set and because we have assumed that $DF(\cdot)$ is continuous. A continuous mapping over a closed bounded convex set is bounded.
Now we can state the sublinear convergence results for $\norm{F}$.}
\begin{corollary} \label{lem:EuclidieanF}
Let
\begin{align}
 \rho(x)  &\eqdef  \min_{v \in {\bf Im}\left(DF(x)\right) / \{0\}}\frac{v^\top DF(x)\E{\mH_{\mS}(x)}DF(x)^\top v}{\norm{v}^2}, \label{eq:defrho} \\
 \rho  &\eqdef  \min_{x \in \{x \mid \norm{x-x^*} \leq \norm{x^0 - x^*}\}} \rho(x). \label{eq:defrho2}
\end{align}
It follows that $0\leq \rho(x)\leq 1$.
\ruilinex{If   
\begin{equation}\label{eq:DFTinEmHass}
 F(x) \; \in \; {\bf Im}(DF(x)^\top) \; \subset \; {\bf Im}(\E{\mH_{\mS}(x)}) \quad \mbox{for all } x \in \R^p,
\end{equation}
then
$\quad \quad \rho(x) = \lambda^+_{\min}\left(DF(x)\E{\mH_{\mS}(x)}DF(x)^\top\right) > 0 \quad \forall x \in \R^p, \quad \mbox{and} \quad \rho >0,$
where $\lambda^+_{\min}$ is the smallest \emph{non-zero} eigenvalue.  
Furthermore, if the star-convexity for each sketching matrix~\eqref{eq:cvx-S-star} holds
, then
\begin{eqnarray} \label{eq:EuclidieanF}
\E{\min_{t=0,\ldots,k-1}\norm{F(x^t)}^2} &\leq& \frac{1}{k}\cdot\frac{L^2\norm{x^0-x^*}^2}{\rho\gamma\left(1 - \gamma\right)}.
\end{eqnarray}}
\end{corollary}
\begin{proof}
First recall that
\begin{eqnarray*} 
DF(x)\mH_{\mS}(x)DF(x)^\top DF(x)\mH_{\mS}(x)DF(x)^\top = DF(x)\mH_{\mS}(x)DF(x)^\top \quad \mbox{for all } x \in \R^p,
\end{eqnarray*}
which is shown in the proof of Lemma~\ref{lem:smooth}.
Thus $DF(x)\mH_{\mS}(x)DF(x)^\top$ is a projection. By Jensen's inequality, the eigenvalues of an expected projection are between $0$ and $1$. Thus by the definition of $\rho(x)$, we have $0\leq\rho(x)\leq 1$.
\ruiline{Next,  
by~\eqref{eq:DFTinEmHass}, we have
${\bf Ker}\left(\E{\mH_{\mS}(x)}\right) \subset {\bf Ker}\left(DF(x)\right)$. Thus, by Lemma~\ref{lem:ker} we have that}
\begin{equation} \label{eq:imDFequalkerperp}
{\bf Im}\left(DF(x)\right)  =  \left({\bf Ker}\left(DF(x)^\top\right)\right)^\perp  =  \left({\bf Ker}\left(DF(x)\E{\mH_{\mS}(x)}DF(x)^\top\right)\right)^\perp,
\end{equation}
where the second equality is obtained by Lemma~\ref{lem:ker}. Now from the definition of $\rho(x)$ in~\eqref{eq:defrho}, we have
\begin{eqnarray}
\rho(x) &\overset{\eqref{eq:imDFequalkerperp}}{=}&\min_{v \in \left({\bf Ker}\left(DF(x)\E{\mH_{\mS}(x)}DF(x)^\top\right)\right)^\perp / \{0\}}\frac{v^\top DF(x)\E{\mH_{\mS}(x)}DF(x)^\top v}{\norm{v}^2} \nonumber \\
&=& \lambda^+_{\min}\left(DF(x)\E{\mH_{\mS}(x)}DF(x)^\top\right)  >  0. \nonumber
\end{eqnarray}
It now follows that $\rho >0$, since the definition of $\rho$ in~\eqref{eq:defrho2} is given by minimizing $\rho(x)$ over the closed bounded set $\{x \mid \norm{x-x^*} \leq \norm{x^0 - x^*}\}.$
Next, \ruiline{given $x \in \{x \mid \norm{x-x^*} \leq \norm{x^0 - x^*}\}$, since $F(x) \in {\bf Im}(DF(x)^\top)$ by $\eqref{eq:DFTinEmHass}$ and notice that ${\bf Im}(DF(x)^\top) = {\bf Im}(DF(x)^\top DF(x))$, there exists $v \in \R^m$ s.t.\ $F(x) = DF(x)^\top DF(x)v$. 

\xline{If $F(x) \neq 0$, then $DF(x)v \in {\bf Im}\left(DF(x)\right) / \{0\}$, we have}
\begin{eqnarray}
\norm{F(x)}^2_{\E{\mH_{\mS}(x)}} &=& v^\top DF(x)^\top DF(x)\E{\mH_{\mS}(x)}DF(x)^\top DF(x)v \nonumber \\
&\overset{\eqref{eq:defrho}}{\geq} & \rho(x)v^\top DF(x)^\top DF(x)v.\label{eq:zsnlo8jh8o} 
\end{eqnarray}
Since $F(x) = DF(x)^\top DF(x)v$ and ${\bf Im}(DF(x)^\top) \oplus {\bf Ker}(DF(x)) = \R^m$,\footnote{\ruiline{The operator $\oplus$ denotes the direct sum of two vector spaces.}} we have that
\[\exists! \ y \in {\bf Ker}(DF(x)) \subset \R^m \mbox{ s.t.\ } v =  (DF(x)^\top DF(x))^{\dagger} F(x) + y.\]
Thus 
\[DF(x) v = DF(x) (DF(x)^\top DF(x))^{\dagger} F(x) = (DF(x)^\top)^{\dagger} F(x).\]
Substituting this in~\eqref{eq:zsnlo8jh8o},} we have that
\begin{eqnarray}
\norm{F(x)}^2_{\E{\mH_{\mS}(x)}} \geq \rho(x)\norm{F(x)}^2_{\left(DF(x)^\top DF(x)\right)^{\dagger}} 
\geq \frac{\rho}{L^2}\norm{F(x)}^2, \label{eq:FEHs}
\end{eqnarray}
where 
on the last inequality, we use that $\sup_{x \in \{x \mid \norm{x-x^*} \leq \norm{x^0 - x^*}\}} \|DF(x)\| \leq L$ and $\rho(x) \geq \rho$ by the definition of $\rho$ in~\eqref{eq:defrho2} 
.

\xline{If $F(x) = 0$,~\eqref{eq:FEHs} still holds. Thus, for all $x \in \{x \mid \norm{x-x^*} \leq \norm{x^0 - x^*}\}$,~\eqref{eq:FEHs} holds.}
Consequently by Thm.~\ref{theo:convex} and~\eqref{eq:bounded-det} under the star-convexity condition~\eqref{eq:cvx-S-star} with $\norm{x^t - x^*} \leq \norm{x^0 - x^*}$ for all $t \in \{0, \cdots, k-1\}$, we have that
\[ \frac{\rho}{L^2} \E{\min_{t=0,\ldots,k-1}\norm{F(x^t)}^2}  \overset{\eqref{eq:FEHs}}{\leq}  \E{\min_{t=0,\ldots,k-1} \norm{F(x^t)}^2_{\E{\mH_{\mS}(x^t)}}}  \overset{\eqref{eq:convminEf}}{\leq} \frac{1}{k} \frac{\norm{x^0-x^*}^2}{\gamma\left(1 - \gamma\right)},\]
which after multiplying through by $\left. L^2 \right/\rho >0$ concludes the proof.
\end{proof}

\ruiline{
Thus with Cor.~\ref{lem:EuclidieanF}, we show that $F(x^t)$ converges to zero. This lemma relies on the inclusion~\eqref{eq:DFTinEmHass}, which in turn imposes some restrictions on the sketching matrix and $F(x)$. In our forthcoming examples in Sec.~\ref{sec:globalNR} and~\ref{sec:Kaczmarz}, we can directly verify the inclusion of~\eqref{eq:DFTinEmHass}. For other examples in Sec.~\ref{sec:SNM} and~\ref{sec:GLM}, we provide the following Lemma~\ref{lem:reformulationnew} where we give sufficient conditions for~\eqref{eq:DFTinEmHass} to hold.
\begin{lemma}\label{lem:reformulationnew}
Let $F(x) \in {\bf Im}(DF(x)^\top)$. Furthermore, we suppose that  $\mS \sim \cD_x$ is \emph{adapted} to $DF(x)$ 
by which we mean
\begin{align}
{\bf Ker} \left(\E{\mS\mS^\top}\right)  \; \subset \; {\bf Ker}(DF(x)) &  \subset \; {\bf Ker}\left(\mS^\top\right) , \quad \mbox{for all } \mS \sim \cD_x. \label{eq:adaptedsketch}
\end{align}
Then it follows that~\eqref{eq:DFTinEmHass} holds for all $x \in \R^p$.
\end{lemma}
\begin{proof}
Since  ${\bf Ker}(DF(x)^\top DF(x)) = {\bf Ker}(DF(x)) \overset{\eqref{eq:adaptedsketch}}{\subset} {\bf Ker}(\mS^\top)$,
we have
\begin{align}
{\bf Ker}\left(\left(\mS^\top DF(x)^\top DF(x)\mS\right)^\dagger\right)  &=  {\bf Ker}\left(\mS^\top DF(x)^\top DF(x)\mS\right)  =  {\bf Ker}(\mS), \label{eq:ker2new}
\end{align}
where 
the last equality is obtained by Lemma~\ref{lem:ker} with ${\bf Ker}(DF(x)^\top DF(x)) \subset {\bf Ker}(\mS^\top)$.
Thus, using Lemma~\ref{lem:ker} again with $\mG = \left(\mS^\top DF(x)^\top DF(x)\mS\right)^\dagger$, $\mW = \mS$ and ${\bf Ker}(\mG) \subset {\bf Ker}(\mW)$ given by~\eqref{eq:ker2new}, we have that
\begin{equation}\label{eq:ker3new}
{\bf Ker}(\mH_{\mS}(x)) \overset{\eqref{eq:Hk}}{=} {\bf Ker}(\mS\left(\mS^\top DF(x)^\top DF(x)\mS\right)^\dagger\mS^\top)  =  {\bf Ker}(\mS^\top)  =  {\bf Ker}(\mS\mS^\top).
\end{equation}
\noindent \ruiline{As $\mH_{\mS}(x)$ is symmetric positive semi-definite $\forall \ \mS\sim\cD_x$, we have that
\begin{align*}
&v \in {\bf Ker}\left(\E{\mH_{\mS}(x)}\right)
 \Longleftrightarrow  \E{\mH_{\mS}(x)} v = 0
 \Longleftrightarrow  \norm{v}_{\E{\mH_{\mS}(x)}}^2 = 0 \ \ \left(\mbox{as } \E{\mH_{\mS}(x)} \succeq 0\right) \\ 
 &\Longleftrightarrow  \E{\norm{v}_{\mH_{\mS}(x)}^2} = 0 \Longleftrightarrow \int_{\mS} \norm{v}_{\mH_{\mS}(x)}^2 d\mathbb{P}_{\cD_x}(\mS) = 0 \\
 &\Longleftrightarrow  \norm{v}_{\mH_{\mS}(x)}^2 = 0 \; \forall \; \mS \sim \cD_x
 \ \ \left(\mbox{as } \norm{v}_{\mH_{\mS}(x)}^2  \geq 0  \; \forall \; \mS\right) \\
 &\overset{\mH_{\mS}(x) \succeq 0}{\Longleftrightarrow} v \in {\bf Ker}\left(\mH_{\mS}(x)\right) \; \forall \; \mS \sim \cD_x 
 \Longleftrightarrow v \in \bigcap_{\mS\sim\cD_x} {\bf Ker}\left(\mH_{\mS}(x)\right),
\end{align*}
where we use $\bigcap_{\mS\sim\cD_x} {\bf Ker}\left(\mH_{\mS}(x)\right)$ to note the intersection of the random subsets ${\bf Ker}\left(\mH_{\mS}(x)\right)$ for all $\mS \sim \cD_x$.
}
Similarly, we have ${\bf Ker}\left(\E{\mS\mS^\top}\right) = \bigcap_{\mS\sim\cD_x} {\bf Ker}\left(\mS\mS^\top\right)$ because $\mS\mS^\top$ is also symmetric, positive semi-definite for all $\mS \sim \cD_x$.
Thus we have
\begin{eqnarray*}
{\bf Ker}(\E{\mH_{\mS}(x)}) &=& \bigcap_{\mS\sim\cD_x} {\bf Ker}(\mH_{\mS}(x)) \\ 
&\overset{\eqref{eq:ker3new}}{=}& \bigcap_{\mS\sim\cD_x}{\bf Ker}(\mS\mS^\top) 
\; = \; {\bf Ker}(\E{\mS\mS^\top}) 
\; \overset{\eqref{eq:adaptedsketch}}{\subset} \; {\bf Ker}(DF(x)).
\end{eqnarray*}
Consequently, by considering the complement of the above, we arrive at~\eqref{eq:DFTinEmHass}.
\end{proof}}

\ruilinex{We refer to a sketching matrix $\mS \sim \cD_x$ that satisfies~\eqref{eq:adaptedsketch} as a sketch that is \emph{adapted to} $DF(x)$. One easy way to design such adapted sketches is 
the following.
\begin{lemma} \label{lem:adaptedsketch}
Let $\hat{\mS} \in \R^{p \times \tau} $ s.t.\ $\hat{\mS} \sim \cD$ a fixed distribution independent to $x$ and ${\bf Ker} (\mathbb{E}[\hat{\mS} \hat{\mS}^\top] ) \subset {\bf Ker}(DF(x)^\top)$. Thus, $\mS = DF(x)^\top \hat{\mS} \in \R^{m \times \tau}$ is adapted to $DF(x)$.
\end{lemma}
\begin{proof}
First, $ {\bf Ker}(DF(x)) \subset {\bf Ker}(\hat{\mS}^\top DF(x)) =  {\bf Ker}(\mS^\top).$ Furthermore, from Lemma~\ref{lem:ker} with  ${\bf Ker} (\mathbb{E}[\hat{\mS} \hat{\mS}^\top] ) \subset {\bf Ker}(DF(x)^\top)$, we conclude the proof with
\[ {\bf Ker} \left(\E{\mS\mS^\top}\right) \; =\;  {\bf Ker} (DF(x)^\top\mathbb{E}[\hat{\mS} \hat{\mS}^\top]DF(x)) \; \subset   {\bf Ker}(DF(x)). \]
\end{proof}}

\ruilinex{
The condition ${\bf Ker} (\mathbb{E}[\hat{\mS} \hat{\mS}^\top] ) \subset {\bf Ker}(DF(x)^\top)$ in Lemma~\ref{lem:adaptedsketch}
holds for many standard sketches including Gaussian and subsampling sketches presented as follows.
\begin{lemma}\label{lem:SSTinv}
For Gaussian and uniform subsampling sketches defined in Sec.~\ref{sec:sketch}, we have that $\E{\mS\mS^\top} = c \mI_m $ with $c>0$ a fixed constant depending on the sketch.
\end{lemma}
\begin{proof}
For Gaussian sketches with $\mS_{ij} \sim \cN(0,\frac{1}{\tau})$, we have that $c=1$. Indeed, since the mean is zero, off-diagonal elements of  $\E{\mS\mS^\top}$ are all zero. We note $\mS_{i:}$ the $i$th row of $\mS$, then the $i$th diagonal element of the matrix $\E{\mS\mS^\top}$  is given by
\[
\E{\mS_{i:}\mS_{i:}^\top} =\sum_{j=1}^{\tau} \E{\mS_{ij}^2} = \sum_{j=1}^{\tau} \frac{1}{\tau} =1.
\]
For the uniform subsampling sketch~\eqref{eq:subsample}, we have again that off-diagonal elements are zero since the rows of $\mS$ are orthogonal. The diagonal elements are constant with
 \begin{align*}
\E{\mS_{i:}\mS_{i:}^\top} 
=\frac{1}{\binom{m}{\tau}}\sum_{\substack{C \subset \{1,\ldots, m\}, |C|=\tau, i\in C}} 1  = \frac{\binom{m-1}{\tau-1}}{\binom{m}{\tau}} = \frac{\tau}{m}, \quad \mbox{for all } i = 1, \cdots, m.
\end{align*}
\end{proof}}

\ruilinex{
From Lemma~\ref{lem:SSTinv}, we know that $\mathbb{E}[\hat{\mS} \hat{\mS}^\top] = c\mI_p$ invertible with $c>0$. Thus ${\bf Ker} (\mathbb{E}[\hat{\mS} \hat{\mS}^\top] ) = \{0\} \subset {\bf Ker}(DF(x)^\top)$ holds for \emph{any} sketch size $\tau$.}

\subsection{Convergence for strongly convex}

Here we establish a global linear convergence of \SNR when assuming that $f_y$ is strongly quasi-convex.

\begin{assumption}[$\mu$-Strongly Quasi-Convexity]\label{ass:strconv}
Let $x^*$ satisfy Asm.~\ref{ass:zero}
and
\begin{equation}\label{eq:strconv}
\exists \ \mu > 0 \mbox{ s.t.\ } \ f_y(x^*) \geq f_y(x) + \dotprod{\nabla f_y(x), x^*-x} + \frac{\mu}{2} \norm{x^*-x}^2 \ \ \forall \ x, y \in \R^p.
\end{equation}
\end{assumption}

This Asm.~\ref{ass:strconv} is strong, so much so, we have the following lemma.

\begin{lemma}\label{lem:assstrconvimplyasskerr}
Asm.~\ref{ass:strconv} implies~\eqref{ass:ker} and that the solution to~\eqref{eq:main} is unique.
\end{lemma}

\begin{proof}
Let $y \in \R^p$ and let $u \in F(\R^p) \cap {\bf Ker}\left(\E{\mH_{\mS}(y)}\right)$. $u \in F(\R^p)$ implies that $\exists x \in \R^p$ s.t.\ $F(x) = u$. Besides, $u \in {\bf Ker}\left(\E{\mH_{\mS}(y)}\right)$ implies that $\E{\mH_{\mS}(y)}F(x) = 0$. Now we apply \eqref{eq:strconv} at point $x$ knowing that $f_y(x^*) = 0$:
\begin{align*}
&\quad \quad 0 \geq f_y(x) + \dotprod{\nabla f_y(x), x^*-x} + \frac{\mu}{2} \norm{x^*-x}^2 \\
&\Longrightarrow 0 \geq 0 + \dotprod{0, x^*-x} + \frac{\mu}{2} \norm{x^*-x}^2 \ \ \left(\mbox{as } \E{\mH_{\mS}(y)}F(x) = 0\right)
\Longleftrightarrow x = x^*.
\end{align*}
Thus $F(x) = u = 0$. We conclude $F(\R^p) \cap {\bf Ker}\left(\E{\mH_{\mS}(y)}\right) = \{0\}$, i.e.~\eqref{ass:ker} holds.

Besides, let $x'$ be a global minimizer of $f_y(\cdot)$. Then $f_y(x') = f_y(x^*) = 0$ and $\nabla f_y(x') = 0$. Similarly, by applying \eqref{eq:strconv} at point $x'$, we obtain $x' = x^*$. Consequently, $x^*$ is the unique minimizer of $f_y(\cdot)$ for all $y$, thus the unique solution to \eqref{eq:main}, according to~\eqref{ass:ker} and Lemma~\ref{lem:equiv}.
\end{proof}

Under 
Asm.~\ref{ass:strconv}
, choosing $\gamma = 1$ guarantees a fast global linear convergence.

\begin{theorem}\label{theo:strconv}
If $x^*$ satisfies Asm.~\ref{ass:strconv} and $\gamma \leq 1$, then \SNR converges linearly:
\begin{eqnarray}
\E{\norm{x^{k+1}-x^*}^2} &\leq& (1-\gamma \mu)^{k+1} \norm{x^{0}-x^*}^2 \quad \quad \mbox{with } \ \mu \leq 1.
\end{eqnarray}
\end{theorem}

\begin{proof}
Let $\delta_k \eqdef x^{k} -x^*.$
By expanding the squares, similarly we have that
\begin{eqnarray}
\EE{k}{\norm{\delta_{k+1}}^2} &=& \norm{\delta_k }^2 - 2\gamma \dotprod{\delta_k , \nabla f_k(x^k)} + \gamma^2\EE{k}{\norm{\nabla f_{\mS_k,k}(x^k)}^2} \nonumber \\
&\overset{\eqref{eq:strconv}}{\leq}& (1-\gamma \mu) \norm{\delta_k }^2 - 2\gamma (f_k(x^k) - f_k(x^*))+ \gamma^2\EE{k}{\norm{\nabla f_{\mS_k,k}(x^k)}^2} \nonumber \\
&\overset{\eqref{eq:1smooth}}{\leq }& (1-\gamma \mu) \norm{\delta_k }^2 - 2\gamma\left(1 - \gamma\right)(f_k(x^k) -f_k(x^*)) \nonumber \\
&\leq& (1-\gamma \mu) \norm{\delta_k }^2. \nonumber \quad \quad \left(\mbox{since } \gamma\left(1 - \gamma\right)(f_k(x^k) -f_k(x^*)) \geq 0\right)
\end{eqnarray}
Now by taking total expectation, we have that
\[\E{\norm{x^{k+1}-x^*}^2} \ \leq \ (1-\gamma \mu) \E{\norm{x^k-x^*}^2} \ \leq \ (1-\gamma \mu)^{k+1}\norm{x^0-x^*}^2.\]

Next, we show that  $\mu \leq 1$. In fact, when we imply \eqref{eq:strconv} at the point $x^k$, it shows
\begin{eqnarray}\label{eq:mu}
\eqref{eq:strconv} &\overset{\eqref{eq:1smooth} }{\Longrightarrow}& f_k(x^*) \geq \frac{1}{2}\EE{k}{\norm{\nabla f_{\mS_k,k}(x^k)}^2} + \dotprod{x^*-x^k, \nabla f_{k}(x^k)} + \frac{\mu}{2}\|x^*-x^k\|^2 \nonumber \\
&\Longleftrightarrow& f_k(x^*) \geq \frac{1}{2}\EE{k}{\norm{x^* - \left(x^k - \nabla f_{\mS_k,k}(x^k)\right)}^2} - \frac{1-\mu}{2}\|x^*-x^k\|^2 \nonumber \\
&\overset{f_k(x^*) = 0}{\Longrightarrow}& (1-\mu)\|x^*-x^k\|^2 \geq \EE{k}{\norm{x^* - \left(x^k - \nabla f_{\mS_k,k}(x^k)\right)}^2} \geq 0. \nonumber
\end{eqnarray}
Thus $\mu \leq 1$.
\end{proof}



\section{New global convergence theory of the \NR method}
\label{sec:globalNR}

As a direct consequence of our general convergence theorems, in this section we develop a new global convergence theory for the original \NR method. We first provide the results in $1$-dimension in Sec.~\ref{sec:globalNR-1D}, then a general result in higher dimensions in the subsequent Sec.~\ref{sec:globalNR-dD} and compare this result to the classic monotone convergence theory in Sec.~\ref{sec:monotoneconv}.

\subsection{A single nonlinear equation} \label{sec:globalNR-1D}

Consider the case where $F(x) = \phi(x) \in \R$ is a one dimensional function and $x \in \R$. This includes common applications of the \NR method such as calculating square roots of their reciprocal\footnote{Used in particular to compute angles of incidence and reflection in games such as quake~\url{https://en.wikipedia.org/wiki/Fast_inverse_square_root}} and finding roots of polynomials.
Even in this simple one dimension case, we find that \ruilinex{our assumptions of global convergence given in Cor.~\ref{lem:EuclidieanF} are \emph{strictly weaker} than the standard assumptions used to guarantee \NR convergence,} as we explain next.

The \NR method in one dimension at every iteration $k$ is given by
\[x^{k+1} \; = \; x^k - \frac{\phi(x^k)}{\phi'(x^k)} \; \eqdef \; g(x^k).\]

To guarantee that this is well defined, we assume that $\phi'(x^k) \neq 0$ for all $k$.
A sufficient condition for this procedure to converge locally is that $|g'(x)| < 1$ with $x \in I$ where $I$ is a given interval containing the solution $x^*.$
See for example Section 1.1 in~\cite{DeuflhardNewton:2011} or Chapter 12 in~\cite{Ortega:2000}. We can extend this to a global convergence by requiring that $|g'(x)|<1$ globally. In the case of \NR, since
$g'(x) \; = \; 1- \frac{\phi'(x)^2 - \phi(x)\phi''(x)}{\phi'(x)^2} \; = \; \frac{ \phi(x)\phi''(x)}{\phi'(x)^2},$
this condition amounts to requiring
\begin{eqnarray}\label{eq:1dconvergence}
\frac{ |\phi(x)\phi''(x)|}{\phi'(x)^2} &<& 1.
\end{eqnarray}
Curiously, this condition~\eqref{eq:1dconvergence} has an interesting connection to convexity. In fact, condition~\eqref{eq:1dconvergence} implies that $\phi^2(x)$ is convex and twice continuously differentiable. To see this, note that
$\frac{d^2}{dx^2}\phi^2(x) \geq 0$ is equivalent to
\begin{equation}\label{eq:twicedifocn}
\frac{d^2}{dx^2}\phi^2(x) \; = \; 2\frac{d}{dx} \phi'(x)\phi (x) \; = \; 2\left(\phi(x)\phi''(x) + \phi'(x)^2\right) \; \geq \; 0.
\end{equation} 
Now it is easy to see that~\eqref{eq:1dconvergence} implies~\eqref{eq:twicedifocn}. Finally~\eqref{eq:twicedifocn} also implies that $\phi^2(x)$ is \emph{star-convex}, which is exactly what is required by our convergence theory in \ruiline{Cor.~\ref{lem:EuclidieanF}}.

Indeed, in this one dimensional setting, \ruilinex{Asm.~\ref{ass:convex} is equivalent to~\eqref{eq:cvx-S-star}} and our reformulation in~\eqref{eq:sgdreform} boils down to minimizing
$f_y(x) = \left(\phi(x) / \phi'(y)\right)^2$.
Thus by \ruiline{Cor.~\ref{lem:EuclidieanF}}, the \NR method converges globally if $f_{x^k}(x)$, or simply if $\phi(x)^2$ is star-convex and $\phi'(x^k)\neq 0$ for all iterates of $\NR,$
\ruilinex{which shows that our condition is strictly weaker than the other conditions, because there exists functions that are star-convex but not convex, e.g. $\phi(x)^2 = |x|(1-\exp(-|x|))$ from~\cite{Nesterov2006cubic,LeeV16}.} 

For future reference and convenience, we can re-write the star-convexity of each $\phi(x)^2$ as
\begin{eqnarray*}
0 &=& \phi(x^*)^2\; \geq \; \phi^2(x) + 2\phi(x) \phi'(x)(x^* - x),
\end{eqnarray*}
where $x^*$ is the global minimum of $\phi(x)^2$, i.e.\ $\phi(x^*) = 0$. This can be re-written as
\begin{eqnarray}\label{eq:unicond}
0 &\geq& \phi(x) \left(\phi(x) + 2\phi'(x)(x^* - x)\right).
\end{eqnarray}
By verifying~\eqref{eq:unicond} and that $\phi'(x^k) \neq 0$ on the iterates of \NR, we can guarantee that the method converges globally.

\subsection{The full \NR} \label{sec:globalNR-dD}

Now let $F(x) \in \R^m$ and consider the full \NR method~\eqref{eq:newton}. \ruilinex{Similarly, since $\mS = \mI_m$, Asm.~\ref{ass:convex} is equivalent to~\eqref{eq:cvx-S-star}. Cor.~\ref{lem:EuclidieanF} sheds some new light on the convergence of \NR.} In this case, our reformulation~\eqref{eq:sgdreform} is given by
\begin{equation} \label{eq:fkxfullnewt}
f_y(x) \; = \; \tfrac{1}{2}F(x)^\top (DF(y)^\top DF(y))^\dagger F(x) \; = \; \tfrac{1}{2}\norm{\left(DF(y)^\top\right)^\dagger F(x)}^2
\end{equation}
and \ruiline{Cor.~\ref{lem:EuclidieanF}} states that \NR converges if $f_{x^k}(x)$ is star-convex for all the iterates $x^k \in \R^p$. This has a curious re-interpretation in this setting. Indeed, let
\begin{eqnarray} \label{eq:nx}
n(x) &\eqdef& - (DF(x)^\top)^\dagger F(x)
\end{eqnarray}
be the Newton direction. From~\eqref{eq:fkxfullnewt} and~\eqref{eq:nx}, we have that
\begin{eqnarray} \label{eq:fxxnx2}
f_x(x) &=& \frac{1}{2} \norm{n(x)}^2.
\end{eqnarray}
Using~\eqref{eq:fxxnx2}, \ruiline{Cor.~\ref{lem:EuclidieanF}} can be stated in this special case as the following corollary.

\begin{corollary}\label{cor:FullNewton}
Consider $x^k$ given by the \NR~\eqref{eq:newton} with $\gamma <1$. If we have
\begin{eqnarray}
F(x) & \in & {\bf Im}(DF(x)^\top), 
\label{eq:KermSI} \\
\frac{1}{2}\norm{n(x)}^2 &\leq& \dotprod{n(x),x^*-x} \label{eq:nxcvx}
\end{eqnarray} 
hold for every $x = x^k$ with solution $x^*$, then 
\ruilinex{it exists $L > 0$ s.t.\ $ \|DF(x^k)\| \leq L$ and
\begin{eqnarray} \label{eq:NREuclidieanFnew}
\min_{t=0,\ldots,k-1}\norm{F(x^t)}^2 &\leq& \frac{1}{k}\cdot\frac{L^2\norm{x^0-x^*}^2}{\gamma\left(1 - \gamma\right)}.
\end{eqnarray}}
\end{corollary}
\begin{proof}
From~\eqref{eq:SGDgrad}, we have that
\begin{align} \label{eq:nablafxxnx}
\nabla f_{x}(x) \; = \; DF(x)(DF(x)^\top DF(x))^{\dagger} F(x) \; = \; (DF(x)^\top)^\dagger F(x) \; = \; - n(x).
\end{align}
Substituting~\eqref{eq:fxxnx2} and~\eqref{eq:nablafxxnx} in~\eqref{eq:cvx-S-star} yields~\eqref{eq:nxcvx}. 
Next, for $\mS =\mI_m$, we have that
\begin{eqnarray*}
{\bf Im} (\E{\mH_{\mS}(x)})) = {\bf Im} ((DF(x)^\top DF(x))^\dagger) = {\bf Im} (DF(x)^\top DF(x)) = {\bf Im}(DF(x)^\top).
\end{eqnarray*}
\ruilinex{Thus,
we have that $F(x) \ \in \; {\bf Im}(DF(x)^\top)  \; \subset \; {\bf Im}(\E{\mH_{\mS}(x)}))$, i.e.~\eqref{eq:DFTinEmHass} holds.
So all the conditions in Cor.~\ref{lem:EuclidieanF} are verified. Since $\mS = \mI_m$, we have that $\rho(x) = 1$ for all $ x$, so $\rho = 1$. 
Furthermore, because we assume that $DF(\cdot)$ is continuous and the iterates $x^k$ are in a closed bounded convex set~\eqref{eq:bounded-det} which is implied by~\eqref{eq:cvx-S-star} from Thm.~\ref{theo:convex},
there exists $L > 0$ s.t.\ $ \|DF(x^k)\| \leq L$ for all the iterates. Finally, by Cor.~\ref{lem:EuclidieanF}, the iterates converge sublinearly according to~\eqref{eq:EuclidieanF} which in this case is given by~\eqref{eq:NREuclidieanFnew}.}
\end{proof}

\ruiline{Condition~\eqref{eq:nxcvx} can be seen as a co-coercivity property of the Newton direction. 
This co-coercivity establishes a curious link with the modern proofs of convergence of gradient descent which rely on the co-coercivity of the gradient direction. 
That is, if $f(x)$ is convex and $L$--smooth, then we have that the gradient is $L$--co-coercive with
$$\frac{1}{L}\norm{\nabla f(x)}^2 \;\leq\; \dotprod{\nabla f(x),x-x^*}.$$
This is the key property for proving convergence of gradient descent, see e.g.\ Section 5.2.4 in~\cite{Bachbook}.
To the best of our knowledge, this is the first time that the co-coercivity of the Newton direction has been identified as a key property for proving convergence of the Newton's method.}
 In particular, global convergence results for the \NR method such as the monotone convergence theories (MCT) only hold for functions $F: \R^p \rightarrow \R^m$ with $p=m$ and rely on a stepsize $\gamma =1$, see~\cite{Ortega:2000,DeuflhardNewton:2011}. 
Cor.~\ref{cor:FullNewton} accommodates ``non-square'' functions $F: \R^p \rightarrow \R^m$.
Excluding the difference in stepsizes and focusing on ``square'' functions $F: \R^p \rightarrow \R^m$ with $p=m$, next we show in Thm.~\ref{theo:strongerthanmono}  that our assumptions are strictly weaker than those used for establishing the global convergence of \NR \ruiline{with constant stepsizes} through the MCT.

\subsection{Comparing to the classic monotone convergence theory of \NR}
\label{sec:monotoneconv}

Consider $m = p$.
Here we show that our \ruilinex{Asm.~\ref{ass:zero},~\eqref{eq:KermSI} and~\eqref{eq:nxcvx}} 
are strictly weaker than the classic assumptions used for establishing the global convergence of \NR~\ruiline{with constant stepsize}. To show this, we take the assumptions used in the MCT in Section 13.3.4 in~\cite{Ortega:2000} and compare with our assumptions in the following theorem.

\begin{theorem}\label{theo:strongerthanmono}
Let  $F:\R^p \rightarrow \R^p$ and 
let $x^k$ be the iterates of the \NR method with stepsize $\gamma =1$, that is
\begin{eqnarray}
x^{k+1} &=& x^k - \left(DF(x^k)^\top\right)^\dagger F(x^k). \label{eq:newton1}
\end{eqnarray}
Consider the following two sets of assumptions
\begin{enumerate}
\item[(I)] $F(x)$ is component wise convex, $(DF(x)^\top)^{-1}$ exists and is element-wise positive $\forall x \in \R^p$. There exist $x$ and $y$ s.t.\ $F(x) \leq 0 \leq F(y)$ element-wise.
\item[(II)] 
\ruilinex{There exists a unique $x^* \in \R^p$ s.t.\ $F(x^*) = 0$,~\eqref{eq:KermSI} and~\eqref{eq:nxcvx} hold for $k \geq 1$.}
\end{enumerate}
If (I) holds, then (II) always holds. Furthermore, there exist problems for which (II) holds and (I) does not hold.
\end{theorem}

\begin{proof}
First, we prove (I) $\implies$ (II). Assume that (I) holds. \ruiline{Since $DF(x)$ is invertible,~\eqref{eq:KermSI} holds trivially.} By 13.3.4 in~\cite{Ortega:2000}, we know that there exists a unique $x^* \in \R^p$  s.t.\ $F(x^*) =0$. It remains to verify if~\eqref{eq:nxcvx} holds for $k \geq 1$. First, note that the invertibility of $DF(x^k)$ gives
\begin{equation} \label{eq:fkxk}
f_{k}(x^k) \; = \; \frac{1}{2}\norm{F(x^k)}^2_{\left(DF_k^\top DF_k\right)^\dagger} \; = \; \frac{1}{2}\norm{(DF_k^{\top})^{-1}F_k}^2 \; \overset{\eqref{eq:newton1}}{=} \; \frac{1}{2}\norm{x^{k+1}-x^k}^2,
\end{equation}
with abbreviations $f_k(x^k) \equiv f_{x^k}(x^k)$, $F_k \equiv F(x^k)$ and $DF_k \equiv DF(x^k)$. Furthermore,
\begin{equation} \label{eq:nablafkxk}
\nabla f_k(x^k) \; = \; DF_k (DF_k^\top DF_k)^{-1}  F(x^k) \; = \; (DF_k^{\top})^{-1} F(x^k) \; \overset{\eqref{eq:newton1}}{=} \; x^k-x^{k+1}.
\end{equation}
Thus we can re-write the right hand side of the star-convexity assumption \eqref{eq:cvx}
as
\begin{align}
& f_k(x^k) + \dotprod{\nabla f_k(x^k), x^* - x^k} 
\overset{\eqref{eq:fkxk}+\eqref{eq:nablafkxk}}{=} \frac{1}{2}\norm{x^{k+1}-x^k}^2+ \dotprod{x^k-x^{k+1}, x^* - x^k} \nonumber\\
& \quad \quad = \frac{1}{2}\norm{x^{k+1}-x^k}^2 + \dotprod{x^k-x^{k+1}, x^{k+1} - x^k + x^* - x^{k+1}} \nonumber\\
& \quad \quad = -\frac{1}{2}\norm{x^{k+1}-x^k}^2+ \dotprod{x^k-x^{k+1}, x^* - x^{k+1}}. \nonumber 
\end{align}
From (I), we induce by Lemma 3.1 in~\cite{DeuflhardNewton:2011} that \NR is component wise monotone with $x^* \leq x^{k+1} \leq x^k$ for $k \geq 1$. Thus $x^k-x^{k+1} \geq 0$ and $x^* -x^{k+1} \leq 0$ component wise and consequently, $\dotprod{x^k-x^{k+1}, x^* - x^{k+1}} \leq 0.$ Thus it follows that
\begin{align*}
f_k(x^k) + \dotprod{\nabla f_k(x^k), x^* - x^k} 
\; &\leq \; 0 \; = \; f_k(x^*).
\end{align*}
Thus~\eqref{eq:nxcvx} holds for $k \geq 1$ and this concludes that (I) $\implies$ (II).

We now prove that (II) does {\bf not} imply (I). Consider the example $F(x) = Ax-b$, where $A \in \R^{p \times p}$ is invertible and $b \in \R^p.$ Thus, $DF(x) = A^\top$ is invertible and~\eqref{eq:KermSI} holds. As for~\eqref{eq:nxcvx}, let $x^*$ be the solution, i.e.\ $Ax^* =b,$ we have that
\begin{align*}
f_k(x) &= \frac{1}{2}\norm{F(x)}^2_{(DF(x_k)^\top DF(x_k))^{-1}} = \frac{1}{2}\norm{A(x-x^*)}^2_{(A A^\top)^{-1}} 
= \frac{1}{2}\norm{x-x^*}^2,
\end{align*}
which is a convex function and so~\eqref{eq:nxcvx} holds and thus (II) holds. However, (I) does not necessarily hold. Indeed, if $A = - \mI_p$, then $DF(x)$ is not element-wise positive.
\end{proof}

We observe that our assumptions are also strictly weaker than the affine covariates formulations of convex functions given in Lemma 3.1 in~\cite{DeuflhardNewton:2011}. The proof is verbatim to the above.

Thm.~\ref{theo:strongerthanmono} only considers the case that the stepsize $\gamma = 1$.
We also investigate the case where the stepsize $\gamma < 1$ in particular in $1$-dimension and show that MCT does not hold in this case in App.~\ref{sec:monotoneconv1}.
Thus we claim that our assumptions are strictly weaker than the assumptions used in MCT~\cite{Ortega:2000,DeuflhardNewton:2011} for establishing the global convergence of \NR, \ruiline{albeit for different step sizes.}

\section{Single row sampling: the nonlinear Kaczmarz method} \label{sec:Kaczmarz}

The \SNR enjoys many interesting instantiations. Among which, we have chosen three to present in the main text: the nonlinear Kaczmarz method in this section, the Stochastic Newton method~\cite{rodomanov16superlinearly,Kovalev2019stochNewt} in Sec.~\ref{sec:SNM} and a new specialized variant for solving GLMs in Sec.~\ref{sec:GLM}.

 Here we present the new nonlinear Kaczmarz method as a variant of \SNR
 .
Consider 
the original problem~\eqref{eq:main}.
We use a single row importance weighted subsampling sketch 
to sample rows of $F(x)=0.$ That is, let $\Prob{\mS = e_i} = p_i$ %
with the $i$th unit coordinate vector $e_i\in\R^m$ for $i = 1, \cdots, m$.
Then the \SNR update~\eqref{eq:update} is given by
\begin{equation}\label{eq:nonkaczmarz}
x^{k+1} = x^k - \gamma \frac{F_i(x^k)}{\norm{\nabla F_i(x^k)}^2} \nabla F_i(x^k).
\end{equation}
We dub~\eqref{eq:nonkaczmarz} the \emph{nonlinear Kaczmarz method}, as it can be seen as an extension of the randomized Kaczmarz method~\cite{Kaczmarz1937,Strohmer2009} for solving linear systems to the nonlinear case\footnote{We note that there exists a nonlinear variant of the Kaczmarz method which is referred to as the Landweber{\textendash}Kaczmarz~method \cite{Leit_o_2016}. Though the  Landweber{\textendash}Kaczmarz is very similar to Kaczmarz, it is not truly an extension since it does not adaptively re-weight the stepsize by $\norm{\nabla F_i(x^k)}^2$. }.
By~\eqref{eq:fk}, this nonlinear Kaczmarz method is simply \SGD applied to minimizing
\xline{
\[f_{x^k}(x) \; = \; \sum_{i=1}^m \Prob{\mS = e_i} f_{e_i, x^k}(x) \; \overset{\eqref{eq:fk}+\eqref{eq:Hk}}{=} \; \frac{1}{2}\sum_{i=1}^m p_i \frac{F_i(x)^2}{\norm{\nabla F_i(x^k)}^2}. \]
} 

A sufficient condition for~\eqref{ass:ker} to hold is that the diagonal matrix
\begin{align}
\EE{e_i}{\mH_{e_i}(x^k)} \overset{\eqref{eq:Hk}}{=} \sum_{i=1}^m p_i \frac{ e_ie_i^\top}{\norm{\nabla F_i(x^k)}^2}  = \Diag{\frac{p_i}{\norm{\nabla F_i(x^k)}^2}}
\end{align}
is invertible. Thus $\E{\mH_{\mS}(x^k)}$ is invertible if $\nabla F_i(x^k) \neq 0$ for all $i \in \{1, \cdots, m \}$ and $x^k \in \R^p$. In which case ${\bf Ker}\left(\E{\mH_{\mS}(y)}\right) =\{0\}$ for all $y \in \R^p$ and~\eqref{ass:ker} holds.

Finally, to guarantee that~\eqref{eq:nonkaczmarz} converges through Thm.~\ref{theo:convex}, we need $f_{x^k}(x)$ to be star-convex on $x^k$ at every iteration.
In this case, it suffices for each $F_i(x)^2$ to be star-convex, since any conic combination of star-convex functions is star-convex~\cite{LeeV16}. This is a straightforward abstraction of the one dimension case, in that, if~\eqref{eq:unicond} holds for every $F_i$ in the place of $\phi$, we can guarantee the convergence of~\eqref{eq:nonkaczmarz}. \ruilinex{This is also equivalent to assuming the star-convexity for each sketching matrix~\eqref{eq:cvx-S-star}. Furthermore, if we have $F(x)\in{\bf Im}(DF(x)^\top)$ hold for all $x$, then~\eqref{eq:DFTinEmHass} holds, as ${\bf Ker}\left(\E{\mH_{\mS}(y)}\right) =\{0\}$. We can guarantee the convergence of~\eqref{eq:nonkaczmarz} through Cor.~\ref{lem:EuclidieanF}.}

\section{The Stochastic Newton method}
\label{sec:SNM}


We now show that the Stochastic Newton method (\SNM)~\cite{rodomanov16superlinearly,Kovalev2019stochNewt} is a special case of \SNR. 
This connection combined with the global convergence theory of \SNR, gives us the first global convergence theory of \SNM, which we detail in Sec.~\ref{sec:globalSNM}.

\SNM~\cite{Kovalev2019stochNewt} is a stochastic second order method that takes a Newton-type step at each iteration to solve optimization problems with a finite-sum structure
\begin{equation} \label{eq:fSNM}
\min_{w\in \R^d} \left[P(w) \eqdef \frac{1}{n}\sum_{i=1}^n\phi_i(w)\right],
\end{equation}
where each $\phi_i:\R^d \rightarrow \R$ is twice differentiable and strictly convex. In brevity, the updates in \SNM at the $k$th iteration are given by
\begin{align}
w^{k+1} \; &= \; \left(\frac{1}{n}\sum_{i=1}^n \nabla^2 \phi_i(\alpha^k_i)\right)^{-1} \left(\frac{1}{n}\sum_{i=1}^n \nabla^2 \phi_i(\alpha^k_i) \alpha^k_i - \frac{1}{n}\sum_{i=1}^n \nabla \phi_i(\alpha^k_i)\right), \label{eq:x} \\
\alpha^{k+1}_i \; &= \; \label{eq:w}
\begin{cases}
w^{k+1} & \quad \mbox{if } i \in B_n \\
\alpha^k_i & \quad \mbox{if } i \notin B_n
\end{cases},
\end{align}
where $\alpha^k_1, \cdots, \alpha^k_n$ are auxiliary variables, initialized in \SNM, and $B_n \subset \{1, \ldots, n\}$ is a subset of size $\tau$ chosen uniformly on average from all subsets of size $\tau$.

\subsection{Rewrite \SNM as a special case of \SNR}

Since $P(w)$ is strictly convex, every minimizer of $P$ satisfies
$\nabla P(w) = \frac{1}{n} \sum_{i=1}^n \nabla \phi_i(w) = 0.$ Our main insight to deducing \SNM is that we can re-write this stationarity condition using a \emph{variable splitting trick}.
That is, by introducing a new variable $\alpha_i \in \R^d$ for each gradient $\nabla \phi_i$, 
and let $p:=(n+1)d$ and 
 $x = \begin{bmatrix}
w \ ;
\alpha_1 \ ;
\cdots \ ;
\alpha_n
 \end{bmatrix}\in \R^{p}$ 
 be the stacking\footnote{\ruiline{In this paper, vectors are columns by default, and given $x_1, \dots, x_n \in \mathbb{R}^q$, we note $[x_1 ; \dots ; x_n]\in \mathbb{R}^{qn}$ the (column) vector stacking the $x_i$'s on top of each other with $q \in \N$.}} of the $w$ and $\alpha_i$ variables,
we have that solving $\nabla P(w) = 0$ is equivalent to 
finding the \emph{roots} of the following nonlinear equations
\begin{equation}\label{eq:Fxwi}
F(x) \; = \; F(w \ ;
\alpha_1 \ ;
\cdots \ ;
\alpha_n) \; \eqdef \;
\left[
\frac{1}{n}\sum_{i=1}^n \nabla \phi_i(\alpha_i) \ ; \
w - \alpha_1 \ ; \
\cdots \ ; \
w - \alpha_n
\right],
\end{equation}
where $F:\R^{(n+1)d} \rightarrow \R^{(n+1)d}$. \ruiline{Our objective now becomes solving $F(x) = 0$ with $p = m = (n+1)d$.} To apply \SNR to \eqref{eq:Fxwi}, we are going to use a structured sketching matrix. But first, we need some additional notations.

Divide $\mI_{nd} \in \R^{nd \times nd}$ into $n$ contiguous blocks of size $nd \times d $ as follows 
\begin{eqnarray*}
\mI_{nd} &\eqdef& [ \ \mI_{nd, 1} \ \mI_{nd, 2} \cdots \mI_{nd, n} \ ]
\end{eqnarray*}
where $\mI_{nd, i}$ is the $i$th block of $\mI_{nd}$. Let $B_n \subset \{1,\ldots, n\}$ with $|B_n| = \tau$ chosen uniformly at average. Let $\mI_{B_n} \in \R^{nd \times \tau d}$ denote the concatenation of the blocks $\mI_{nd, i}$ such that the indices $i \in B_n$.

At the $k$th iteration of \SNR, denote $x^k = [w^k; \ \alpha^k_1; \ \cdots; \ \alpha^k_n]$, we define our sketching matrix $\mS \sim \cD_{x^k}$ as
\begin{eqnarray} \label{eq:SNMsketch}
\mS &=& 
\begin{bmatrix}
\mI_d & 0 \\
	\begin{matrix}
	\frac{1}{n}\nabla^2 \phi_1(\alpha^k_1) \\
	\vdots \\
	\frac{1}{n}\nabla^2 \phi_n(\alpha^k_n)
	\end{matrix}
& \mI_{B_n}
\end{bmatrix} \in \R^{(n+1)d \times (\tau+1)d}.
\end{eqnarray}
Here the distribution $\cD_{x^k}$ depends on the iterates $x^k$. \ruiline{The sketch size of $\mS$ is $(\tau+1)d$ with \emph{any} $\tau \in \{1, \cdots, n\}$.} Now we can state the following lemma.

\begin{lemma} \label{lem:SNMdetail}
Let $\phi_i$ be strictly convex for $i=1,\ldots, n.$
At each iteration $k$, the updates of \SNR \eqref{eq:update} with $F$ defined in \eqref{eq:Fxwi}, the sketching matrix $\mS_k$ defined in \eqref{eq:SNMsketch} , and stepsize $\gamma = 1$, are equal to the updates \eqref{eq:x} and \eqref{eq:w} of \SNM.
\end{lemma}

In our upcoming proof of Lemma~\ref{lem:SNMdetail}, we still need the following lemma.

\begin{lemma} \label{lem:DFinv}
Let $\phi_i$ be twice differentiable and strictly convex for $i=1,\ldots, n.$
The Jacobian $DF(x)^\top$ of $F(x)$ defined in \eqref{eq:Fxwi} is invertible for all $x \in \R^{(n+1)d}$.
\end{lemma}

\begin{proof}
Let $x \in \R^{(n+1)d}$. Let $y \eqdef (u; v_1; \cdots; v_n) \in \R^{(n+1)d}$ with $u, v_1, \cdots, v_n \in \R^d$ such that $DF(x)y = 0$.
The transpose of the Jacobian of $F(x)$ is given by
\begin{eqnarray} \label{eq:DFxwi}
DF(x) &=&
\begin{bmatrix}
0 & \mI_d & \cdots & \mI_d \\
	\begin{matrix}
	\frac{1}{n}\nabla^2\phi_1(\alpha_1) \\
	\vdots \\
	\frac{1}{n}\nabla^2\phi_n(\alpha_n)
	\end{matrix}
& & -\mI_{nd}
\end{bmatrix}.
\end{eqnarray}
From $DF(x)y = 0$ and \eqref{eq:DFxwi}, we obtain
\begin{align*}
\sum_{i=1}^n v_i = 0 , \quad \mbox{ and } \quad
\frac{1}{n}\nabla^2\phi_i(\alpha_i)u = v_i \quad \mbox{for all } i = 1, \cdots, n.
\end{align*}
Plugging the second equation in the first one gives $\left(\frac{1}{n}\sum_{i=1}^n\nabla^2\phi_i(\alpha_i)\right)u = 0$. Since every $\phi_i$ is twice differentiable and strictly convex, we have $\nabla^2\phi_i(\alpha_i) > 0$. This implies $\frac{1}{n}\sum_{i=1}^n\nabla^2\phi_i(\alpha_i) > 0$, and is thus invertible. Consequently $u = 0$ and $v_i = 0$  from which we conclude that the Jacobian $DF(x)^\top$ is invertible.
\end{proof}

Now we can give the proof of Lemma~\ref{lem:SNMdetail}.
\begin{proof}
Consider an update of \SNR \eqref{eq:update} with $F$ defined in \eqref{eq:Fxwi}, the sketching matrix $\mS_k$ defined in \eqref{eq:SNMsketch} , and stepsize $\gamma = 1$ at the $k$th iteration. By Lemma \ref{lem:DFinv}, we have that $DF(x)$ is invertible and thus Asm.~\ref{ass:existence} holds. By \eqref{eq:sketch}, the \SNR update \eqref{eq:update} can be re-written as
\begin{align}
x^{k+1} 
= \argmin \norm{w - w^k}^2 + \sum_{i=1}^n\norm{\alpha_i - \alpha^k_i}^2 \ \mbox{s.t.} \ \mS_k^\top DF(x^k)^\top
(x - x^k)
= - \mS_k^\top F(x^k). \label{eq:newtonprojF}
\end{align}
Plugging \eqref{eq:Fxwi}, \eqref{eq:SNMsketch} and \eqref{eq:DFxwi} into the constraint in~\eqref{eq:newtonprojF} 
gives
\begin{align*}
&\quad \quad \begin{bmatrix}
\mI_d & 0 \\
	\begin{matrix}
	\frac{1}{n}\nabla^2 \phi_1(\alpha^k_1) \\
	\vdots \\
	\frac{1}{n}\nabla^2 \phi_n(\alpha^k_n)
	\end{matrix}
& \mI_{B_n}
\end{bmatrix}^\top
\begin{bmatrix}
0 & \mI_d & \cdots & \mI_d \\
	\begin{matrix}
	\frac{1}{n}\nabla^2\phi_1(\alpha_1) \\
	\vdots \\
	\frac{1}{n}\nabla^2\phi_n(\alpha_n)
	\end{matrix}
& & -\mI_{nd}
\end{bmatrix}^\top
\begin{bmatrix}
w - w^k \\ \alpha_1 - \alpha^k_1 \\ \vdots \\ \alpha_n - \alpha^k_n
\end{bmatrix} \\
&= -
\begin{bmatrix}
\mI_d & 0 \\
	\begin{matrix}
	\frac{1}{n}\nabla^2 \phi_1(\alpha^k_1) \\
	\vdots \\
	\frac{1}{n}\nabla^2 \phi_n(\alpha^k_n)
	\end{matrix}
& \mI_{B_n}
\end{bmatrix}^\top
\begin{bmatrix}
\frac{1}{n}\sum_{i=1}^n \nabla \phi_i(\alpha^k_i) \\ w^k - \alpha^k_1 \\ \vdots \\ w^k - \alpha^k_n
\end{bmatrix}.
\end{align*}
After 
 simplifying the above matrix multiplications, we have that~\eqref{eq:newtonprojF} is given by
\begin{align}
x^{k+1} = [w^{k+1}; \alpha^{k+1}_1; \cdots; \alpha^{k+1}_n] &= \argmin \norm{w-w^k}^2 +\sum_{i=1}^n\norm{\alpha_i-\alpha^k_i}^2  \nonumber\\
& \mbox{ s. t.} \quad \frac{1}{n}\sum_{i=1}^n \nabla^2 \phi_i(\alpha^k_i)(w - \alpha^k_i) = -\frac{1}{n}\sum_{i=1}^n \nabla \phi_i(\alpha^k_i), \nonumber\\
&\phantom{\mbox{ s. t.}} \quad w = \alpha_j, \quad \mbox{ for } j \in B_n. \label{eq:projxwiNR}
\end{align}
To solve \eqref{eq:projxwiNR}, first note that $\alpha^{k+1}_i = \alpha^k_i$ for $i \not\in B_n$, since there is no constraint on the variable $\alpha_i$ in this case.
Furthermore, by the invertibility of $\frac{1}{n}\sum_{i=1}^n \nabla^2 \phi_i(\alpha^k_i)$, we have that \eqref{eq:projxwiNR} has a unique solution s.t.\ $\alpha_j = w$ for all $j \in B_n$ and
\begin{align*}
w = \left(\frac{1}{n}\sum_{i=1}^n \nabla^2 \phi_i(\alpha^k_i)\right)^{-1}\left(\frac{1}{n}\sum_{i=1}^n \nabla^2 \phi_i(\alpha^k_i) \alpha^k_i-\frac{1}{n}\sum_{i=1}^n \nabla \phi_i(\alpha^k_i)\right).
\end{align*}
Concluding, we have that the \SNR update~\eqref{eq:projxwiNR} is given by
\begin{align*}
w^{k+1} \; &= \; \left(\frac{1}{n}\sum_{i=1}^n \nabla^2 \phi_i(\alpha^k_i)\right)^{-1} \left(\frac{1}{n}\sum_{i=1}^n \nabla^2 \phi_i(\alpha^k_i) \alpha^k_i - \frac{1}{n}\sum_{i=1}^n \nabla \phi_i(\alpha^k_i)\right),  \\
\alpha^{k+1}_i \; &= \; 
\begin{cases}
w^{k+1} & \quad \mbox{if } i \in B_n \\
\alpha^k_i & \quad \mbox{if } i \notin B_n
\end{cases},
\end{align*}
which is exactly the Stochastic Newton method's updates \eqref{eq:x} and \eqref{eq:w} in~\cite{Kovalev2019stochNewt}.
\end{proof}
\robline{Thus we conclude that \SNM is a special case of \SNR. However, in practice for solving GLMs, instead of sampling $\mS \sim \cD_{x^k}$ provided in~\eqref{eq:SNMsketch}, we only sample $B_n$ and we execute the efficient updates as suggested in~\cite{Kovalev2019stochNewt}.}

\subsection{Global convergence theory of \SNM} \label{sec:globalSNM}

Let $x' \eqdef (w'; \alpha'_1; \cdots; \alpha'_n) \in \R^{(n+1)d}$ and $\mS\sim\cD_{x}$ defined in \eqref{eq:SNMsketch}.
By applying the global convergence theory of \SNR, we can now provide the first global convergence theory for \SNM.


\begin{corollary} \label{cor:SNM1}
Let $w^*$ be a solution to $\nabla P(w) = 0$. Consider the iterate $x^k = (w^k; \alpha^k_1; \cdots; \alpha^k_n)$ given by \SNM \eqref{eq:x} and \eqref{eq:w} and note $x^* \eqdef (w^*; w^*; \cdots; w^*) \in \R^{(n+1)d}$. If there exists $\mu > 0$ such that for all $x$, $x' \in \R^{(n+1)d}$,
\begin{align} \label{eq:strconvSNM}
f_{x'}(x^*) \; &\geq \; f_{x'}(x) + \dotprod{\nabla f_{x'}(x), x^* - x} + \frac{\mu}{2} \norm{x^* - x}^2 \\
\; &= \; f_{x'}(x) + \dotprod{\nabla f_{x'}(x), x^* - x} + \frac{\mu}{2} \bigg(\norm{w^*-w}^2 + \sum_{i=1}^n\norm{w^*-\alpha_i}^2\bigg), \nonumber
\end{align}
then the iterates $\{x^k\}$ of \SNM converge linearly according to
\begin{equation}\label{eq:linearconvSNM}
\E{\norm{x^{k+1}-x^*}^2}  \leq (1 - \mu)^{k+1} \norm{x^0-x^*}^2.
\end{equation}
\end{corollary}
\begin{proof}
As $\nabla P(w^*) = 0$, this implies immediately that $x^*$ is a solution of $F$. Besides, \eqref{eq:strconvSNM} satisfies Asm. \ref{ass:strconv}. Thus by Thm. \ref{theo:strconv}, we get \eqref{eq:linearconvSNM}. 
\end{proof}

\ruilinex{Even though~\eqref{eq:strconvSNM} is a strong assumption,} this is the first global convergence theory of \SNM, since only local convergence results of \SNM are addressed in~\cite{Kovalev2019stochNewt}.

\ruiline{
As a by-product, we find that the function $F(x)$ in \eqref{eq:Fxwi} and the sketch $\mS$ defined in \eqref{eq:SNMsketch} actually satisfy~\eqref{eq:DFTinEmHass} through Lemma~\ref{lem:reformulationnew}, namely as the following lemma.
\begin{lemma} \label{lem:adaptedsketchSNM} 
Consider the function $F$ defined in \eqref{eq:Fxwi} and the sketching matrix $\mS$ defined in \eqref{eq:SNMsketch}, then
we have the condition~\eqref{eq:DFTinEmHass} hold.
\end{lemma}
}
\ruiline{
\begin{proof}
First, we show that $\E{\mS\mS^\top}$ is invertible $\forall x \in \R^{(n+1)d}$. By the definition of $\mS$ in \eqref{eq:SNMsketch},
\begin{equation*}
\mS\mS^\top \; = \;
\begin{bmatrix}
\mI_d & \frac{1}{n}\nabla^2\phi_1(\alpha_1) & \cdots & \frac{1}{n}\nabla^2\phi_n(\alpha_n) \\
    \begin{matrix}
    \frac{1}{n}\nabla^2\phi_1(\alpha_1) \\
    \vdots \\
    \frac{1}{n}\nabla^2\phi_n(\alpha_n)
    \end{matrix}
& & \mI_{B_n}\mI_{B_n}^\top + \mM
\end{bmatrix},
\end{equation*}
where $\mM = \{\mM_{ij}\}_{1 \leq i \leq n, 1 \leq j \leq n}$ is divided into $n \times n$ contiguous blocks of size $d \times d$ with each block $\mM_{ij}$ defined as the following
\begin{equation*}
\mM_{ij} \eqdef \frac{1}{n}\nabla^2\phi_i(\alpha_i)\cdot\frac{1}{n}\nabla^2\phi_j(\alpha_j) \in \R^{d \times d} \quad \quad \mbox{ and } \quad \quad \mM \in \R^{nd \times nd}.
\end{equation*}
Taking the expectation over $\mS$ w.r.t.\ the distribution $\cD_{(w, \alpha_1, \cdots, \alpha_n)}$ gives
\begin{align} \label{eq:SST}
\E{\mS\mS^\top} \; &= \;
\begin{bmatrix}
\mI_d & \frac{1}{n}\nabla^2\phi_1(\alpha_1) & \cdots & \frac{1}{n}\nabla^2\phi_n(\alpha_n) \\
    \begin{matrix}
    \frac{1}{n}\nabla^2\phi_1(\alpha_1) \\
    \vdots \\
    \frac{1}{n}\nabla^2\phi_n(\alpha_n)
    \end{matrix}
& & \frac{\tau}{n}\mI_{nd} + \mM
\end{bmatrix} \nonumber \\
\; &= \;
\begin{bmatrix}
\mI_d \\ \frac{1}{n}\nabla^2\phi_1(\alpha_1) \\ \vdots \\ \frac{1}{n}\nabla^2\phi_n(\alpha_n)
\end{bmatrix}
\begin{bmatrix}
\mI_d \\ \frac{1}{n}\nabla^2\phi_1(\alpha_1) \\ \vdots \\ \frac{1}{n}\nabla^2\phi_n(\alpha_n)
\end{bmatrix}^\top + \frac{\tau}{n}
\begin{bmatrix}
0 & 0 \\
0 & \mI_{nd}.
\end{bmatrix}
\end{align}
$\E{\mS\mS^\top}$ is symmetric, positive semi-definite. Let $(u; v_1; \cdots; v_n) \in \R^{(n+1)d}$ s.t.\
\[(u; v_1; \cdots; v_n)^\top\E{\mS\mS^\top}(u; v_1; \cdots; v_n) = 0.\]
From \eqref{eq:SST}, we obtain
\begin{align*}
\left\|\left[\mI_d ; \frac{1}{n}\nabla^2\phi_1(\alpha_1) ; \cdots ; \frac{1}{n}\nabla^2\phi_n(\alpha_n)\right]^\top [u ; v_1; \cdots; v_n]\right\|^2 
+ \frac{\tau}{n}\sum_{i=1}^n\norm{v_i}^2 = 0.
\end{align*}
Since both terms are non negative, we obtain $\sum_{i=1}^n\norm{v_i}^2 = 0 \implies \forall \ i, v_i = 0$, and then $u = 0$. This confirms that $\E{\mS\mS^\top}$ is positive-definite, thus invertible and ${\bf Ker}(\E{\mS\mS^\top}) = \{0\}$.
Besides, from Lemma~\ref{lem:DFinv}, we get $DF(x)$ invertible. Thus $F(x) \in {\bf Im}(DF(x)^\top)$ and ${\bf Ker}(DF(x)) = \{0\}$. We have~\eqref{eq:adaptedsketch} hold. By Lemma~\ref{lem:reformulationnew}, we have that~\eqref{eq:DFTinEmHass} holds for all $x \in \R^{(n+1)d}$.
\end{proof}}


\ruiline{From Lemma~\ref{lem:adaptedsketchSNM}, we know that for \emph{any} size of the subset sampling $|B_n| = \tau \in \{1, \cdots, n\}$, the condition~\eqref{eq:DFTinEmHass} holds. The corresponding sketch size of $\mS$ is $(\tau+1)d$.}

Furthermore, using Lemma~\ref{lem:SNMdetail}, we can also provide the global convergence of \SNM with stepsizes $\gamma < 1$ (\SNM with relaxation) by using the weaker star-convexity assumption 
 and Thm.~\ref{theo:convex}.
We expand on this comment in App.~\ref{sec:SNMwRelax} in the appendix. 


\section{Applications to GLMs -- \emph{tossing-coin-sketch} method} \label{sec:GLM}

Consider the problem of training a generalized linear model
\begin{equation} \label{eq:P}
w^* = \arg\min_{w \in \R^d} P(w) \, \eqdef \, \frac{1}{n}\sum_{i=1}^n \phi_i(a_i^\top w) + \frac{\lambda}{2}\|w\|^2,
\end{equation}
where $\phi_i:\R \rightarrow \R^+$ is a convex and continuously twice differentiable loss function, $a_i \in \R^d$ are data samples and $w \in \R^d$ is the parameter to optimize. As the objective function is strongly convex, the unique minimizer satisfies $\nabla P(w) = 0$, that is
\begin{align}\label{eq:nablaP}
\nabla P(w) &= \frac{1}{n}\sum_{i=1}^n\phi_i'(a_i^\top w)a_i + \lambda w = 0.
\end{align}
Let
$
\Phi(w) \;\eqdef\;
\begin{bmatrix}
\phi_1'(a_1^\top w) \ \cdots \ \phi_n'(a_n^\top w)
\end{bmatrix}^\top \;\in \;\R^{n}
$
and
$ 
\mA  \;\eqdef\;
\begin{bmatrix}
a_1 \ \cdots \ a_n
\end{bmatrix} \;\in \; \R^{d \times n}.
$
By introducing auxiliary variables $\alpha_i \in \R$ s.t.\  $\alpha_i \eqdef -\phi_i'(a_i^\top w)$, we can re-write~\eqref{eq:nablaP} as
\begin{equation} 
w \; = \; \frac{1}{\lambda n}\mA\alpha,
\quad \quad \mbox{and } \quad \quad 
\alpha \; = \; -\Phi(w). \label{eq:fixedPhi} 
\end{equation}
Note $x = [\alpha; w] \in \R^{n+d}$.
The objective of finding the minimum of \eqref{eq:P} is now equivalent to finding \emph{zeros} for the function 
\begin{equation}\label{eq:F}
F(x) \; = \; F(\alpha; w) \; \eqdef \;
\begin{bmatrix}
\frac{1}{\lambda n}\mA\alpha - w \\
\alpha + \Phi(w)
\end{bmatrix},
\end{equation}
where $F:\R^{n+d} \rightarrow \R^{n+d}$. \ruiline{Our objective now becomes solving $F(x) = 0$ with $p = m = n+d$.} For this, we will use a variant of the \SNR. The advantage in representing~\eqref{eq:nablaP} as the nonlinear system~\eqref{eq:F} is that we now have one row per data point (see the second equation in~\eqref{eq:fixedPhi}). This allows us to use sketching to \emph{subsample} the data. 

Since the function $F$ has a block structure, we will use a structured sketching matrix which we refer to as a \emph{Tossing-coin-sketch}. But first, we need the following definition of a block sketch.

\begin{definition}[($n,\tau$)--block sketch]
Let $B_n \subset \{1,\ldots, n\}$ be a subset of size $\tau$ uniformly sampling at random. We say that $\mS \in \R^{n\times \tau}$ is a $(n,\tau)$--block sketch if $\mS = \mI_{B_n}$ where $\mI_{B_n}$ denotes the column concatenation of the columns of the identity matrix $\mI_n \in \R^{n\times n}$ whose indices are in  $B_n$.
\end{definition}

Our Tossing-coin-sketch is a sketch that alternates between two blocks depending on the result of a coin toss.

\begin{definition}[Tossing-coin-sketch]\label{def:TCS}
Let $\mS_d \in \R^{d \times \tau_d}$ and $\mS_n \in \R^{n \times \tau_n}$ be a $(d,\tau_d)$--block sketch and a  $(n,\tau_n)$--block sketch, respectively. Let $b \in (0; 1)$. 
Now each time we sample $\mS$, we ``toss a coin'' to determine the structure of $\mS  \in \R^{(d+n)\times(\tau_d+\tau_n)}$. That is, $\mS = \begin{bmatrix}
\mS_d & 0 \\
0 & 0
\end{bmatrix}$ with probability $1-b$ and $\mS = \begin{bmatrix}
0 & 0 \\
0 & \mS_n
\end{bmatrix}$ with probability $b$.
\end{definition}

By applying the \SNR method with a tossing-coin-sketch for solving~\eqref{eq:F}, we arrive at an efficient method for solving~\eqref{eq:P} that we call the \emph{TCS} method.
By using a tossing-coin-sketch, we can alternate between solving a linear system based on the first $d$ rows of~\eqref{eq:F} and a nonlinear system based on the last $n$ rows of~\eqref{eq:F}. 

\ruilinex{TCS is inspired by the first-order stochastic dual ascent methods~\cite{SDCA,Shalev-Shwartz2015,Qu2015b}. Indeed, eq.~\eqref{eq:fixedPhi} can be seen as primal-dual systems with primal variables $w$ and dual variables $\alpha$. Stochastic dual ascent methods are efficient to solve~\eqref{eq:fixedPhi}. At each iteration, they update alternatively the primal and the dual variables $w$ and $\alpha$ with the first-order informations. Thus, by sketching alternatively the primal and the dual systems and updating accordingly with the Newton-type steps, TCS's updates can be seen as the second-order stochastic dual ascent methods.}

We show in the next section that the TCS method verifies~\eqref{eq:DFTinEmHass}. 
Using sketch sizes s.t.\ $\tau_n \ll n$,
the TCS method has the same cost as \SGD in the case $d \ll n$. 
\ruiline{The low computational cost per iteration is thus another advantage of the TCS method.} 
See App.~\ref{sec:detailexp} the cost per iteration analysis.
For a detailed derivation of the TCS method, see App.~\ref{sec:equivalencealgorithms} and a detailed implementation in Algorithm~\ref{algo:tau-TCS} in the appendix.

\subsection{The condition~\texorpdfstring{\eqref{eq:DFTinEmHass}}{} in the case of TCS method}
\label{sec:kerTCS}

In this section, we show that the TCS method verifies~\eqref{eq:DFTinEmHass} through Lemma~\ref{lem:reformulationnew} in the following.
\ruiline{
\begin{lemma} \label{lem:adaptTCS} 
Consider the function $F$ defined in \eqref{eq:F} and the tossing-coin-sketch $\mS$ defined in Definition~\ref{def:TCS}, then~\eqref{eq:DFTinEmHass} holds.
\end{lemma}}

\begin{proof}
First, we show that $\E{\mS\mS^\top}$ is invertible.
By Definition~\ref{def:TCS}, it is straightforward to verify that
\begin{eqnarray*}
\E{\mS\mS^\top} &=& \begin{bmatrix} \frac{(1-b)\tau_d}{n}\mI_d & 0 \\ 0 & \frac{b\tau_n}{n}\mI_n \end{bmatrix}
\end{eqnarray*}
is invertible and ${\bf Ker}(\E{\mS\mS^\top}) = \{0\}$.
Now we show the Jacobian $DF^\top(x)$ invertible. Let $x = [\alpha; w] \in \R^{n+d}$ with $\alpha \in \R^n$ and $w \in \R^d$. Then $DF(x)$ is written as
\begin{eqnarray} \label{eq:DFalphax}
DF(x)^\top &=&
\begin{bmatrix}
\frac{1}{\lambda n}\mA & - \mI_d \\
\mI_n & \nabla \Phi(w)^\top
\end{bmatrix},
\end{eqnarray}
where $\nabla \Phi(w)^\top = \Diag{\phi_1''(a_1^\top w), \ldots, \phi_n''(a_n^\top w)} \,\mA^\top \in \R^{n \times d}.$ 
Denote the diagonal matrix $D(w) \eqdef \Diag{\phi_1''(a_1^\top w), \ldots, \phi_n''(a_n^\top w)}$. Since $\phi_i$ is continuously twice differentiable and convex, $\phi_i''(a_i^\top w) \geq 0$ for all $i$. Thus, $D(w) \geq 0$.

Let $(u; v) \in \R^{n+d}$ with $u \in \R^n$ and $v \in \R^d$ such that $DF(x)^\top [u; v] = 0$. We have
\begin{align}
DF(x)^\top \begin{bmatrix}u \\ v\end{bmatrix} = 0 \overset{\eqref{eq:DFalphax}}{\Longleftrightarrow}
\begin{bmatrix}
\frac{1}{\lambda n}\mA & - \mI_d \\
\mI_n & \nabla \Phi(w)^\top
\end{bmatrix} \begin{bmatrix}u \\ v\end{bmatrix} = 0
\Longrightarrow \left(\mI_n + \frac{1}{\lambda n}D(w)\mA^\top \mA\right) u = 0. \label{eq:IDAA}
\end{align}

If $D(w)$ is invertible, \eqref{eq:IDAA} becomes
\begin{eqnarray}
D(w)\left(D(w)^{-1} + \frac{1}{\lambda n} \mA^\top \mA\right) u = 0
&\Longleftrightarrow& \left(D(w)^{-1} + \frac{1}{\lambda n} \mA^\top \mA\right) u = 0. \label{eq:DAA}
\end{eqnarray}
Since $D(w)$ is invertible, i.e.\ $D(w) > 0$, we obtain $D(w)^{-1} > 0$. As $\frac{1}{\lambda n} \mA^\top \mA \geq 0$, we get $D(w)^{-1} + \frac{1}{\lambda n} \mA^\top \mA > 0$, thus invertible. From \eqref{eq:DAA}, we get $u = 0$.

Otherwise, $D(w)$ is not invertible. Without losing generality, we assume that $\phi_1''(a_1^\top w) \geq \phi_2''(a_2^\top w) \geq \cdots \geq \phi_n''(a_n^\top w) = 0$. Let $j$ be the largest index for which $\phi_j''(a_j^\top w) > 0$. If $j$ does not exist, then $D(w) = 0$. From \eqref{eq:IDAA}, we get $u=0$ directly. If $j$ exists, we have $1 \leq j < n$ and
\begin{align}
D(w)\mA^\top \mA \; &= \; \Diag{\phi_1''(a_1^\top w), \ldots, \phi_j''(a_j^\top w), 0, \ldots, 0} \mA^\top \mA \nonumber \\
\; &= \; \begin{bmatrix}
\Diag{\phi_1''(a_1^\top w), \ldots, \phi_j''(a_j^\top w)}\mA_{1:j}^\top \mA_{1:j} & 0 \\
0 & 0 \end{bmatrix}, \label{eq:DAAj}
\end{align}
where $\mA_{1:j} \eqdef [a_1 \ \cdots \ a_j] \in \R^{d \times j}$. Note $u = [u_1; \cdots; u_n] \in \R^n$. Plugging~\eqref{eq:DAAj} into~\eqref{eq:IDAA}, we get
\begin{eqnarray}
& & \left(\mI_n + \frac{1}{\lambda n}
\begin{bmatrix}
\Diag{\phi_1''(a_1^\top w), \ldots, \phi_j''(a_j^\top w)}\mA_{1:j}^\top \mA_{1:j} & 0 \\
0 & 0
\end{bmatrix} \right) u = 0 \nonumber \\
&\Longleftrightarrow& \begin{cases} \left(\mI_j + \frac{1}{\lambda n}\Diag{\phi_1''(a_1^\top w), \ldots, \phi_j''(a_j^\top w)}\mA_{1:j}^\top \mA_{1:j}\right)u_{1:j} = 0 \\ u_{(j+1):n} = 0 \end{cases} \label{eq:IDAAj},
\end{eqnarray}
where $u_{1:j} \eqdef \left[u_1; \cdots; u_j\right] \in \R^j$ and $u_{(j+1):n} \eqdef \left[u_{j+1}; \cdots; u_n\right] \in \R^{n-j}$. From~\eqref{eq:IDAAj}, $u_{(j+1):n} = 0$. Now $\Diag{\phi_1''(a_1^\top w), \ldots, \phi_j''(a_j^\top w)}$ is invertible in the subspace $\R^j$ as every coordinate in the diagonal $\phi_i''(a_i^\top w)$ is strictly positive for all $1 \leq i \leq j$. Similarly, we obtain $u_{1:j} = 0$ from the first equation of \eqref{eq:IDAAj}. Overall we get $u=0$.

Thus, in all cases, $u=0$, then $v = \frac{1}{\lambda n}Au = 0$. We can thus induce that $DF(\alpha; w)^\top$ is invertible for all $\alpha$ and $w$.
\ruiline{Similar to Lemma~\ref{lem:adaptedsketchSNM}, we have~\eqref{eq:adaptedsketch} hold, and by Lemma~\ref{lem:reformulationnew}, we have that~\eqref{eq:DFTinEmHass} holds.}
\end{proof}


\ruiline{From Lemma~\ref{lem:adaptTCS}, we know that for \emph{any} size of the block sketch $\tau_d \in \{1, \cdots, d\}$ and $\tau_n \in \{1, \cdots, n\}$,~\eqref{eq:DFTinEmHass} holds. The corresponding sketch size of $\mS$ is $\tau_d + \tau_n$.}

\subsection{Experiments for TCS method applied for GLM} \label{sec:experiments}

We consider the logistic regression problem with $8$ datasets\footnote{All datasets except for the artificial dataset can be found downloaded on \url{https://www.csie.ntu.edu.tw/~cjlin/libsvmtools/datasets/}. 
Some of the datasets can be found originally in~\cite{a9a,covtype,phishing,webspam,UCI}.
} taken from LibSVM~\cite{libsvm}, except for one artificial dataset. Table~\ref{tab:datasets} provides the details of these datasets, including the \emph{condition number} (C.N.) of the model and the smoothness constant $L$ of the model.
C.N.\ of the logistic regression problem is given by
\[\mbox{C.N.} \eqdef \frac{\lambda_{\max}(\mA\mA^\top)}{4n\lambda} + 1,\]
where  $\lambda_{\max}(\cdot)$ is the largest eigenvalue operator. 
The smoothness constant $L$ is given by
\[L \eqdef \frac{\lambda_{\max}(\mA\mA^\top)}{4n} + \lambda.\]
As for the logistic regression problem, we consider the loss function $\phi_i$ in \eqref{eq:P} in the form
\[\phi_i(t) = \ln(1+\mathrm{e}^{-y_it})\]
where $y_i$ are the target values for $i = 1, \cdots, n$.

\paragraph{The artificial dataset}
The artificial dataset $\mA^\top \in \R^{n \times d}$ in Table~\ref{tab:datasets} is of size $10000 \times 50$ and generated by a Gaussian distribution whose mean is zero and covariance is  a Toeplitz matrix. Toeplitz matrices are completely determined by their diagonal. We  set the diagonal of our Toeplitz matrix  as
\[[c^0; c^1; \cdots; c^{d-1}] \in \R^d\]
where $c \in \R^+$ is a parameter. We choose $c = 0.9$ (closed to $1$) which results in $\mA$ having highly correlated columns, which in turn makes $\mA$ an  ill-conditioned data set. We set the ground truth coefficients of the model
\[{\bf w} = [(-1)^0 \cdot \mathrm{e}^{-\frac{0}{10}}; \cdots; (-1)^{d-1} \cdot \mathrm{e}^{-\frac{d-1}{10}}] \in \R^d\]
and the target values of the dataset
\[{\bf y} = \text{sgn}\left(\mA^\top {\bf w} + {\bf r}\right) \in \R^n\]
where ${\bf r} \in \R^n$ is the noise generated from a standard normal distribution.

\begin{table}
  \caption{Details of the data sets for binary classification}
  \label{tab:datasets}
  \centering
  \begin{tabular}{ ccccc }
  	 \toprule
     dataset     & dimension ($d$)  & samples ($n$)  & C.N. of the model        & $L$                  \\
     \midrule
     covetype    & $54$             & $581012$       & $7.45 \times 10^{12}$  & $1.28 \times 10^7$          \\
     a9a         & $123$            & $32561$        & $5.12 \times 10^4$     & $1.57$                      \\
     fourclass   & $2$              & $862$          & $4.86 \times 10^6$     & $5.66 \times 10^3$          \\
     artificial  & $50$             & $10000$        & $3.91 \times 10^4$     & $3.91$                      \\
     ijcnn1      & $22$             & $49990$        & $2.88 \times 10^3$     & $5.77 \times 10^{-2}$       \\
     webspam     & $254$            & $350000$       & $7.47 \times 10^4$     & $2.13 \times 10^{-1}$       \\
     epsilon     & $2000$           & $400000$       & $3.51 \times 10^4$     & $8.76 \times 10^{-2}$       \\
     phishing    & $68$             & $11055$        & $1.04 \times 10^3$     & $9.40 \times 10^{-2}$       \\
     \bottomrule
  \end{tabular}
\end{table}

\begin{figure}
\centering
\begin{tabular}{cccc}
\includegraphics[width=.21\linewidth]{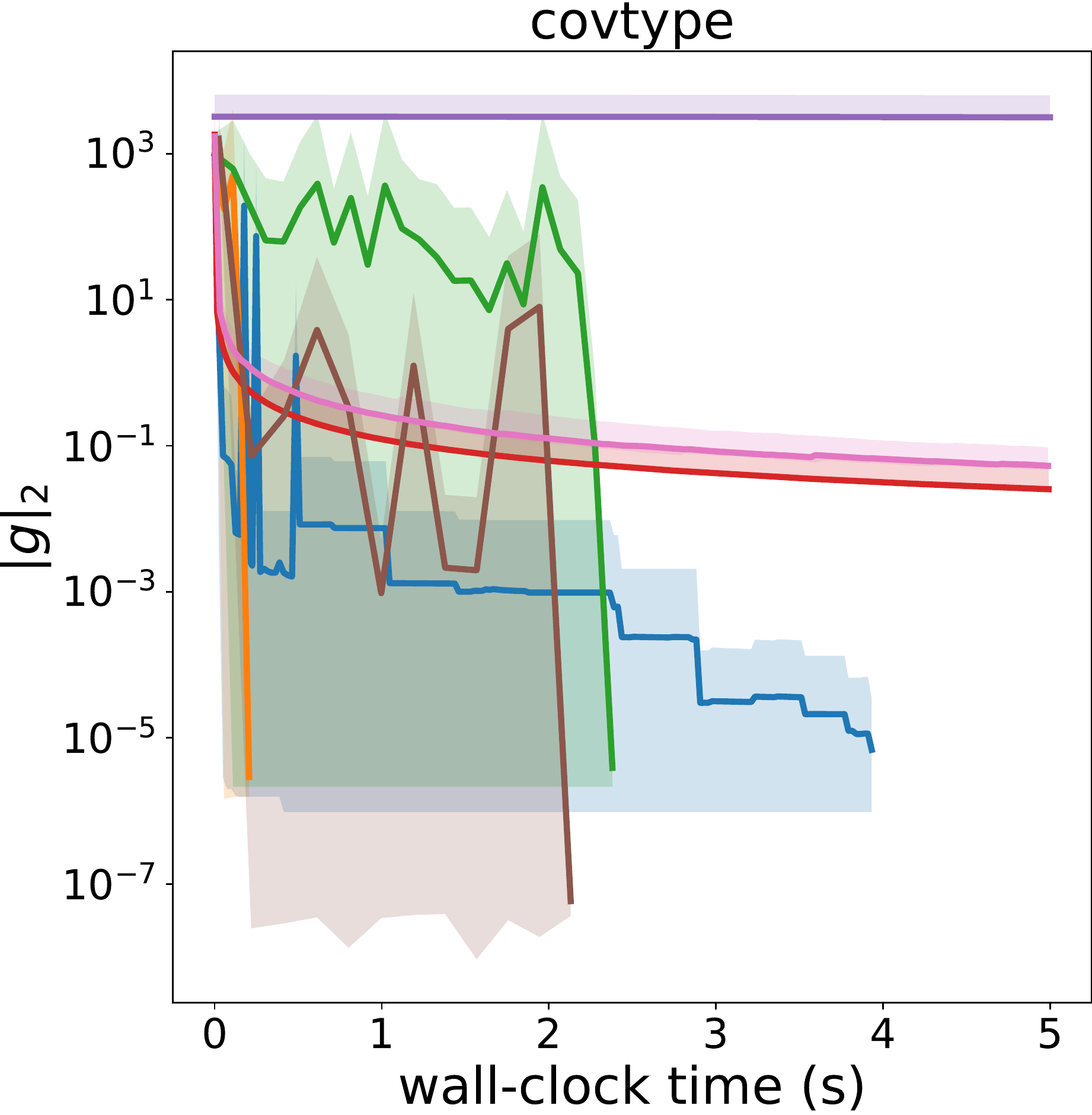}&
\includegraphics[width=.21\linewidth]{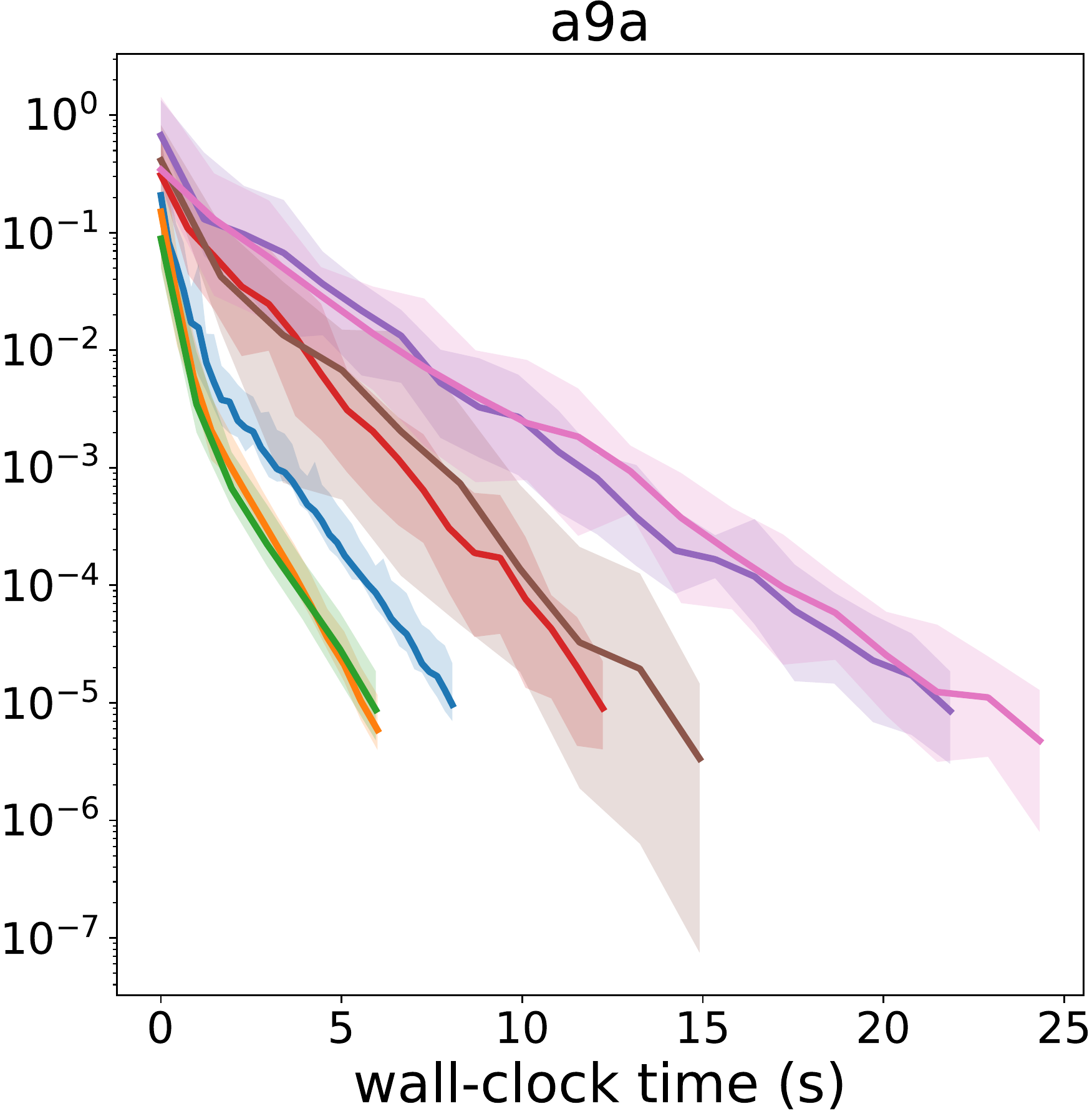}&
\includegraphics[width=.21\linewidth]{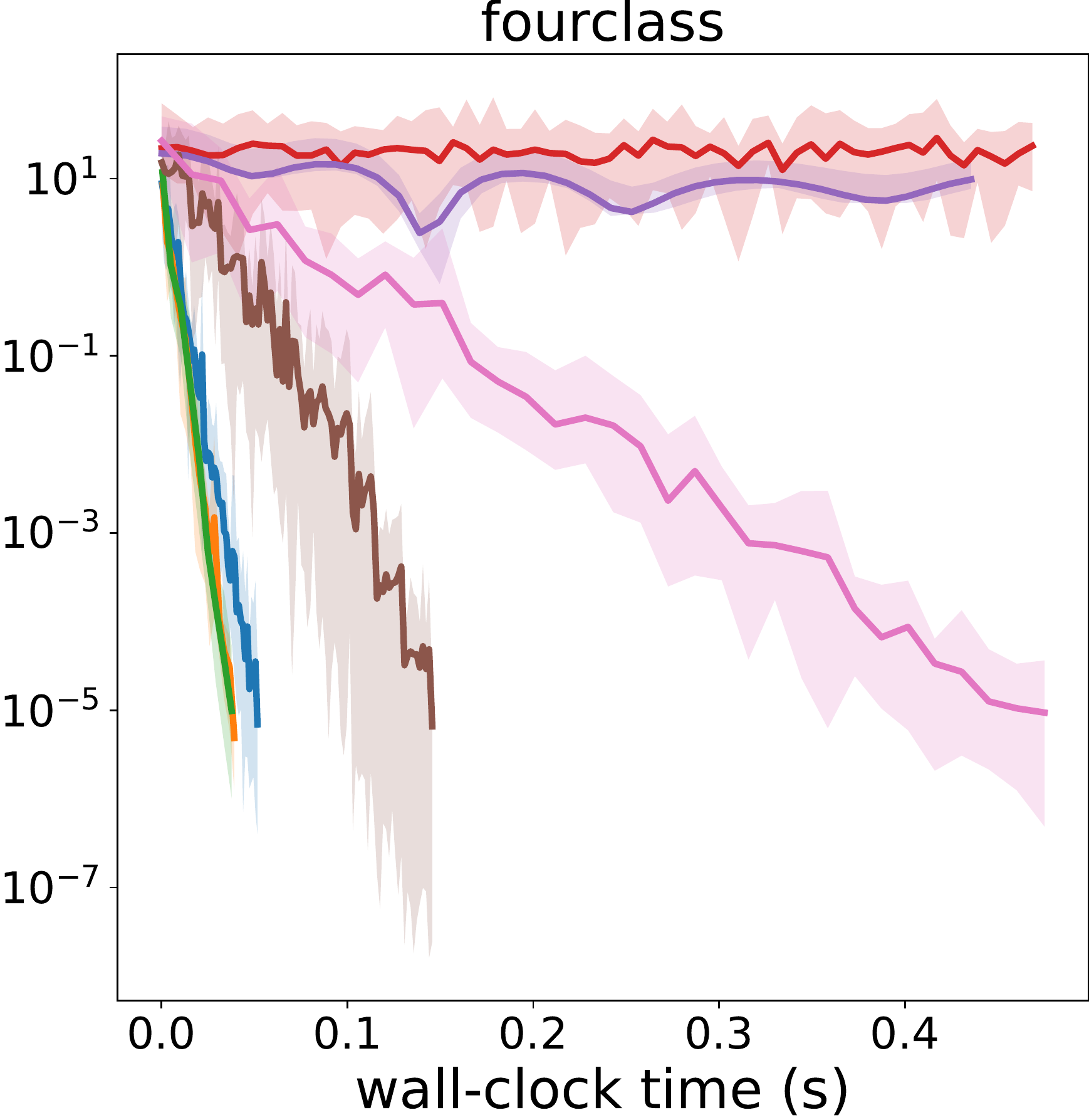}&
\includegraphics[width=.21\linewidth]{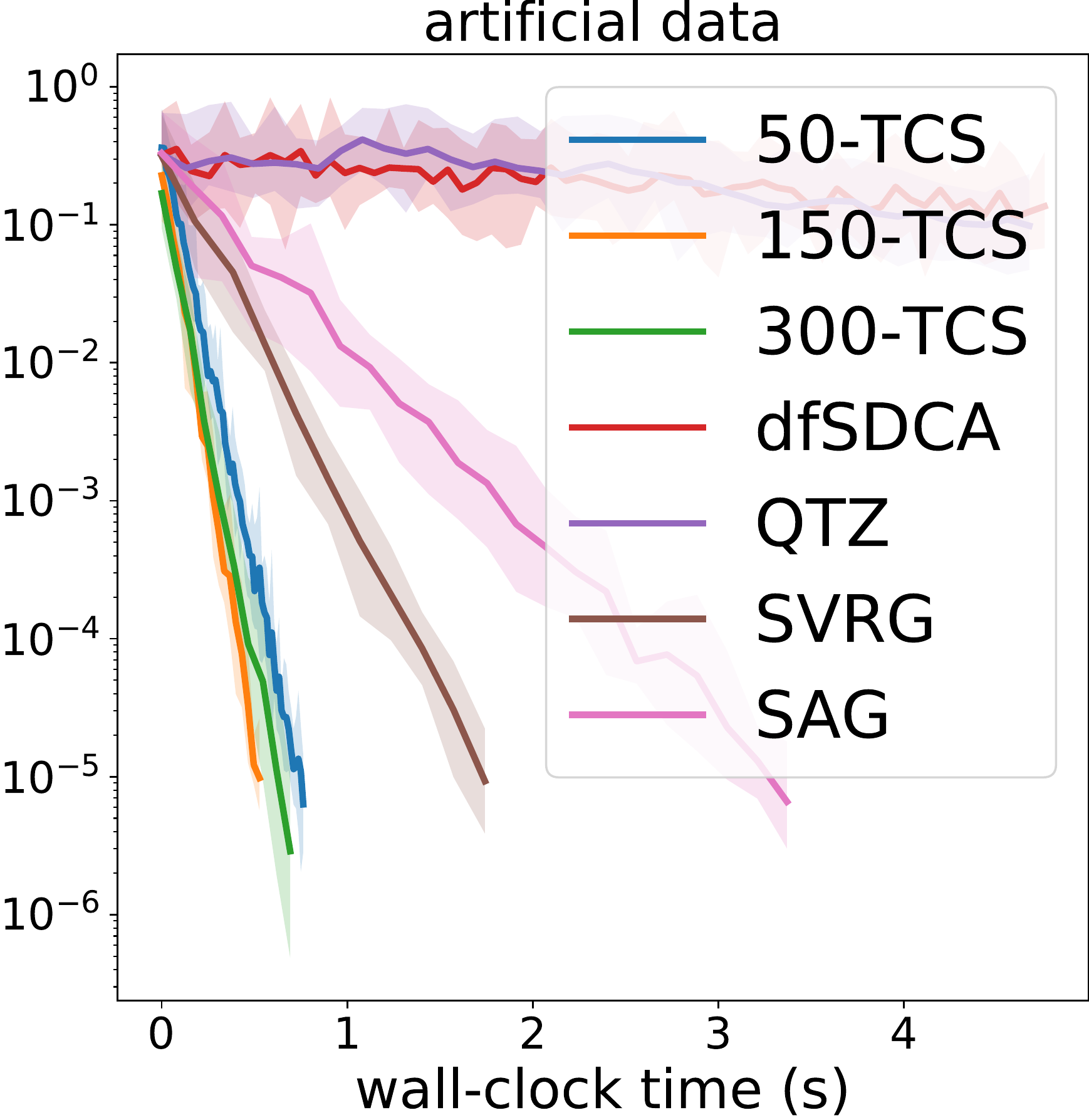} \\
\includegraphics[width=.21\linewidth]{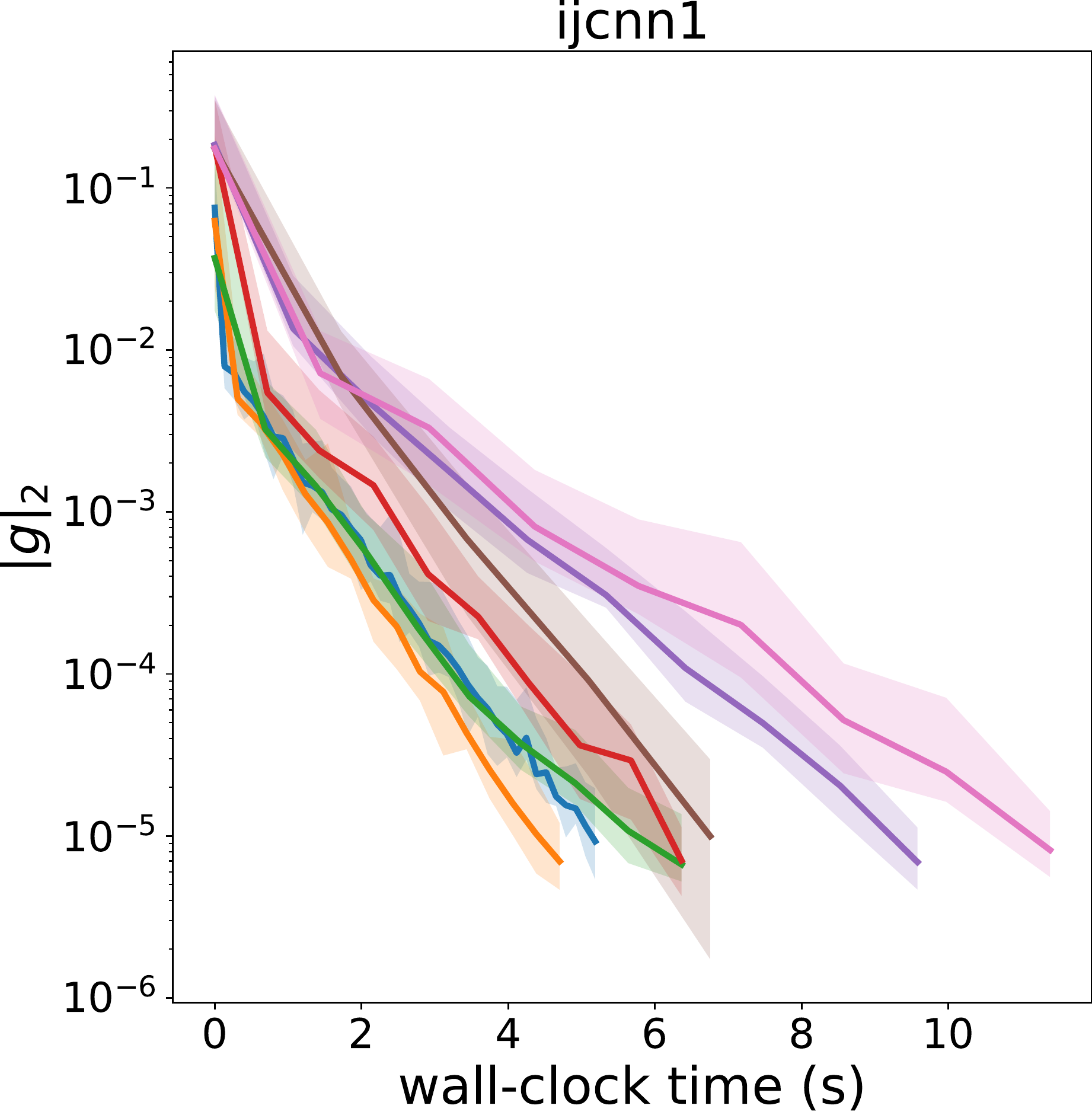}&
\includegraphics[width=.21\linewidth]{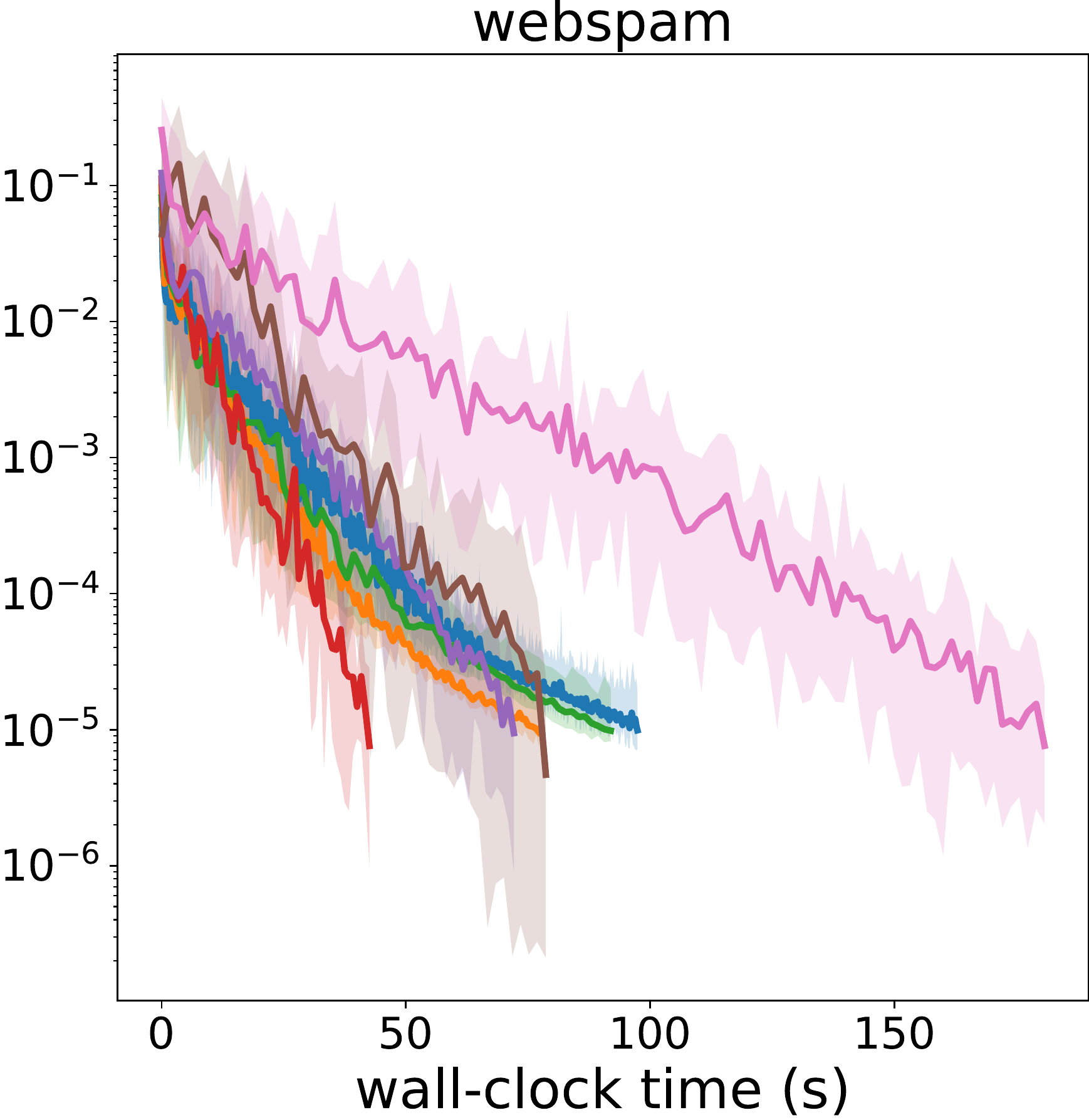}&
\includegraphics[width=.21\linewidth]{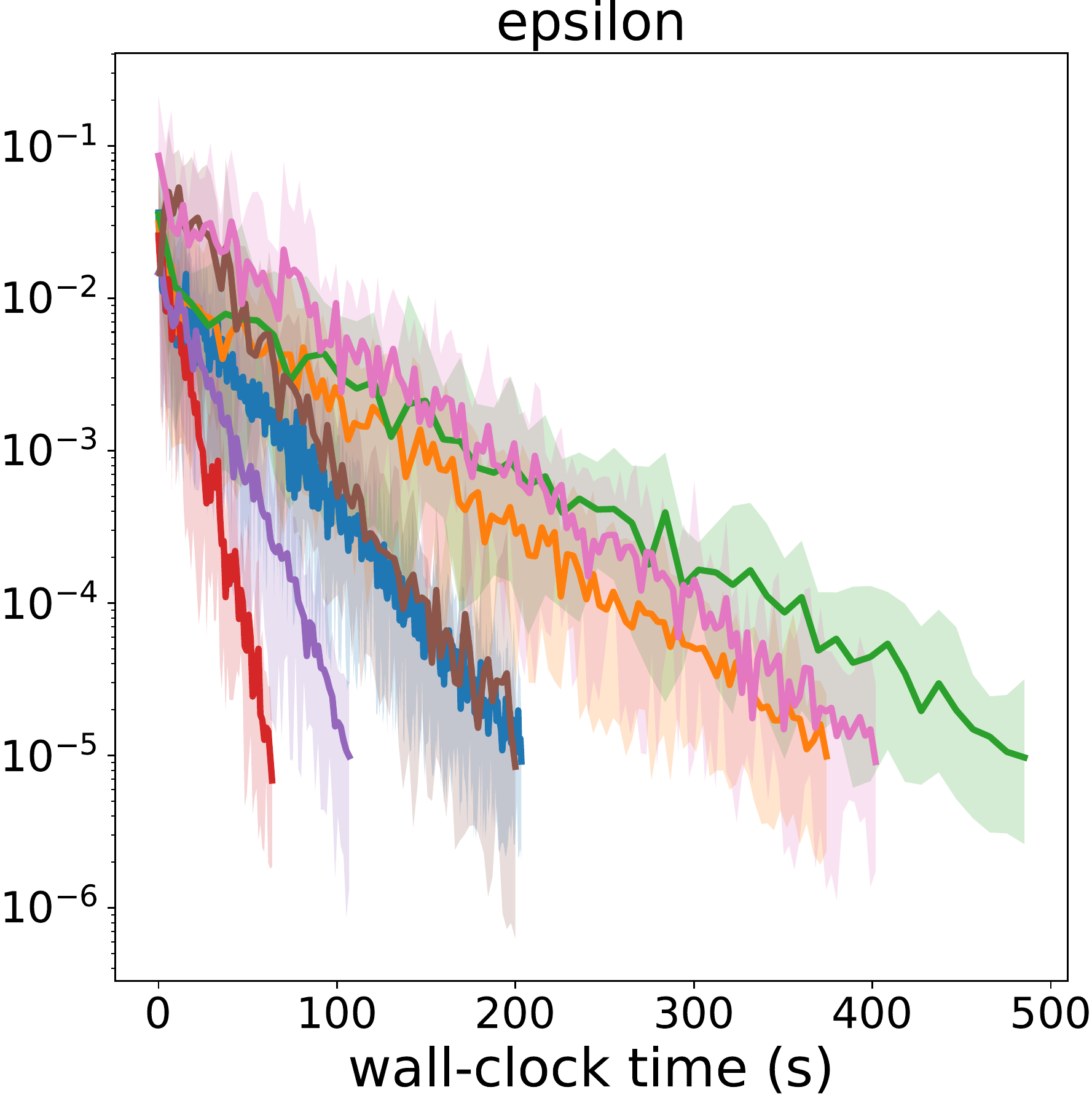}&
\includegraphics[width=.21\linewidth]{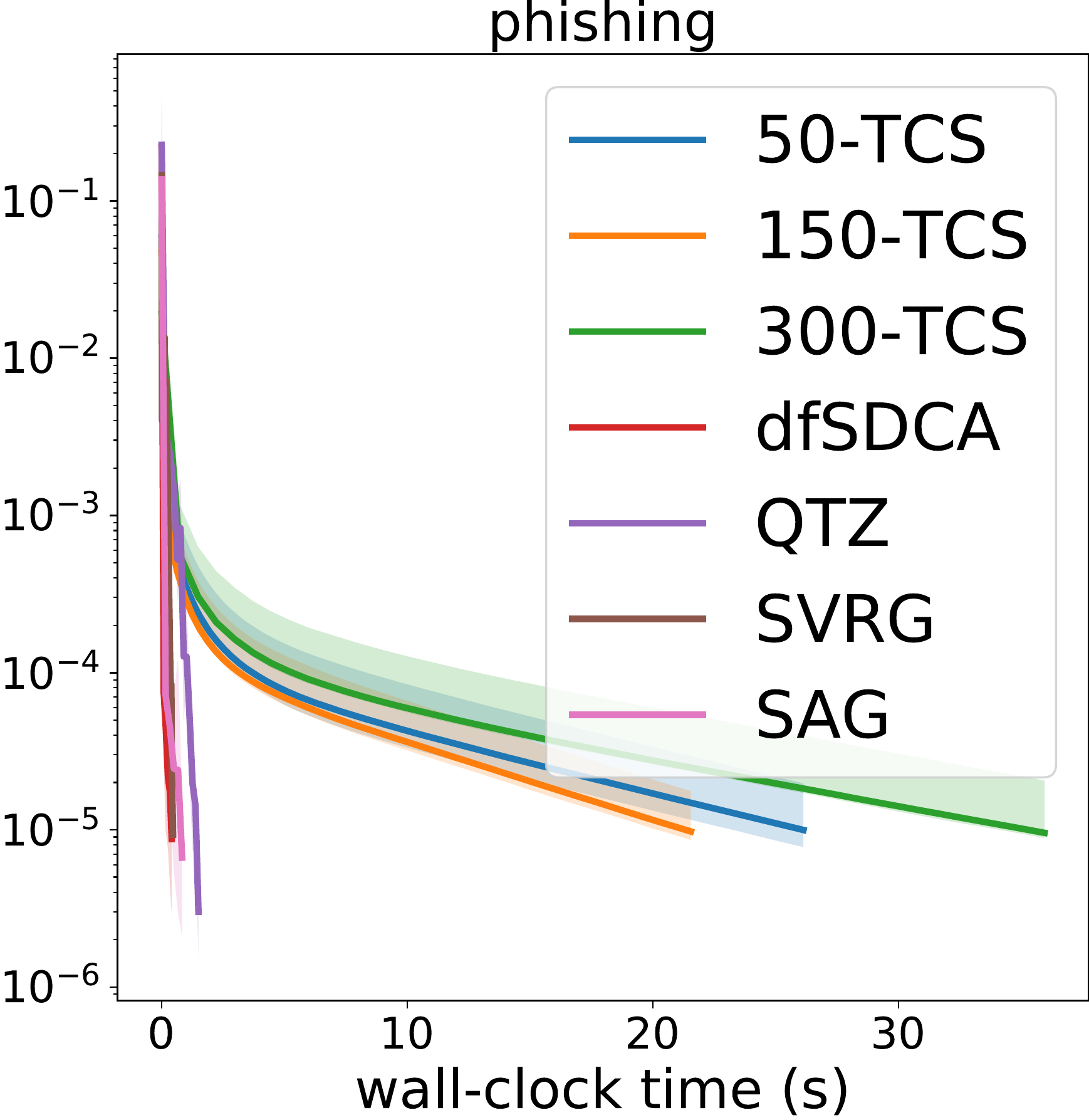}
\end{tabular}
\caption{Experiments for TCS method applied for generalized linear model.}
\label{figure}
\end{figure}

We compare the TCS method with SAG~\cite{SAG}, SVRG~\cite{Johnson2013}, dfSDCA~\cite{Shalev-Shwartz2015} and Quartz~\cite{Qu2015b}. All experiments were initialized at $w^0 = 0 \in \R^d$ (and/or $\alpha^0 = 0 \in \R^n$ for TCS/dfSDCA methods)
 and were performed in Python 3.7.3 on a laptop with an Intel Core i9-9980HK CPU and 32 Gigabyte of DDR4 RAM running OSX 10.14.5. 
For all methods, we used the stepsize that was shown to work well in practice. For instance, the common rule of thumb for SAG and SVRG is to use a stepsize $\frac{1}{L}$, where $L$ is the smoothness constant. This rule of thumb stepsize is not supported by theory. Indeed for SAG, the theoretical stepsize is $\frac{1}{16L}$ and it should be even smaller for SVRG depending on the condition number. For dfSDCA and Quartz’s, we used the stepsize suggested in the experiments in~\cite{Shalev-Shwartz2015} and~\cite{Qu2015b} respectively.
For TCS, we used two types of stepsize, related to the C.N.\ of the model. If the condition number is big (Fig.~\ref{figure} top row), we used $\gamma = 1$ except for a9a with $\gamma = 1.5$. If the condition number is small (Fig.~\ref{figure} bottom row), we used $\gamma = 1.8$.
We also set the Bernoulli parameter $b$ (probability of the coin toss) depending on the size of the dataset (see Table~\ref{tab:parameters} in Sec.~\ref{sec:detailexp}), and $\tau_d = d$. We tested three different sketch sizes $\tau_n = 50, 150, 300$. More details of the parameter settings are presented in Sec.~\ref{sec:detailexp}.

We used $\lambda = \frac{1}{n}$ regularization parameter, evaluated each method $10$ times and stopped once the gradient norm\footnote{\xline{We evaluated the true gradient norm every $1000$ iterations. We also paused the timing when computing the performance evaluation of the gradient norm.}} was below $10^{-5}$ or some maximum time had been reached. In Fig.~\ref{figure}, we plotted the central tendency as a solid line and all other executions as a shaded region for the wall-clock time vs gradient norm.

From Fig.~\ref{figure}, TCS outperforms all other methods on ill-conditioned problems (Fig.~\ref{figure} top row), but not always the case on well-conditioned problems (Fig.~\ref{figure} bottom row). This is because in ill-conditioned problems, the curvature of the optimization function is not uniform over directions and varies in the input space. Second-order methods effectively exploit information of the local curvature to converge faster than first-order methods.
To further illustrate the performance of TCS on ill-conditioned problems, we compared the performance of TCS on the artificial dataset
in the top right of Fig.~\ref{figure}.
Note as well that for reaching an approximate solution at early stage (i.e.\ $tol = 10^{-3}, 10^{-4}$), TCS is very competitive on all problems. TCS also has the smallest variance compared to the first-order methods based on eye-balling the shaded error bars in Fig.~\ref{figure}, especially compared to SVRG. Among the three tested sketch sizes, $150$ performed the best except on the \textit{epsilon} dataset.


\section{Conclusion and future work}

We introduced the \SNR method, for which we provided strong convergence guarantees. We also developed several promising applications of \SNR to show that \SNR is very flexible and tested one of these specialized variants for training GLMs. \SNR is flexible by the fact that its primitive goal is to solve efficiently nonlinear equations. Since there are many ways to re-write an optimization problem as nonlinear equations, each re-write leads to a distinct method, thus leads to a specific implementation in practice (e.g.\ \SNM, TCS methods) when using \SNR. Besides, the convergence theories presented in Sec.~\ref{sec:convergence} guarantee a large variety of choices for the sketch. This flexibility allows us to discover many applications of \SNR and their induced consequences, especially providing new global convergence theories. As such, we believe that \SNR and its global convergence theory will open the way to designing and analyzing a host of new stochastic second order methods. Further venues of investigation include exploring the use of adaptive norms for projections and leveraging efficient sketches (e.g.\ the fast Johnson-Lindenstrauss sketch~\cite{Pilanci2015a}, sketches with determinantal sampling~\cite{mutny2020convergence}) to design even faster variants of \SNR or cover other stochastic second-order methods. Since \SNR can be seen as \SGD, it might be possible to design and develop efficient accelerated \SNR or \SNR with momentum methods. On the experimental side, it would be interesting to apply our method to the training of deep neural networks.


\appendix

\section{Other viewpoints of \SNR} \label{sec:viewpoints}

Beside the connection between \SNR and \SGD, in the next section we reformulate \SNR as a stochastic Gauss-Newton (\GN) method and a stochastic fixed point method in the subsequent Appendix~\ref{sec:fixedpoint}.

\subsection{Stochastic Gauss-Newton method} \label{sec:SGN}

The \GN method is a method for solving nonlinear least-squares problems such as
\begin{equation}\label{eq:leastsquares}
\min_{x \in \R^p}\norm{F(x)}^2_{\mG}, 
\end{equation}
where $\mG$ is a symmetric positive-definite matrix. Like the Newton-Raphson method, at each step of the \GN method, the function $F(x)$ is replaced by its linearization in~\eqref{eq:leastsquares} and then solved to give the next iterate. That is
\begin{align}\label{eq:originalGN}
x^{k+1} &\in  \argmin_{x \in \R^p} \norm{DF(x^k)^\top (x-x^{k}) + \gamma F(x^k)}_{\mG}^2,
\end{align}
where $x^{k+1}$ is the least-norm solution to the above.

Now consider the \GN method where the matrix that defines the norm in~\eqref{eq:originalGN} changes at each iteration as is given by  $\mG \equiv \mG^k \eqdef \E{\mH_{\mS}(x^k)}$ and let $d \eqdef x - x^k$. Since $\mG^k$ is an expected matrix, we can write
\begin{eqnarray*}
\norm{DF(x^k)^\top d + \gamma F(x^k)}_{\E{\mH_{\mS}(x^k)}}^2 &=& \E{\norm{DF(x^k)^\top d + \gamma F(x^k)}_{\mH_{\mS}(x^k)}^2}.
\end{eqnarray*}
This suggests a stochastic variant of the \GN where we use the unbiased estimate $\mH_{\mS}(x^k)$ instead of $\mG^k$. This stochastic \GN method is in fact equivalent to \SNR, as we show next.

\begin{lemma} \label{lem:SGNequiv}
Let $x^0 \in \R^p$ and consider the following \emph{Stochastic Gauss-Newton} method
\begin{align}
d^{k} &\in  \underset{d \in \R^p}{\argmin} \norm{DF(x^k)^\top d + \gamma F(x^k)}_{\mH_{\mS_k}(x^k)}^2 \nonumber\\
x^{k+1} & = x^k +d^k\label{eq:SGN}
\end{align}
where $\mS_k$ is sampled from $\cD_{x^k}$ at $k$th iteration and $d^{k}$ is the least-norm solution. If Assumption~\ref{ass:existence} holds, then the iterates~\eqref{eq:SGN} are equal to the iterates of \SNR~\eqref{eq:update}.
\end{lemma}

\begin{proof}
Differentiating~\eqref{eq:SGN} in $d$, we find that $d^k$ is a solution to
\[DF(x^k)\mH_{\mS_k}(x^k)DF(x^k)^\top d^k =  - \gamma DF(x^k)\mH_{\mS_k}(x^k)F(x^k).\]
Let $\mA \eqdef DF(x^k)\mH_{\mS_k}(x^k)DF(x^k)^\top.$ Taking the least-norm solution to the above gives
\begin{align*}
d^k \; &= \; -\gamma \mA^\dagger DF(x^k)\mH_{\mS_k}(x^k)F(x^k) \; = \; - \gamma \mA^\dagger \mA v \\
\; &= \; - \gamma \mA^\dagger \mA \mA v \; = \; - \gamma \mA v \\
\; &= \; - \gamma DF(x^k)\mH_{\mS_k}(x^k)F(x^k),
\end{align*}
where on the first line, we used that Assumption~\ref{ass:existence} shows there exists $v \in \R^p$ such that $F(x^k) = DF(x^k)^\top v$. 
On the second line, we used that $\mA = \mA \mA$ which is shown in the proof of Corollary~\ref{lem:EuclidieanF}
. Then we used $\mA^\dagger \mA \mA = \mA$ which is a property of the pseudoinverse operator that holds for all symmetric matrices. Consequently $x^{k+1} = x^k +d^k$ which is exactly the update given in~\eqref{eq:update}. 
\end{proof}

Thus our sketched Newton-Raphson method can also be seen as a stochastic Gauss-Newton method. Furthermore, if $\mS =\mI$, then~\eqref{eq:SGN} is no longer stochastic and is given by
\begin{align}
d^{k} &\in  \underset{d \in \R^p}{\argmin} \norm{DF(x^k)^\top d + \gamma F(x^k)}_{(DF(x^k)^\top DF(x^k))^\dagger}^2 \nonumber\\
x^{k+1} & = x^k +d^k. \label{eq:sm8es8j4}
\end{align}
Thus as a consequence of Lemma~\ref{lem:SGNequiv}, we have that this variant~\eqref{eq:sm8es8j4} of \GN  is in fact the Newton-Raphson method.

\subsection{Stochastic fixed point method} \label{sec:fixedpoint}

In this section, we reformulate \SNR as a stochastic fixed point method. Such interpretation is inspired from \cite{Reformulation}'s stochastic fixed point viewpoint. We extend their results from the linear case to the nonlinear case.

Assume that Assumption~\ref{ass:existence} holds and re-consider the sketch-and-project viewpoint $\eqref{eq:sketch}$ in Section~\ref{sec:sketchproj}. Note the zeros of the function $F$
\[\cL \eqdef \left\{ x \mid F(x) = 0 \right\}\]
and the sketched Newton system based on $y$
\begin{eqnarray*}
\cL_{\mS,y} &\eqdef& \left\{ x \in \R^p \mid \mS^\top DF(y)^\top(x-y) = -\mS^\top F(y) \right\}
\end{eqnarray*}
with $y \in \R^p$ and $\mS \sim \cD_y$. For a closed convex set $\cY \subseteq \R^d$, let $\Pi_{\cY}$ denote the projection operator onto $\cY$. That is
\begin{eqnarray} \label{eq:projectiony}
\Pi_{\cY}(x) &\eqdef& \argmin_{y \in \R^p} \left\{\norm{y-x} : y \in \cY\right\}.
\end{eqnarray}
Then, from \eqref{eq:sketch} by plugging $\cY = \cL_{\mS,y}$ and $y=x$ into \eqref{eq:projectiony}, we have
\begin{eqnarray}
\Pi_{\cL_{\mS,x}}(x) &=& x - DF(x)\mH_{\mS}(x)F(x). \label{eq:projection}
\end{eqnarray}
Now we can introduce the fixed point equation as follows
\begin{eqnarray} \label{eq:chi}
\chi &\eqdef& \left\{ x \mid x = \EE{\mS \sim \cD_x}{\Pi_{\cL_{\mS,x}}(x)}\right\}.
\end{eqnarray}
Assumption~\ref{ass:existence} guarantees that finding fixed points of \eqref{eq:chi} is equivalent to the reformulated optimization problem \eqref{eq:sgdreform} with $y=x$, as we show next.

\begin{lemma} \label{lem:equivfix}
If Assumption~\ref{ass:existence} holds, then
\begin{eqnarray} \label{eq:chieqf}
\chi &=& \argmin_{x \in \R^p} \frac{1}{2}\norm{F(x)}_{\EE{\mS \sim \cD_x}{\mH_{\mS}(x)}}^2.
\end{eqnarray}
\end{lemma}

\begin{proof}
Let $\chi_{\mS} \eqdef \left\{ x \mid x = \Pi_{\cL_{\mS,x}}(x)\right\}$ with $\mS \sim \cD_x$. First, we show that
\begin{eqnarray} \label{eq:chifS}
\chi_{\mS} &=& \argmin_{x\in\R^d}\frac{1}{2}\norm{F(x)}_{\mH_{\mS}(x)}^2.
\end{eqnarray}
In fact,
\begin{eqnarray*}
&\quad& x \in \chi_{\mS} 
\overset{\eqref{eq:projection}}{\Longleftrightarrow} DF(x)\mH_{\mS}(x)F(x) = 0 \nonumber \\
&\overset{\mbox{Assumption~\ref{ass:existence}}}{\Longleftrightarrow}& \exists v \in \R^d \mbox{ s.t. } F(x) = DF(x)^\top v \mbox{ and } DF(x)\mH_{\mS}(x)DF(x)^\top v = 0 \nonumber \\
&\Longleftrightarrow& \mH_{\mS}(x)F(x) = 0 \quad \quad (\mbox{as } DF(x)\mH_{\mS}(x)DF(x)^\top \succeq 0) \nonumber \\
&\Longleftrightarrow& \frac{1}{2}\norm{F(x)}_{\mH_{\mS}(x)}^2 = 0.
\end{eqnarray*}
So we induce \eqref{eq:chifS}. Finally \eqref{eq:chieqf} follows by taking expectations with respect to $\mS$ in \eqref{eq:chifS}.
\end{proof}

To solve the fixed point equation \eqref{eq:chi}, the natural choice of method is the stochastic fixed point method with relaxation. That is, we pick a relaxation parameter $\gamma > 0$, and consider the following equivalent fixed point problem
\begin{eqnarray*}
x &=& \EE{\mS \sim \cD_x}{\gamma \Pi_{\cL_{\mS,x}}(x) + (1-\gamma)x}.
\end{eqnarray*}
Using relaxation is to improve the contraction properties of the map. Then at $k$th iteration,
\begin{eqnarray}
x^{k+1} &=& \gamma \Pi_{\cL_{\mS,x^k}}(x^k) + (1-\gamma)x^k \label{eq:sfp},
\end{eqnarray}
where $\mS \sim \cD_{x^k}$. Consequently, it is straight forward to verify that \eqref{eq:sfp} is exactly the update given in \eqref{eq:update}.

\section{Sufficient conditions for reformulation assumption~\texorpdfstring{\eqref{ass:ker}}{}} \label{sec:suffcondassker}

To give sufficient conditions for~\eqref{ass:ker} to hold, we need to study the spectra of $\E{\mH_{\mS}(x)}$ .
The expected matrix $\E{\mH_{\mS}(x)}$ has made an appearance in several references~\cite{RSN_nips,mutny2020convergence,derezinski2020precise-expressions} in different contexts and with different sketches. We build upon some of these past results and adapt them to our setting.

First note that~\eqref{ass:ker} holds if $\E{\mH_{\mS}(x)}$ is invertible. 
The invertibility of $\E{\mH_{\mS}(x)}$ was already studied in detail in the linear setting in Thm.~3 in~\cite{Gower2015c} when $\mS$ is sampled from a discrete distribution.
Here we can state a sufficient condition of~\eqref{ass:ker} for sketching matrices that have a continuous distribution.

\begin{lemma}\label{lem:inv}
For every $x \in \R^p$, 
if $\EE{\mS\sim\cD_x}{\mS\mS^\top}$ and $DF(x)^\top DF(x)$  are invertible, then $\EE{\mS\sim\cD_x}{\mH_{\mS}(x)}$
is invertible.
\end{lemma}

\begin{proof}
Let $x \in \R^p$ and $\mS \sim \cD_x$. Let $\mG = DF(x)^\top DF(x)$ which is thus symmetric positive definite and $\mW = \mS^\top$. In this case, since $\mG$ is invertible we have that ${\bf Ker}(\mG) = \{0\} \subset {\bf Ker}(\mW)$ verified, by Lemma~\ref{lem:ker}, we have that
\begin{equation} \label{eq:ker2}
{\bf Ker}\left(\left(\mS^\top DF(x)^\top DF(x)\mS\right)^\dagger\right)  =  {\bf Ker}\left(\mS^\top DF(x)^\top DF(x)\mS\right)  =  {\bf Ker}(\mS),
\end{equation}
Consequently, using Lemma \ref{lem:ker} again with $\mG = \left(\mS^\top DF(x)^\top DF(x)\mS\right)^\dagger$, $\mW = \mS$ and ${\bf Ker}(\mG) \subset {\bf Ker}(\mW)$ given by~\eqref{eq:ker2}, we have that
\begin{equation}\label{eq:ker}
{\bf Ker}(\mH_{\mS}(x)) = {\bf Ker}\left(\mS\left(\mS^\top DF(x)^\top DF(x)\mS\right)^\dagger\mS^\top\right)  =  {\bf Ker}(\mS^\top)  =  {\bf Ker}(\mS\mS^\top).
\end{equation}
\ruiline{Following the same steps in the proof of Lemma~\ref{lem:reformulationnew} right after~\eqref{eq:ker3new}, we obtain}
\begin{align*}
{\bf Ker}(\E{\mH_{\mS}(x)}) = \bigcap_{\mS\sim\cD_x} {\bf Ker}(\mH_{\mS}(x)) 
\overset{\eqref{eq:ker}}{=} \bigcap_{\mS\sim\cD_x}{\bf Ker}(\mS\mS^\top) = {\bf Ker}(\E{\mS\mS^\top}) = \{0\},
\end{align*}
where the last equality follows as $\E{\mS\mS^\top}$ is invertible,
which concludes the proof.
\end{proof}

The invertibility of $\E{\mS\mS^\top}$ states that the sketching matrices need to ``span every dimension of the space'' in expectation. 
\ruilinex{
This is the case for Gaussian and subsampling sketches which are shown in Lemma~\ref{lem:SSTinv}. This is also the case for our applications \SNM and TCS which are shown in the proofs of Lemma~\ref{lem:adaptedsketchSNM} and Lemma~\ref{lem:adaptTCS}, respectively.}

\ruiline{As for the invertibility of $DF(x)^\top DF(x) \in \R^{m \times m}$, this imposes that $DF(x)$ has full-column rank for all $x \in \R^p$, thus $m \leq p$. This excludes the regime of solving $F(x) = 0$ with $m > p$.} \ruilinex{However, our applications \SNM and TCS also satisfy this condition which are again shown in the proofs of Lemma~\ref{lem:DFinv} and Lemma~\ref{lem:adaptTCS}, respectively.}


\ruilinex{
Consequently, by Lemma~\ref{lem:inv}, we 
have that \SNM and TCS satisfy~\eqref{ass:ker}.}




\section{The monotone convergence theory of \NR with stepsize \texorpdfstring{$\gamma < 1$}{g < 1}}
\label{sec:monotoneconv1}

The MCT of the \NR in both~\cite{Ortega:2000} and~\cite{DeuflhardNewton:2011} need to have the stepsize $\gamma = 1$. If $\gamma < 1$ which is the case in our convergence Theorem~\ref{theo:convex} and Corollary~\ref{lem:EuclidieanF}, the iterates $\{x^k\}_{k \geq 1}$ under the set of assumptions (I) proposed in Theorem~\ref{theo:strongerthanmono} are no longer guaranteed to be component wise monotonically decreasing. Here we investigate alternatives. In particular, we consider the case in $1$-dimension for function $F = \phi: \R \rightarrow \R$.

\begin{lemma}\label{lem:strongerthanmono1D}
Let $x^k$ be the iterate of the \NR method with stepsize $\gamma < 1$ for solving $\phi(x) = 0$, that is
\begin{align}
x^{k+1} &= x^k - \gamma\frac{\phi(x^k)}{\phi'(x^k)}. \label{eq:newton_gamma}
\end{align}
If $\phi$ satisfies the set of assumptions (I) proposed in Theorem \ref{theo:strongerthanmono}, then
\begin{enumerate}
\item[(a)] The iterates of the ordinary \NR method \eqref{eq:newton1} are necessarily monotonically decreasing.
\item[(b)] The iterates of the \NR method \eqref{eq:newton_gamma} with $\gamma < 1$ are not necessarily monotonically decreasing.
\item[(c)] Assumption~\eqref{eq:KermSI} holds; for $\frac{1}{2} \leq \gamma < 1$, there exists $k' \geq 0$ such that for all $k \neq k'$, there exists a unique $x^*$ that satisfies Assumption~\ref{ass:zero} and the iterates $x^k$ and the optimum $x^*$ satisfy~\eqref{eq:nxcvx}.
\item[(d)] The iterates $x^k$ following the \NR method \eqref{eq:newton_gamma} with $\frac{1}{2} \leq \gamma < 1$ converge sublinearly to a zero of $\phi$.
\end{enumerate}
\end{lemma}

\paragraph*{Remark of (a)} Even though this result is known and generalized in $d$-dimension in \cite{Ortega:2000} and \cite{DeuflhardNewton:2011}, we stress it here to highlight the impact of the stepsize $\gamma$ in the \NR method and leverage the analysis of (a) in the special $1$-dimensional case to prove (b).

\begin{proof}
If $\phi$ satisfies (I), then $\phi$ is convex and $\phi'^{-1} > 0$, which implies $\phi'' \geq 0$ and $\phi' > 0$. From $\phi' > 0$, we obtain that $\phi$ is strictly increasing. Besides, from (I), $\exists \ x, y \in \R$ such that $\phi(x) \leq 0 \leq \phi(y)$. This with the strictly increase of $\phi$ induces that $\exists! \ x^*$ such that $\phi(x^*) = 0$, i.e. $x^*$ satisfies Assumption~\ref{ass:zero}. So $\forall x < x^*$, $\phi(x) < 0$ and $\forall x > x^*$, $\phi(x) > 0$.
Now consider the following two functions
\begin{equation*}
u(x) \; \eqdef \; x - \frac{\phi(x)}{\phi'(x)} \quad \quad \mbox{and } \quad \quad U(x) \; \eqdef \; x - \gamma\frac{\phi(x)}{\phi'(x)}
\end{equation*}
which are exactly the updates of the ordinary \NR \eqref{eq:newton1} and the \NR \eqref{eq:newton_gamma} with a stepsize $\gamma \in (0,1)$, respectively. We first analyze the behaviour of the function $u$ and show (a), which can be formulated as $x^* \leq x^{k+1} \leq x^k$ for all $k \geq 1$. The derivative of $u$ is
\begin{equation*}
u'(x) = \frac{\phi(x)\phi''(x)}{\phi'(x)^2}.
\end{equation*}
By the sign of functions $\phi$, $\phi'$ and $\phi''$, we know that if $x > x^*$, then $u'(x) \geq 0$ and if $x < x^*$, then $u'(x) \leq 0$. This implies that the function $u$ is increasing in $[x^*, +\infty[$ and decreasing in $]-\infty, x^*]$. Overall, we have
\begin{subnumcases}{}
\min_{x \in \R} u(x) = u(x^*) \overset{\phi(x^*)=0}{=} x^*, \label{eq:cases_u1} \\
u(x) < x \mbox{ and $u$ increasing}, \quad \mbox{when } x > x^*, \label{eq:cases_u2} \\
u(x) > x \mbox{ and $u$ decreasing}, \quad \mbox{when } x < x^*. \label{eq:cases_u3}
\end{subnumcases}
Consequently, $x^* \leq x^{k+1}$ is obtained by $x^* \overset{\eqref{eq:cases_u1}}{=} \min u(x) \leq u(x^k) = x^{k+1}$. As for the inequality $x^{k+1} \leq x^k$, $x^* \overset{\eqref{eq:cases_u1}}{=} \min u(x) \leq u(x^{k-1}) = x^k$ for $k \geq 1$ and $x^{k+1} = u(x^k) \overset{\eqref{eq:cases_u2}}{\leq} x^k$ as $x^k \geq x^*$.

To show (b), we analyze the behavior of the function $U$. Consider its derivative
\[U'(x) = (1 - \gamma) + \gamma\frac{\phi(x)\phi''(x)}{\phi'(x)^2}.\]
By the sign of functions $\phi$, $\phi'$ and $\phi''$, if $x > x^*$, $U'(x) > 0$. However, $U'(x^*) \overset{\phi(x^*)=0}{=} 1 - \gamma > 0$, which implies $\min U(x) < x^*$. Here we include the case where $\min U(x) = -\infty$. Also by the sign of functions $\phi$ and $\phi'$ and $\gamma < 1$, when $x > x^*$, we have $u(x) < U(x)$ and $u(x) > U(x)$ for $x < x^*$. In summary, we have
\begin{subnumcases}{}
\min_{x \in \R} U(x) < U(x^*) \overset{\phi(x^*)=0}{=} u(x^*) \overset{\eqref{eq:cases_u1}}{=} x^*, \label{eq:case_U1} \\
u(x) < U(x) < x \mbox{ and $U$ increasing}, \quad \mbox{when } x > x^*, \label{eq:case_U2} \\
u(x) > U(x) > x \mbox{ when } x < x^*. \label{eq:case_U3}
\end{subnumcases}
In \NR with stepsize $\gamma < 1$, consider $x^0 \in \R$. We discuss different cases based on the comparison between $x^0$ and $x^*$.

If $x^0 \geq x^*$ named as case {\bf (i)}, by induction, we get $x^* \overset{\eqref{eq:case_U1}}{=} U(x^*) \overset{\eqref{eq:case_U2}}{\leq} = U(x^k) = x^{k+1} \overset{\eqref{eq:case_U2}}{\leq} x^k$ for all $k \geq 0$. In this case, the iterates decrease monotonically.

If $x^0 < x^*$, there are two cases, named as case {\bf (ii)}, for all $k \in \N$, $U(x^k) \leq x^*$, and case {\bf (iii)}, $\exists k' \in \N$, $U(x^{k'}) > x^*$.

If {\bf (ii)} holds, we have that the iterates increase monotonically. Indeed, by {\bf (ii)} and by induction, we get $x^k \overset{\eqref{eq:case_U3}}{\leq} U(x^k) = x^{k+1} \overset{\mbox{{\bf (ii)}}}{\leq} x^*$ for all $k \in \N$.

Otherwise, we are in case {\bf (iii)}. Let $k'$ be the smallest index that $U(x^{k'}) > x^*$. Then we conclude that the iterates $\{x^k\}_{k \geq 0}$ increase monotonically when $k \leq k'$ and $\{x^k\}_{k \geq k' + 1}$ decrease monotonically. In fact, by the definition of $k'$, we know that for $k \in [\![ 0, k' - 1 ]\!]$, $U(k) \leq x^*$. Then by induction as in case {\bf (ii)} but for $k \leq k'$, we get $\{x^k\}_{k \geq 0}$ increase monotonically when $k \leq k'$. When $k \geq k'+1$, by induction as in case {\bf (i)} but for $U(x^{k'}) = x^{k'+1} > x^*$, we get $\{x^k\}_{k \geq k' + 1}$ decrease monotonically. We thus observe (b).

Statement (c) follows from the proof of Theorem \ref{theo:strongerthanmono} in $1$-dimension in taking account the stepsize $\gamma < 1$. Then~\eqref{eq:KermSI} holds and \eqref{eq:nxcvx} becomes
\begin{eqnarray} \label{eq:sni848sh2}
&\quad& f_k(x^k) + \dotprod{\nabla f_k(x^k), x^* - x^k} \\
&=& \underbrace{\frac{1-2\gamma}{2\gamma^2}(x^{k+1}-x^k)^2}_{\leq 0 \mbox{ as } \frac{1}{2} \leq \gamma \mbox{ in (c)}} + \frac{1}{\gamma}\underbrace{(x^k - x^{k+1})(x^* - x^{k+1})}_{\eqdef (*)} \nonumber
\end{eqnarray}
in considering
\begin{eqnarray*}
f_k(x^k) &=& \frac{1}{2}\left(\frac{\phi(x^k)}{\phi'(x^k)}\right)^2 \overset{\eqref{eq:newton_gamma}}{=}\frac{1}{2\gamma^2}(x^{k+1}-x^k)^2, \\  
\nabla f_k(x^k) &=& \frac{\phi(x^k)}{\phi'(x^k)} \overset{\eqref{eq:newton_gamma}}{=} \frac{1}{\gamma}\left(x^k-x^{k+1}\right),
\end{eqnarray*}
with $\phi = F$ and $\phi' = DF$. To get~\eqref{eq:nxcvx} hold, from \eqref{eq:sni848sh2}, it suffices to prove $(*) \leq 0$.

By the analysis of (b), we know: in case {\bf (i)}, $(*) \leq 0$ for all $k \geq 0$ as $x^* \leq x^{k+1} \leq x^k$; in case {\bf (ii)}, $(*) \leq 0$ for all $k \geq 0$ as $x^* \geq x^{k+1} \geq x^k$; finally in case {\bf (iii)}, for $k \neq k'$, $(*) \leq 0$ as $x^* \geq x^{k+1} \geq x^k$ for $k \leq k' -1$ and $x^* \leq x^{k+1} \leq x^k$ for $k \geq k'+1$. So in all cases, $(*) \leq 0$ for all $k$ or for $k \neq k'$. We thus obtain (c).

It remains to show (d), which is simply obtained by (c) and Corollary \ref{lem:EuclidieanF}, as~\eqref{eq:nxcvx} holds for all iterates $x^k$ except for just one iterate $x^{k'}$ potentially.
\end{proof}

The monotone convergence theory is based on assumptions (I) with stepsize $\gamma = 1$. Under the same assumptions with $\gamma < 1$, such theory may not hold. Indeed, following the analysis in Lemma~\ref{lem:strongerthanmono1D} in $1$-dimension case, by \eqref{eq:case_U3} we do not have the monotone property for the function $U$ when $x < x^*$. That is the reason why (b) happens but not (a) in Lemma~\ref{lem:strongerthanmono1D}. In $d$-dimension case, without such monotone property for the function $U$, $\{x^k\}$ is not guaranteed to be monotone, which is the main clue in their theory's proof. However, with stepsize $\gamma < 1$, assumptions (I) can still imply our Assumptions~\ref{ass:zero},~\eqref{eq:KermSI} and~\eqref{eq:nxcvx} under constraint $\frac{1}{2} \leq \gamma < 1$ in $1$-dimension case. In addition, though our theory does not either require any constraint for stepsize $\gamma < 1$ or guarantee that the \NR method is monotonic in terms of the iterates component wisely, we still guarantee the sublinear global convergence. We thus conclude that Assumptions~\ref{ass:zero},~\eqref{eq:KermSI} and~\eqref{eq:nxcvx} are strictly weaker than the assumptions used in the monotone convergence theory in \cite{Ortega:2000} and \cite{DeuflhardNewton:2011}, albeit for different step sizes.

\section{Stochastic Newton method with relaxation} \label{sec:SNMwRelax}

Consider the function $P(\cdot)$ defined in~\eqref{eq:fSNM}. By the analysis of Lemma \ref{lem:SNMdetail}, we can even develop a variant of \SNM in the case stepsize $\gamma < 1$ and we call the method \emph{Stochastic Newton method with relaxation}. The updates are the following
\begin{align}
w^{k+1} \; &= \; \gamma \bigg(\frac{1}{n}\sum_{i=1}^n \nabla^2 \phi_i(\alpha^k_i)\bigg)^{-1} \bigg(\frac{1}{n}\sum_{i=1}^n \nabla^2 \phi_i(\alpha^k_i) \alpha^k_i - \frac{1}{n}\sum_{i=1}^n \nabla \phi_i(\alpha^k_i)\bigg) + (1-\gamma) w^k, \label{eq:xgamma} \\
\alpha^{k+1}_i \; &= \; \label{eq:wgamma}
\begin{cases}
w^{k+1} - (1-\gamma)(w^k - \alpha^k_i) & \quad \mbox{if } i \in B_n \\
\alpha^k_i & \quad \mbox{if } i \notin B_n
\end{cases},
\end{align}
In the rest of Section~\ref{sec:SNMwRelax}, we use the shorthand $x \eqdef (w; \alpha_1; \cdots; \alpha_n) \in \R^{(n+1)d}$ and $x^k \eqdef (w^k; \alpha^k_1; \cdots; \alpha^k_n)$ the iterates of \SNR in Lemma~\ref{lem:SNMdetail} with stepsize $\gamma < 1$ at the $k$th iteration.

\begin{lemma} \label{lem:SNMgamma}
At each iteration $k$, the updates of \SNR $x^k$ are equal to the updates~\eqref{eq:xgamma},~\eqref{eq:wgamma} of \SNM with relaxation.
\end{lemma}

\begin{proof}
Following the proof of Lemma~\ref{lem:SNMdetail} and taking account the stepsize $\gamma$, by~\eqref{eq:newtonprojF} and~\eqref{eq:sketch}, the updates of \SNR $x^{k+1}$ at $(k+1)$th iteration are given by
\begin{align}
x^{k+1} = \argmin \norm{w - w^k}^2 + \sum_{i=1}^n\norm{\alpha_i - \alpha^k_i}^2 \ \mbox{s.t.} \ \mS^\top DF(x^k)^\top (x - x^k) = - \gamma \mS^\top F(x^k), \label{eq:newtonsketchprojFgamma}
\end{align}
where the sketching matrix $\mS\sim\cD_{x^k}$ is defined in~\eqref{eq:SNMsketch}. Similar to \eqref{eq:projxwiNR}, \eqref{eq:newtonsketchprojFgamma} can be re-written as
\begin{align}
&x^{k+1} = \argmin \norm{w - w^k}^2 + \sum_{i=1}^n\norm{\alpha_i - \alpha^k_i}^2  \label{eq:newtonsketchprojFgamma2} \\
& \ \ \ \ \mbox{s.t.} \ \frac{1}{n}\sum_{i=1}^n \nabla^2 \phi_i(\alpha^k_i)(w - w^k) = - \gamma \bigg(\frac{1}{n}\sum_{i=1}^n \nabla \phi_i(\alpha^k_i) + \frac{1}{n}\sum_{i=1}^n \nabla^2 \phi_i(\alpha^k_i)(w^k - \alpha^k_i)\bigg), \nonumber\\
&\phantom{ \ \ \ \ \mbox{s.t.}} \ w - \alpha_i = (1-\gamma)(w^k - \alpha^k_i), \quad \mbox{ for } i \in B_n. \nonumber
\end{align}
Similarly, note that if $i \not\in B_n$, then $\alpha^{k+1}_i = \alpha^k_i$, since there is no constraint on the variable $\alpha_i$ in this case. Then by the invertibility of $\frac{1}{n}\sum_{i=1}^n \nabla^2 \phi_i(\alpha^k_i)$, we have the unique solution of \eqref{eq:newtonsketchprojFgamma2}, which is
\begin{align*}
w^{k+1} \; &= \; \gamma \bigg(\frac{1}{n}\sum_{i=1}^n \nabla^2 \phi_i(\alpha^k_i)\bigg)^{-1} \bigg(\frac{1}{n}\sum_{i=1}^n \nabla^2 \phi_i(\alpha^k_i) \alpha^k_i - \frac{1}{n}\sum_{i=1}^n \nabla \phi_i(\alpha^k_i)\bigg) + (1-\gamma) w^k, \\
\alpha^{k+1}_i \; &= \; w^{k+1} - (1-\gamma)(w^k - \alpha^k_i), \quad \mbox{ for } i \in B_n.
\end{align*}
Overall, we have
\begin{align*}
w^{k+1} \; &= \; \gamma \bigg(\frac{1}{n}\sum_{i=1}^n \nabla^2 \phi_i(\alpha^k_i)\bigg)^{-1} \bigg(\frac{1}{n}\sum_{i=1}^n \nabla^2 \phi_i(\alpha^k_i) \alpha^k_i - \frac{1}{n}\sum_{i=1}^n \nabla \phi_i(\alpha^k_i)\bigg) + (1-\gamma) w^k, \\
\alpha^{k+1}_i \; &= \; 
\begin{cases}
w^{k+1} - (1-\gamma)(w^k - \alpha^k_i) & \quad \mbox{if } i \in B_n \\
\alpha^k_i & \quad \mbox{if } i \notin B_n
\end{cases},
\end{align*}
which is exactly the updates \eqref{eq:xgamma} and \eqref{eq:wgamma} in \SNM with relaxation. 
\end{proof}

Notice that both the original \SNM and \SNM with relaxation have the same complexity. Consequently, Theorem \ref{theo:convex} allows us to develop the following global convergence theory of \SNM with relaxation $\gamma$.

\begin{corollary} \label{cor:SNMgamma}
Consider the iterate $x^k = \left(w^k; \alpha^k_1; \cdots; \alpha^k_n\right)$ given by~\eqref{eq:xgamma} and~\eqref{eq:wgamma}. Note $x^* \eqdef (w^*; w^*; \cdots; w^*) \in \R^{(n+1)d}$ where $w^*$ is the stationary point of $\nabla P(w)$ that satisfies
\begin{equation} \label{eq:cvxSNM}
f_{x^k}(x^*) \geq f_{x^k}(x^k) + \dotprod{\nabla f_{x^k}(x^k), x^* - x^k}, \quad \mbox{for all } k \in \N,
\end{equation}
then
\begin{equation*}
\min_{t=0, \ldots, k-1} \E{f_{x^t}(x^t)}
\quad \leq \quad \frac{1}{k}\sum_{t=0}^{k-1}\E{f_{x^t}(x^t)} \quad  \leq \quad \frac{1}{k}\frac{\norm{x^0-x^*}^2}{2\gamma\left(1 - \gamma\right)}.
\end{equation*}
\end{corollary}

\begin{proof}
Equation~\eqref{eq:cvxSNM} with the assumption that $w^*$ is the stationary point of $\nabla P(w)$ implies that Assumption~\ref{ass:zero} and~\ref{ass:convex} hold. Then we conclude the proof by Theorem~\ref{theo:convex}.
\end{proof}

\section{Extension of \SNR and Randomized Subspace Newton} \label{sec:RSN}

In the \SNR method in \eqref{eq:update}, we only consider a projection under the standard Euclidean norm. If we allow \SNR and \eqref{eq:update} for a changing norm that depends on the iterates, we find that the Randomized Subspace Newton~\cite{RSN_nips} (\RSN) method is in fact a special case of \SNR under this extension.

The changing norm projection of \SNR is that, at $k$th iteration of \SNR, instead of applying \eqref{eq:update}, we can apply the following update
\begin{eqnarray} \label{eq:updateW}
x^{k+1} &=& x^k - \gamma\mW_k^{-1} DF(x^k)\mS_k \left(\mS_k^\top  DF(x^k)^\top \mW_k^{-1} DF(x^k)\mS_k \right)^\dagger \mS_k^\top F(x^k)
\end{eqnarray}
where $\mW_k \equiv \mW(x^k)$ with $\mW(x^k)$ a certain symmetric positive-definite matrix associated with the $k$th iterate $x^k \in \R^p$.

The interpretation of using the matrix $\mW(x^k)$ is that, assuming Assumption~\ref{ass:existence} holds, then instead of considering \eqref{eq:sketch}, we apply the following updates
\begin{eqnarray}
x^{k+1} &=& \argmin_{x \in \R^p} \norm{x-x^k}^2_{\mW_k} \nonumber \\
&\quad&  \mbox{ s. t.}\quad \mS_k^\top  DF(x^k)^\top(x-x^k) \; = \; -\gamma\mS_k^\top F(x^k), \label{eq:sketchW}
\end{eqnarray}
using the projection $\|\cdot\|_{\mW_k}$ which changes at each iteration. One can verify easily that \eqref{eq:sketchW} is equivalent to \eqref{eq:updateW} under Assumption~\ref{ass:existence}, even though this assumption is not necessary and the update \eqref{eq:updateW} is still available.

Now we can show that \RSN is a special case of \SNR with a changing norm projection.
The \RSN method~\cite{RSN_nips} is a stochastic second order method that takes a Newton-type step at each iteration to solve the minimization problem
\[\min_{x \in \R^p}P(x)\]
where $P : \R^p \rightarrow \R$ is a twice differentiable and convex function. In brevity, the updates in \RSN at the $k$th iteration are given by
\begin{eqnarray} \label{eq:RSN}
x^{k+1} &=& x^k - \frac{1}{\hat{L}}\mS_k\left(\mS_k^\top\nabla^2P(x^k)\mS_k\right)^\dagger\mS_k^\top\nabla P(x^k)
\end{eqnarray}
where $\mS_k$ is sampled i.i.d from a fixed distribution $\cD$ and $\hat{L}>0$ is the \emph{relative smoothness} constant~\cite{RSN_nips}.

Since $P(x)$ is convex, it suffices to find a stationary point $x$ such that $\nabla P(x) =0.$ 
We can recover the exact same  update~\eqref{eq:RSN} by applying  \SNR to solve  $\nabla P(x) =0$  with an adaptive changing norm.
That is, let $F(x) = \nabla P(x)$ and $DF(x) = \nabla^2P(x)$. At the $k$th iteration, let $\mW_k = \nabla^2P(x^k)$. Then \eqref{eq:updateW} is exactly the \RSN update \eqref{eq:RSN} with $\gamma = \frac{1}{\hat{L}}$.

\section{Explicit formulation of the TCS method}
\label{sec:equivalencealgorithms}

Here we provide details about how TCS method presented in Section~\ref{sec:GLM} is obtained from the general \SNR method \ref{algo:SkeNeR}.

Consider the \SNR method~\eqref{eq:update} applied for the nonlinear equations $F(\alpha; w) = 0$ with $F$ defined in \eqref{eq:F}
and the Jacobian of $F(\alpha; w)$ in \eqref{eq:DFalphax}.

At $k$th iteration $(\alpha^k, w^k) \in \R^n \times \R^d$, let 
\[\colvec{\alpha^{k+1} \\ w^{k+1}} \eqdef \colvec{\alpha^k \\ w^k} + \gamma \cdot \colvec{\Delta \alpha^k \\ \Delta w^k}\] 
and $\mS_k \in \R^{(d+n)\times(\tau_d+\tau_n)}$ the random sketching matrix. By \eqref{eq:update}, we obtain the closed form update
\begin{eqnarray}
\colvec{\Delta \alpha^k \\ \Delta w^k} &=& - DF(\alpha^k; w^k)\mS_k \left( \mS_k^\top DF(\alpha^k; w^k)^\top DF(\alpha^k; w^k)\mS_k \right)^\dagger \mS_k^\top \nonumber \\ 
&\quad& \colvec{\frac{1}{\lambda n}\mA\alpha^k - w^k \\ \alpha^k + \Phi(w^k)}. \label{eq:formula}
\end{eqnarray}
As for the tossing-coin-sketch, consider a Bernoulli parameter $b$ with $b \in (0, 1)$. There is a probability $1-b$ that the random sketching matrix has the type $\mS_k =
\begin{bmatrix}
\mS_d & 0 \\
0 & 0
\end{bmatrix}$ with $\mS_d \in \R^{d \times \tau_d}$, a $(d,\tau_d)$--block sketch, and a probability $b$ that the random sketching matrix has the type $\mS_k =
\begin{bmatrix}
0 & 0 \\
0 & \mS_n
\end{bmatrix} $ with $\mS_n \in \R^{n \times \tau_n}$, a $(n,\tau_n)$--block sketch. So $\mS_d = \mI_{B_d}$ and $\mS_n = \mI_{B_n}$. 

Let $\mA_{B_d, :} \equiv \mI_{B_d}^\top \mA \in \R^{\tau_d \times n}$  denote a row subsampling of $\mA$ and $\mA_{:, B_n} \equiv \mA\mI_{B_n} \in \R^{d \times \tau_n}$  denote a column subsampling of $\mA$. Let $\nabla\Phi^k_{B_n} \equiv \nabla \Phi^k\mI_{B_n} \in \R^{d \times \tau_n} $ denote a column subsampling of $\nabla \Phi^k$ with $\nabla \Phi^k \equiv \nabla \Phi(w^k)$ and $\Phi^k \equiv \Phi(w^k).$ We also use the shorthands $v_{B_n} \equiv \mI_{B_n}^\top v \in \R^{\tau_n}$ with $v \in \R^n$ and $v_{B_d} \equiv \mI_{B_d}^\top v \in \R^{\tau_d}$ with $v \in \R^d$.

If $\mS_k =
\begin{bmatrix}
\mS_d & 0 \\
0 & 0
\end{bmatrix} $, the update \eqref{eq:formula} applied for the function \eqref{eq:F} and its Jacobian \eqref{eq:DFalphax} becomes
\begin{eqnarray}\label{eq:formulaS_d}
\colvec{\Delta \alpha^k \\ \Delta w^k} &=& -\colvec{\frac{1}{\lambda n}\mA^\top \mS_d \\ -\mS_d}\left(\mS_d^\top\left(\frac{1}{\lambda^2n^2}\mA\mA^\top + \mI_d\right)\mS_d\right)^\dagger\mS_d^\top\left(\frac{1}{\lambda n}\mA\alpha^k - w^k\right) \nonumber \\
&=& -\colvec{\frac{1}{\lambda n}\mA_{B_d, :}^\top \\ -\mI_{B_d}}\left(\frac{\mA_{B_d, :}\mA_{B_d, :}^\top}{\lambda^2n^2} + \mI_{\tau_d}\right)^\dagger\left(\frac{\mA_{B_d, :}\alpha^k}{\lambda n} - w^k_{B_d}\right).
\end{eqnarray} 

Similarly, if $\mS_k =
\begin{bmatrix}
0 & 0 \\
0 & \mS_n
\end{bmatrix} $, the update \eqref{eq:formula} becomes
\begin{eqnarray}\label{eq:formulaS_n}
\colvec{\Delta \alpha^k \\ \Delta w^k} &=& -\colvec{\mS_n \\ \nabla\Phi^k\mS_n}\left(\mS_n^\top\left([\nabla\Phi^k]^\top \nabla\Phi^k + \mI_n\right)\mS_n\right)^\dagger\mS_n^\top\left(\alpha^k + \Phi^k\right) \nonumber \\
&=& -\colvec{\mI_{B_n} \\ \nabla\Phi^k_{B_n}}\left([\nabla\Phi^k_{B_n}]^\top\nabla\Phi^k_{B_n} + \mI_{\tau_n}\right)^\dagger\left(\alpha^k_{B_n} + \Phi^k_{B_n}\right).
\end{eqnarray}
Then we update $\colvec{\alpha^{k+1} \\ w^{k+1}} = \colvec{\alpha^k \\ w^k} + \gamma \cdot \colvec{\Delta \alpha^k \\ \Delta w^k}$.

See Algorithm \ref{algo:tau-TCSsimple} the pseudocode for the updates \eqref{eq:formulaS_d} and \eqref{eq:formulaS_n}.

\begin{algorithm}
\caption{$\tau$--TCS}\label{algo:tau-TCSsimple}
\begin{algorithmic}[1] \small
\State Choose $(\alpha^0; w^0) \in \R^{n+d}$, $\gamma >0$, $\tau_d,\tau_n\in \N$ and $b \in (0,1).$
\State Let $v \sim B(b)$ be a Bernoulli random variable (the coin toss)
\For{$k = 0, 1, \cdots $}
	\State Sample $v \in \{0,1\}$
	\If {$v=0$}
	\State Sample  $B_d \subset \{1,\ldots, d\}$ with $|B_d | = \tau_d$ uniformly.
	\State Compute $y_d \in \R^{\tau_d}$ the least norm solution to 
		\State $\left(\frac{\mA_{B_d, :}\mA_{B_d, :}^\top}{\lambda^2n^2} + \mI_{\tau_d}\right)y_d  =\frac{\mA_{B_d, :}\alpha^k}{\lambda n} - w_{B_d}^k$
		\State Compute the updates
		\State $\colvec{\Delta \alpha^k \\ \Delta w^k}   = -\colvec{\frac{1}{\lambda n}\mA_{B_d, :}^\top  \\ -\mI_{B_d}} y_d $
	\Else
	\State Sample  $B_n \subset \{1,\ldots, n\}$ with $|B_n | = \tau_n$ uniformly.	
		\State Compute $y_n \in \R^{\tau_n}$ the least norm solution to 
		\State $\left([\nabla \Phi^k_{B_n}]^{\top} \nabla \Phi^k_{B_n} + \mI_{\tau_n}\right)y_n  =\alpha_{B_n}^k + \Phi^k_{B_n}$
		\State Compute the updates
		\State $\colvec{\Delta \alpha^k \\ \Delta w^k}   = -\colvec{\mI_{B_n} \\ \nabla \Phi^k_{B_n}} y_n $
	\EndIf
	\State $w^{k+1}  = w^k + \gamma \Delta w^k $
	\State $\alpha^{k+1}  = \alpha^k +\gamma  \Delta \alpha^k$
\EndFor
\State \textbf{return:} last iterate $\alpha^k$, $w^k$
\end{algorithmic}
\end{algorithm}

\section{Pseudo code and implementation details for GLMs}
\label{sec:pseudocodeextras}

We also provide a more efficient and detailed implementation of Algorithm~\ref{algo:tau-TCSsimple} in Algorithm~\ref{algo:tau-TCS} in this section.

It is beneficial to first understand Algorithm~\ref{algo:tau-TCSsimple} in the simple setting where $\tau_d= \tau_n =1.$  We refer to this setting as the Kaczmarz--TCS method.

\subsection{Kaczmarz--TCS}

Let $f_j \in \R^d$ ($e_i \in \R^n$) be the $j$th (the $i$th) unit coordinate vector in $\R^d$ (in $\R^n$, respectively). For $\mS_k =
\begin{bmatrix}
\mS_d & 0 \\
0 & 0
\end{bmatrix} $ with $\mS_d = f_j$, from \eqref{eq:formulaS_d} we get
\begin{eqnarray}\label{eq:K_S_d}
\colvec{\Delta \alpha^k \\ \Delta w^k} &=& - \frac{\frac{1}{\lambda n}\sum_{l=1}^na_{lj}\alpha_l^k - w_j^k}{\frac{1}{\lambda^2n^2}\sum_{l=1}^na_{lj}^2+1}
\begin{bmatrix}
\frac{1}{\lambda n} \left( \begin{array}{c}
a_{1j} \\
\vdots \\
a_{nj}
\end{array} \right) \\
- f_j
\end{bmatrix}.
\end{eqnarray}
For $\mS_k =
\begin{bmatrix}
0 & 0 \\
0 & \mS_n
\end{bmatrix} $ with $\mS_n = e_i$, from \eqref{eq:formulaS_n} we get
\begin{eqnarray}\label{eq:K_S_n}
\colvec{\Delta \alpha^k \\ \Delta w^k} &=& - \frac{\alpha_i^k + \phi_i'(a_i^\top w^k)}{\|a_i\|_2^2\phi_i''(a_i^\top w^k)^2 + 1}
\begin{bmatrix}
e_i \\
\phi_i''(a_i^\top w^k)a_i
\end{bmatrix}.
\end{eqnarray}

See Algorithm~\ref{algo:Kaczmarz-TCS} an efficient implementation of Algorithm~\ref{algo:tau-TCSsimple} in a single row sampling case. Notice that we introduce an auxiliary variable $\overline{\alpha}^k$ to update the term $\frac{1}{\lambda n}\sum_{l=1}^na_{lj}\alpha_l^k$ for $j = 1, \cdots, d$ in \eqref{eq:K_S_d} and we store a $d \times d$ matrix \texttt{cov} which can be seen as the covariance matrix of the dataset $\mA$ to update the term $\frac{1}{\lambda^2n^2}\sum_{l=1}^na_{lj}^2$ for $j = 1, \cdots, d$ in \eqref{eq:K_S_d} (see Algorithm~\ref{algo:Kaczmarz-TCS} Line 12). We also store a vector \texttt{sample} $\in \R^n$ to update the term $\norm{a_i}_2^2$ for $i = 1, \cdots, n$ in \eqref{eq:K_S_n} (see Algorithm~\ref{algo:Kaczmarz-TCS} Line 21).

\paragraph*{Cost per iteration analysis of Algorithm~\ref{algo:Kaczmarz-TCS}}
From Algorithm~\ref{algo:Kaczmarz-TCS}, the cost of computing \eqref{eq:K_S_d} is $\cO(n)$ with $n$ coordinates' updates of the auxiliary variable $\alpha$ (see Algorithm~\ref{algo:Kaczmarz-TCS} Line 15). This is affordable as the cost of each coordinate's update is $1$. Besides $\mA$ is often sparse. The update in this case can be much cheaper than $n$. Besides, the cost of computing \eqref{eq:K_S_n} is $\cO(d)$. If we choose the Bernoulli parameter $b = n/(n+d)$ which selects one row of $F$ uniformly, the total cost of the updates TCS in expectation with respect to the Bernoulli distribution will be
\[\mbox{Cost(update TCS)} = \cO(n) * (1-b) + \cO(d) * b = \cO(nd/(n+d)) = \cO(\min(n,d)).\]


So the TCS method can have the same cost per iteration as the stochastic first-order methods in the case $d < n$, such as SVRG~\cite{Johnson2013}, SAG~\cite{SAG}, dfSDCA~\cite{Shalev-Shwartz2015} and Quartz~\cite{Qu2015b}.

\begin{algorithm}
\caption{Kaczmarz-TCS}\label{algo:Kaczmarz-TCS}
\begin{algorithmic}[1]
\State \textbf{parameters:} $\cD =$ distribution over random matrices
\State \textbf{store in memory:}
\State \quad \quad \texttt{sample}: $(\|a_i\|_2^2)_{1 \leq i \leq n} \in \R^n$
\State \quad \quad \texttt{cov}: $\frac{1}{\lambda^2 n^2}\mA\mA^\top \in \R^{d \times d}$
\State \textbf{initialization:}
\State \quad \quad Choose $(\alpha^0, w^0) \in \R^n\times\R^d$ and a step size $\gamma \in \R^{++}$
\State \quad \quad Set $\overline{\alpha}^0 = \frac{1}{\lambda n}\mA\alpha^0$
\For{$k = 0, 1, \cdots $}
	\State sample a fresh tossing-coin sketching matrix: $\mS_k \sim \cD$
	\If{$\mS_k =
			\begin{bmatrix}
			\mS_d & 0 \\
			0 & 0
			\end{bmatrix} $ with $\mS_d = f_j$}
		\State \textbf{update} \eqref{eq:K_S_d}: \Comment{Sketch a linear system based on the first $d$ rows of the Jacobian}
		\State \quad \quad $\Delta w_j^k = \frac{\overline{\alpha}_j^k - w_j^k}{\mbox{\texttt{cov}}[j, \ j] + 1}$
		\State \quad \quad $\Delta \alpha^k = - \Delta w_j^k \cdot \frac{1}{\lambda n} \colvec{a_{1j} \\ \vdots \\ a_{nj}}$
		\State \quad \quad $w_j^{k+1} = w_j^k + \gamma \cdot \Delta w_j^k$ \Comment{$j$th coordinate's update of the variable $w^k$}
		\State \quad \quad $\alpha^{k+1} = \alpha^k + \gamma \cdot \Delta \alpha^k$ \Comment{full vector's update of the auxiliary variable $\alpha^k$}
		\State \quad \quad $\overline{\alpha}^{k+1} = \overline{\alpha}^k - \gamma \cdot \Delta w_j^k \cdot \mbox{\texttt{cov}}[:\ , \ j]$
	\Else
		\State $\mS_k =
			\begin{bmatrix}
			0 & 0 \\
			0 & \mS_n
			\end{bmatrix} $ with $\mS_n = e_i$
		\State \textbf{update} \eqref{eq:K_S_n}: \Comment{Sketch a system based on the last $n$ rows of the Jacobian}
		\State \quad \quad $\mbox{\texttt{temp}} = a_i^\top w^k$\Comment{temporal scalar}
		\State \quad \quad $\Delta \alpha_i^k = - \frac{\alpha_i^k + \phi_i'(\mbox{\texttt{temp}})}{\mbox{\texttt{sample}}[i] \cdot \phi_i''(\mbox{\texttt{temp}})^2+1}$
		\State \quad \quad $\Delta w^k = \Delta \alpha_i^k \cdot \phi_i''(\mbox{\texttt{temp}}) \cdot a_i$
		\State \quad \quad $w^{k+1} = w^k + \gamma \cdot \Delta w^k$ \Comment{full vector's update of the variable $w^k$}
		\State \quad \quad $\alpha_i^{k+1} = \alpha_i^k + \gamma \cdot \Delta \alpha_i^k$ \Comment{$i$th coordinate's update of the auxiliary variable $\alpha^k$}
		\State \quad \quad $\overline{\alpha}^{k+1} = \overline{\alpha}^k + \gamma \cdot \Delta \alpha_i^k \cdot \frac{1}{\lambda n} a_i$
	\EndIf
\EndFor
\State \textbf{return:} last iterate $\alpha^k$, $w^k$
\end{algorithmic}
\end{algorithm}

\subsection{\texorpdfstring{$\tau$}{t}--Block TCS}

Here we provide Algorithm~\ref{algo:tau-TCS} which is a detailed implementation of Algorithm~\ref{algo:tau-TCSsimple} in a more efficient way. Similar to Algorithm~\ref{algo:Kaczmarz-TCS} but with sketch sizes $\tau_d$ and $\tau_n$, we also store a $d \times d$ matrix \texttt{cov}, but not a vector \texttt{sample}. We refer to Algorithm~\ref{algo:tau-TCS} as the $\tau$-block TCS method.

\begin{algorithm}
\caption{$\tau$--Block TCS}\label{algo:tau-TCS}
\begin{algorithmic}[1]
\State \textbf{parameters:} $\cD =$ distribution over random matrices
\State \textbf{store in memory:}
\State \quad \quad \texttt{cov}: $\frac{1}{\lambda^2 n^2}\mA\mA^\top \in \R^{d \times d}$
\State \textbf{initialization:}
\State \quad \quad Choose $(\alpha^0, w^0) \in \R^n\times\R^d$ and a step size $\gamma \in \R^{++}$
\State \quad \quad Set $\overline{\alpha}^0 = \frac{1}{\lambda n}\mA\alpha^0$
\For{$k = 0, 1, \cdots $}
	\State sample a fresh tossing-coin sketching matrix: $\mS_k \sim \cD$
	\If {$\mS_k =
			\begin{bmatrix}
			\mS_d & 0 \\
			0 & 0
			\end{bmatrix} $ with $\mS_d = \mI_{B_d}$ and $|B_d| = \tau_d$} 
		\State \textbf{Update} \eqref{eq:formulaS_d}: \Comment{Sketch a linear system based on the first $d$ rows of the Jacobian}
		\State \quad \quad Compute $y_d \in \R^{\tau_d}$ the least norm solution to the $\tau_d \times \tau_d$ linear system
		\State \quad \quad $\left(\mbox{\texttt{cov}}[B_d, \ B_d] + \mI_{\tau_d}\right) y_d = - \left(\overline{\alpha}^k_{B_d} - w_{B_d}^k\right)$
		\State \quad \quad Compute the updates
		\State \quad \quad $w^{k+1}_{B_d} = w^k_{B_d} - \gamma \cdot y_d$ \Comment{$\tau_d$ coordinates' update of the variable $w^k$}
		\State \quad \quad $\alpha^{k+1} = \alpha^k + \gamma \cdot \frac{1}{\lambda n} \mA_{B_d, :}^\top y_d$ \Comment{full vector's update of the auxiliary variable $\alpha^k$}
		\State \quad \quad $\overline{\alpha}^{k+1} = \overline{\alpha}^k + \gamma \cdot \mbox{\texttt{cov}}[:\ , \ B_d] y_d$
	\Else
		\State $\mS_k =
			\begin{bmatrix}
			0 & 0 \\
			0 & \mS_n
			\end{bmatrix} $ with $\mS_n = \mI_{B_n}$ and $|B_n| = \tau_n$
		\State \textbf{Update} \eqref{eq:formulaS_n}: \Comment{Sketch a system based on the last $n$ rows of the Jacobian}
		\State \quad \quad $\mbox{\texttt{temp}} = \mA_{:, B_n}^\top w^k \in \R^{\tau_n}$ \Comment{Temporal vector}
		\State \quad \quad $\mD_{B_n}^k = \Diag{\phi_{B_n}''(\mbox{\texttt{temp}})} \in \R^{\tau_n \times \tau_n}$ \Comment{Compute $\phi_i''(a_i^\top w^k)$ element-wise $\forall i \in B_n$}
		\State \quad \quad $\nabla \Phi^k_{B_n} = \mA_{:, B_n}\mD_{B_n}^k \in \R^{d \times \tau_n}$
		\State \quad \quad $\Phi^k_{B_n} = \phi_{B_n}'(\mbox{\texttt{temp}})$ \Comment{Compute $\phi_i'(a_i^\top w^k)$ element-wise $\forall i \in B_n$}
		\State \quad \quad Compute $y_n \in \R^{\tau_n}$ the least norm solution to the $\tau_n \times \tau_n$ linear system
		\State \quad \quad $\left([\nabla \Phi^k_{B_n}]^\top \nabla \Phi^k_{B_n} + \mI_{\tau_n}\right) y_n = - \left(\alpha_{B_n}^k + \Phi^k_{B_n}\right) $
		\State \quad \quad Compute the updates
		\State \quad \quad $w^{k+1} = w^k + \gamma \cdot \nabla \Phi^k_{B_n} y_n$ \Comment{full vector's update of the variable $w^k$}
		\State \quad \quad $\alpha^{k+1}_{B_n} = \alpha^k_{B_n} + \gamma \cdot y_n$ \Comment{$\tau_n$ coordinates' update of the auxiliary variable $\alpha^k$}
		\State \quad \quad $\overline{\alpha}^{k+1} = \overline{\alpha}^k + \gamma \cdot \frac{1}{\lambda n} \mA_{:, B_n} y_n$
	\EndIf
\EndFor
\State \textbf{return:} last iterate $\alpha^k$, $w^k$
\end{algorithmic}
\end{algorithm}

From Algorithm~\ref{algo:tau-TCS}, the cost of solving the $\tau_d \times \tau_d$ system (see Algorithm~\ref{algo:tau-TCS} Line 12) is $\cO(\tau_d^3)$ for a direct solver and the cost of updating $\alpha$ and $\overline{\alpha}$ (see Algorithm~\ref{algo:tau-TCS} Line 15 and Line 16) are $\cO(\tau_d n)$ and $\cO(\tau_d d)$ respectively. Overall, this implies that the cost of executing the sketching of the first $d$ rows is
\begin{eqnarray} \label{eq:c_d}
c_d &\eqdef& \cO(\max(\tau_d^3, \tau_d n, \tau_d d)).
\end{eqnarray}
Similarly, the dominant cost of executing the last $n$ rows sketch comes from forming the $\tau_n \times \tau_n$ linear system or solving such system (see Line 25), which gives
\begin{eqnarray} \label{eq:c_n}
c_n \eqdef \cO(\max(\tau_n^3, \tau_n^2d)).
\end{eqnarray}
In average, which means taking the Bernoulli parameter $b$ into account, the total cost per iteration of the TCS updates in expectation is
\begin{eqnarray} \label{eq:c_avg}
c_{avg} &\eqdef& c_d \times (1-b) + c_n \times b \nonumber \\
&=& \cO(\max(\tau_d^3, \tau_d n, \tau_d d)) \times (1-b) + \cO(\max(\tau_n^3, \tau_n^2d)) \times b.
\end{eqnarray}
Depending on the sketch sizes $(\tau_d, \tau_n)$ and the Bernoulli parameter $b$, the nature of $c_{avg}$ can be different from $\cO(d)$ (see Kaczmarz-TCS in Algorithm~\ref{algo:Kaczmarz-TCS}) to $\cO(d^2)$ (see the cost per iteration analysis paragraph in the next section). We discuss the total cost per iteration of the TCS method in practice in different cases in the next section.

\section{Additional experimental details} \label{sec:detailexp}


All the sampling of the methods was pre-computed before starting counting the wall-clock time for each method and each dataset.
We also paused the timing when the algorithms were under process of the performance evaluation of the gradient norm or of the logistic regression loss that were necessary to generate the plots.

In the following, from the experimental results for GLM in Section~\ref{sec:experiments}, we discuss the parameters' choices for TCS in practice, including the sketching sizes $(\tau_d, \tau_n)$, the Bernoulli parameter $b$, the stepsize $\gamma$ and the analysis of total cost per iteration. See Table~\ref{tab:parameters} for the parameters we chose for TCS in the experiments in Figure~\ref{figure}. Such choices are due to TCS's cost per iteration.

\begin{table}
  \caption{Details of the parameters' choices ($\gamma$ and $b$) for $50$-TCS, $150$-TCS and $300$-TCS}
  \label{tab:parameters}
  \centering
  \begin{tabular}{ ccccc }
  	 \toprule
  	             &            & $50$-TCS                        & $150$-TCS                      & $300$-TCS             \\
     \cmidrule(r){3-5}
     dataset     & stepsize   & Bernoulli              & Bernoulli             & Bernoulli    \\
     \midrule
     covetype    & $1.0$      & $\frac{n}{n + \tau_n * 3}$      & $\frac{n}{n + \tau_n * 3}$     & $\frac{n}{n + \tau_n * 3}$       \\
     a9a         & $1.5$      & $\frac{n}{n + \tau_n} - 0.03$   & $\frac{n}{n + \tau_n} - 0.03$  & $\frac{n}{n + \tau_n} - 0.11$    \\
     fourclass   & $1.0$      & $\frac{n}{n + \tau_n} - 0.11$   & $\frac{n}{n + \tau_n} - 0.11$  & $\frac{n}{n + \tau_n} - 0.11$    \\
     artificial  & $1.0$      & $\frac{n}{n + \tau_n} - 0.03$   & $\frac{n}{n + \tau_n} - 0.11$  & $\frac{n}{n + \tau_n} - 0.11$    \\
     ijcnn1      & $1.8$      & $\frac{n}{n + \tau_n} - 0.03$   & $\frac{n}{n + \tau_n} - 0.11$  & $\frac{n}{n + \tau_n} - 0.11$    \\
     webspam     & $1.8$      & $\frac{n}{n + \tau_n * 3}$      & $\frac{n}{n + \tau_n * 3}$     & $\frac{n}{n + \tau_n * 3}$       \\
     epsilon     & $1.8$      & $\frac{n}{n + \tau_n * 3}$      & $\frac{n}{n + \tau_n * 3}$     & $\frac{n}{n + \tau_n * 3}$       \\
     phishing    & $1.8$      & $\frac{n}{n + \tau_n} - 0.03$   & $\frac{n}{n + \tau_n} - 0.11$  & $\frac{n}{n + \tau_n} - 0.11$    \\
     \bottomrule
  \end{tabular}
\end{table}

\paragraph*{Choice of the sketch size $\tau_d$} For all of our experiments, 
 $\tau_d = d$ performs always the best in time and in number of iterations. That means, when $\mS_k = \begin{bmatrix} \mS_d & 0 \\ 0 & 0\end{bmatrix}$ at $k$th iteration, we choose $\mS_d = \mI_d$. 
 Note also that the first $d$ rows $\eqref{eq:F}$ are linear, thus using $\mS_d = \mI_d$ gives an exact solution to these first $d$ equations. We found that such an exact solution from the linear part induces a fast convergence when $d < n$.  We did not test datasets for which $d > n$ with $d$ very large.

\begin{figure}
\centering
\begin{tabular}{cc}
\includegraphics[width=.46\linewidth]{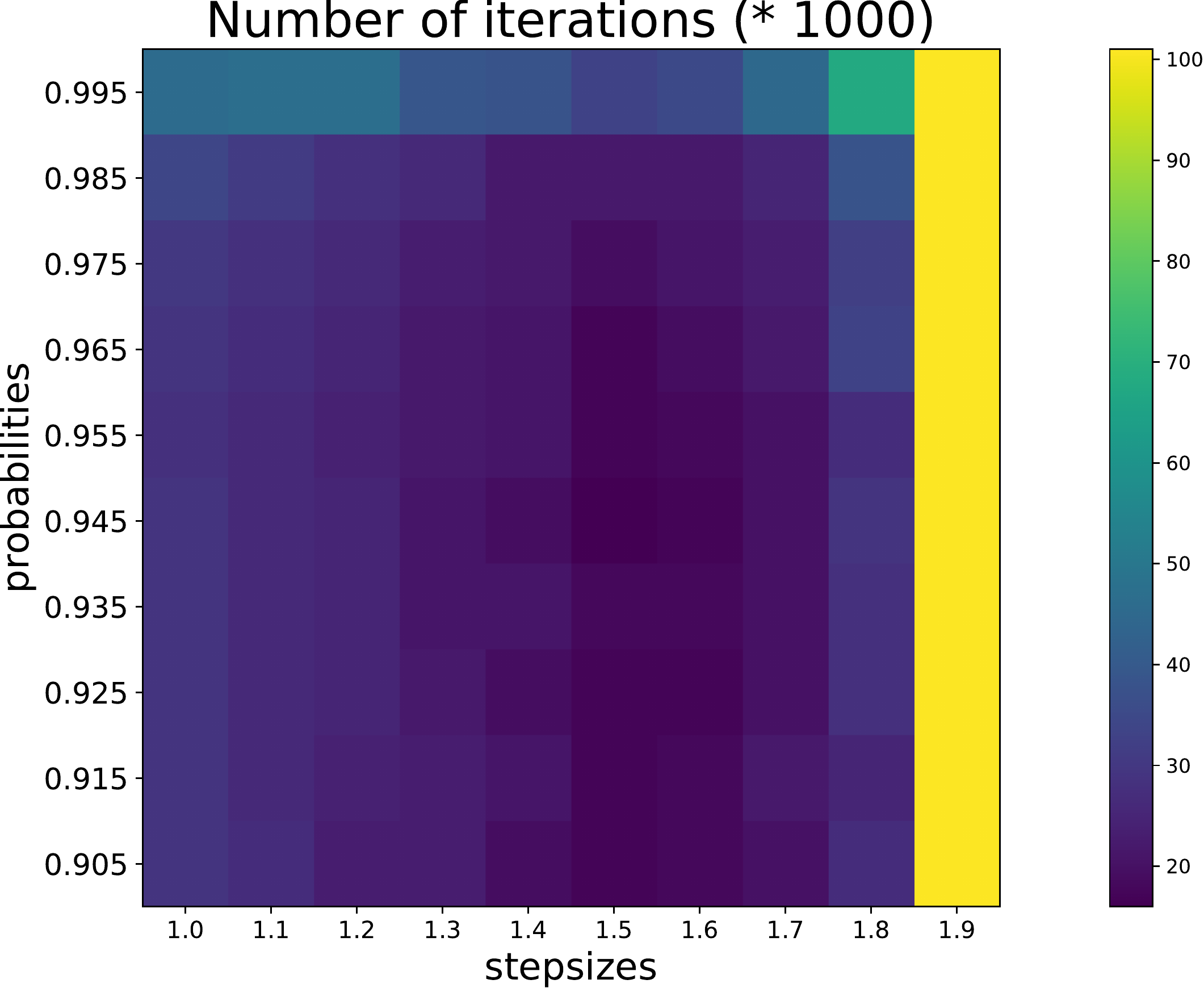}&
\includegraphics[width=.46\linewidth]{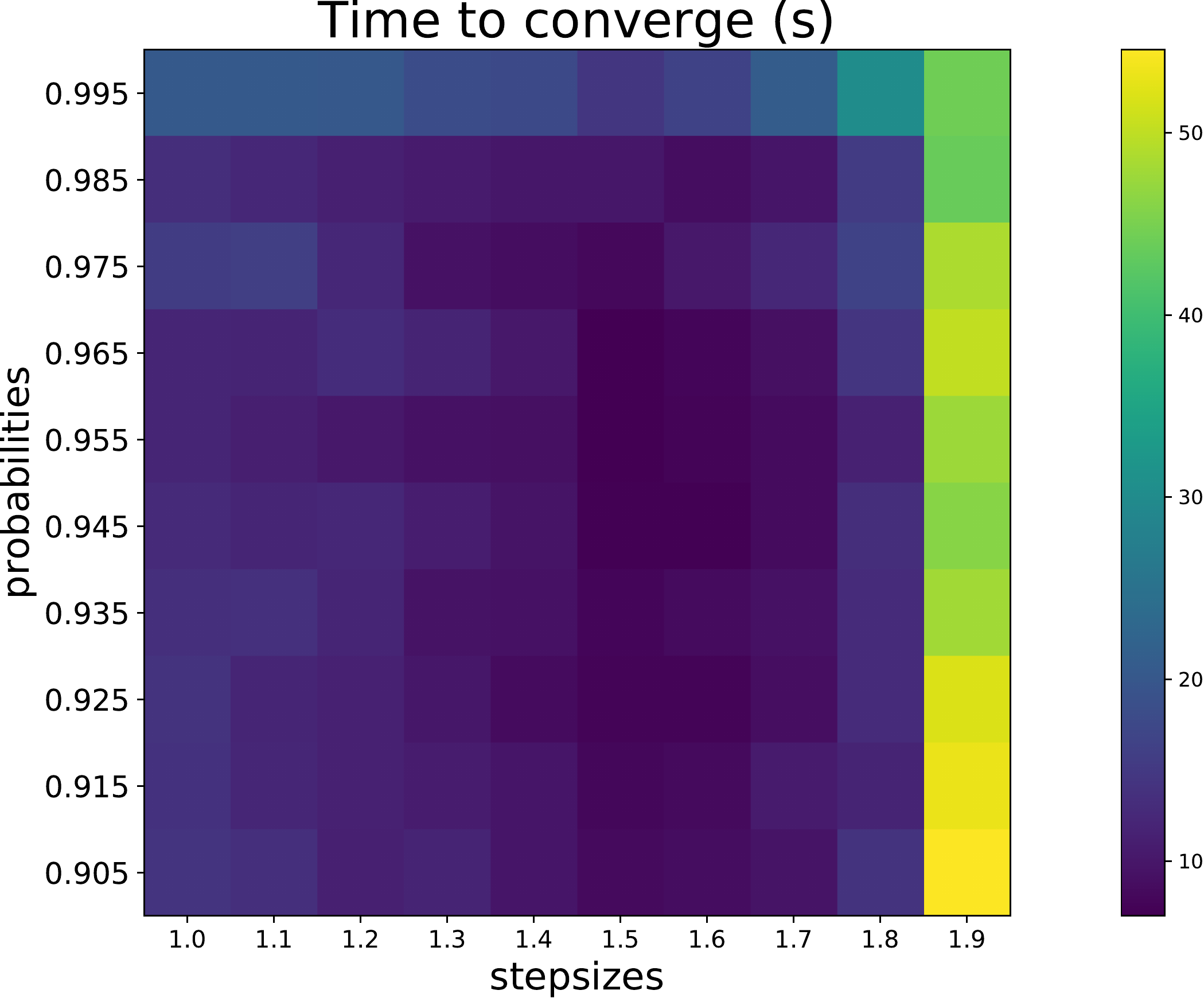}
\end{tabular}
\caption{a9a dataset: Grid search of the Bernoulli parameter $b$ and the stepsize $\gamma$ with $150$-TCS method.The darker colors correspond to a resulting small gradient norm and thus a better solution.}
\label{fig:a9a}
\end{figure}

\paragraph*{Choice of the Bernoulli parameter $b$ for uniform sampling} First, we calculate the probability of sampling one row of the function $F$ \eqref{eq:F}. Since there exists two types of sketching for TCS method depending on the \emph{coin toss}, we address both of them. The probability of sampling one specific row of the first block (fist $d$ rows of \eqref{eq:F}) is
\[p_d = \frac{\tau_d}{d} \times (1-b)\]
and the one of the second block is
\[p_n = \frac{\tau_n}{n} \times b.\]
It is natural to choose $b$ such that the uniform sampling of the whole system, i.e. $p_d = p_n$, is guaranteed. This implies to set
\[p_{uniform} \eqdef \frac{\tau_dn}{\tau_dn+\tau_nd}.\]
As we choose $\tau_d = d$, this implies
\[p_{uniform} = \frac{n}{n + \tau_n}.\]
However, we found through multiple experiments that when setting $b$ slightly smaller than $p_{uniform}$ (e.g.\ $- 1\%$), this reduces significantly the number of iterations to get convergence. See in Figure~\ref{fig:a9a} for a grid search of the Bernoulli parameter $b$ and the stepsize $\gamma$ for a9a dataset with $b = p_{uniform} = 0.995$ in the first line of the figure. Before giving details about how to choose $b$ in practice, we first provide the cost per iteration analysis of the TCS method in detail.

\paragraph*{Total cost per iteration in expectation analysis for TCS in different cases}
Recall two types of costs per iteration $c_d$ \eqref{eq:c_d} and $c_n$ \eqref{eq:c_n}. In our cases, consider $b = \frac{n}{n + \tau_n}$ and $d = \tau_d < n$. To summarize, the cost per iteration in expectation can be one of the three following cases followed with their bounds:
\begin{itemize}
	\item[1.] If $\tau_n < \sqrt{n} < d < n$ such as \emph{epsilon} dataset, then $c_d = \cO(d^3) > c_n = \cO(\tau_n^2d)$, and
	\begin{align} \label{eq:c_avg1}
	c_{avg1} &= \cO(d^3) \times \left(1 - \frac{n}{n + \tau_n}\right) + \cO(\tau_n^2d) \times \frac{n}{n + \tau_n} \nonumber \\
	&= \cO(\frac{\tau_nd}{n + \tau_n}(d^2 + \tau_nn)) \nonumber \\
	\Longrightarrow \cO(\tau_n^2d) \; \leq \; &c_{avg1} \; \leq \; \cO(\tau_nd^2);
	\end{align}
	\item[2.] if $\tau_n < d < \sqrt{n}$ such as \emph{webspam} dataset with $50$-TCS and $150$-TCS, \emph{a9a}, \emph{phishing} and \emph{covtype} datasets with $50$-TCS method, then $c_d = \cO(dn) > c_n = \cO(\tau_n^2d)$, and
	\begin{align} \label{eq:c_avg2}
	c_{avg2} &= \cO(dn) \times \left(1 - \frac{n}{n + \tau_n}\right) + \cO(\tau_n^2d) \times \frac{n}{n + \tau_n} \nonumber \\
	&= \cO(\frac{\tau_nn}{n + \tau_n}(d + \tau_nd)) = \cO(\tau_n^2d) \quad \mbox{as } \frac{1}{2} \leq \frac{n}{n + \tau_n} < 1;
	\end{align}
	\item[3.] if $d < \sqrt{n}$ and $d < \tau_n$ such as all the other experiments for TCS methods in Figure~\ref{figure}, then $c_d = \cO(dn)$, $c_n = \cO(\tau_n^3)$, and
	\begin{align} \label{eq:c_avg3}
	c_{avg3} &= \cO(dn) \times \left(1 - \frac{n}{n + \tau_n}\right) + \cO(\tau_n^3) \times \frac{n}{n + \tau_n} \nonumber \\
	&= \cO(\frac{\tau_nn}{n + \tau_n}(d + \tau_n^2)) \nonumber \\
	&= \cO(\tau_n^3) > \cO(d^3) \quad \mbox{as } \frac{1}{2} \leq \frac{n}{n + \tau_n} < 1.
	\end{align}
\end{itemize}
Notice that $c_{avg1}, c_{avg2}, c_{avg3} \ll \cO(dn)$ in general for large scale datasets with large $n$. For example, $c_{avg1} < \cO(dn)$ when $\tau_n d < n$. This justifies that TCS method is cheaper than the first-order method which requires evaluating the full gradient and thus has a cost per iteration of at least $\cO(dn)$. From $c_{avg1}$ \eqref{eq:c_avg1}, we know that TCS method can have the same cost per iteration as the stochastic first-order methods which is $\cO(d)$ in practice, such as SVRG~\cite{Johnson2013}, SAG~\cite{SAG}, dfSDCA~\cite{Shalev-Shwartz2015} and Quartz~\cite{Qu2015b}.

Furthermore, from the above analysis of computational cost, we can easily obtain the comparisons between $c_d$ and $c_n$ for different datasets and different sketch sizes in Table~\ref{tab:complex}. These comparisons  helped us to choose $b$ as we detail in the following.

\begin{table}
  \caption{Cost per iteration for different datasets and different sketch sizes.}
  \label{tab:complex}
  \centering
  \begin{tabular}{ cccc }
     \toprule
     dataset     & $50$-TCS    & $150$-TCS   & $300$-TCS   \\
     \midrule
     covetype    & $c_d > c_n$ & $c_d > c_n$ & $c_d > c_n$ \\
	 a9a         & $c_d > c_n$ & $c_d > c_n$ & $c_d < c_n$ \\
	 fourclass   & $c_d < c_n$ & $c_d < c_n$ & $c_d < c_n$ \\
	 artificial  & $c_d > c_n$ & $c_d < c_n$ & $c_d < c_n$ \\
	 ijcnn1      & $c_d > c_n$ & $c_d < c_n$ & $c_d < c_n$ \\
	 webspam     & $c_d > c_n$ & $c_d > c_n$ & $c_d > c_n$ \\
	 epsilon     & $c_d > c_n$ & $c_d > c_n$ & $c_d > c_n$ \\
	 phishing    & $c_d > c_n$ & $c_d < c_n$ & $c_d < c_n$ \\
     \bottomrule
  \end{tabular}
\end{table}

\paragraph*{Choice of the Bernoulli parameter $b$ in practice}
From the above discussion about $p_{uniform}$, heuristically, we decrease $b$ from $p_{uniform}$ to achieve faster convergence. For a large range of choices $b$, TCS converges. However, $b$ affects directly the computational cost per iteration. From \eqref{eq:c_avg}, we know that if $c_d > c_n$, decreasing $b$ will increase the average cost of the method. In this case, there is a trade-off between the number of iterations and the average cost to achieve the fastest convergence in time (see Figure~\ref{fig:a9a}). For a large dataset with $n$ large such as \emph{epsilon, webspam and covtype}, we decrease $b$ slightly, as for a small dataset, we make a relatively big decrease for $b$. If $c_d < c_n$, decreasing $b$ will also decrease the average cost. In this case, we tend to decrease $b$ even further. See Table~\ref{tab:parameters} the choices of $b$.

\begin{figure}
\centering
\begin{tabular}{cccc}
\includegraphics[width=.21\linewidth]{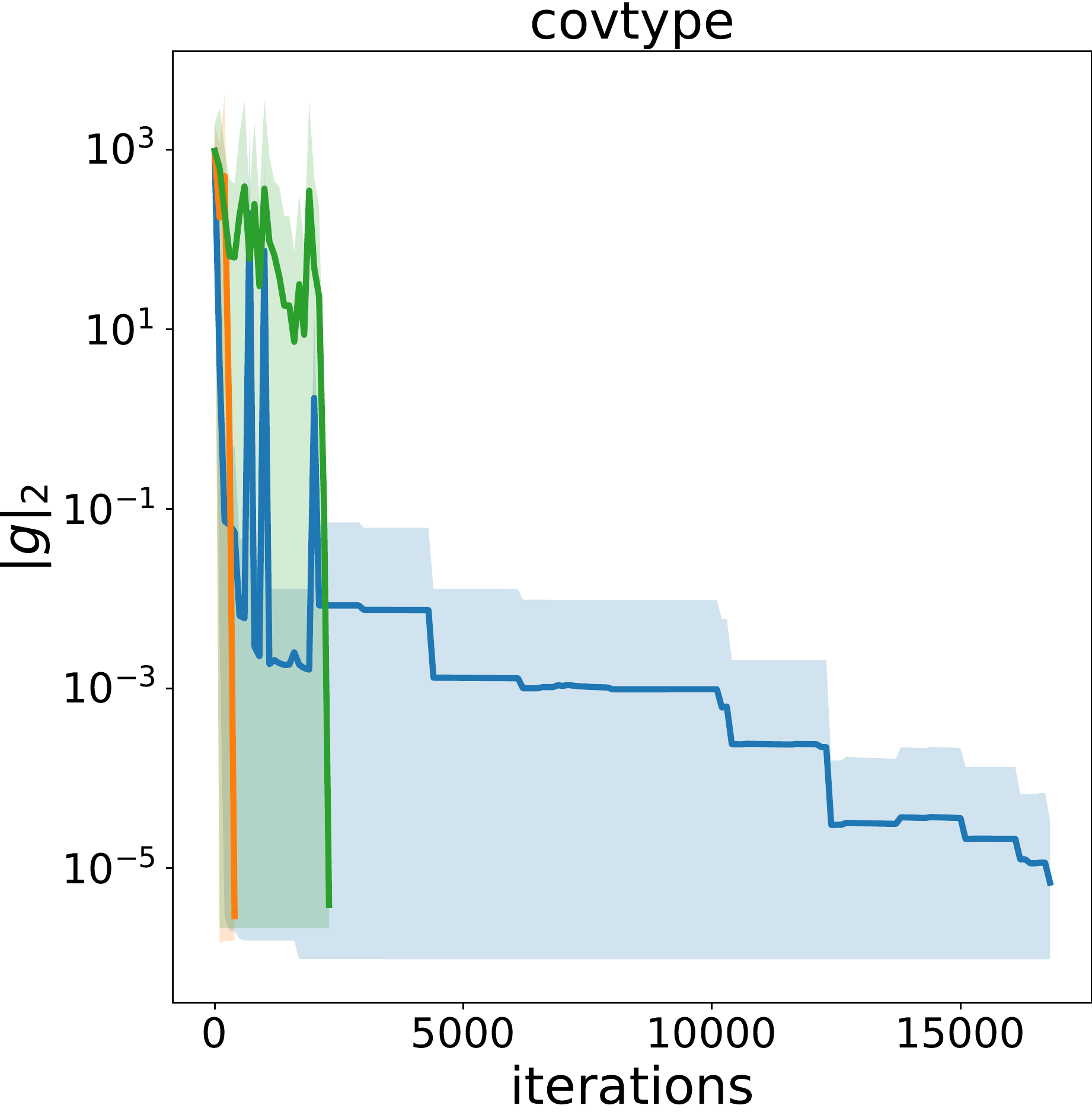}&
\includegraphics[width=.21\linewidth]{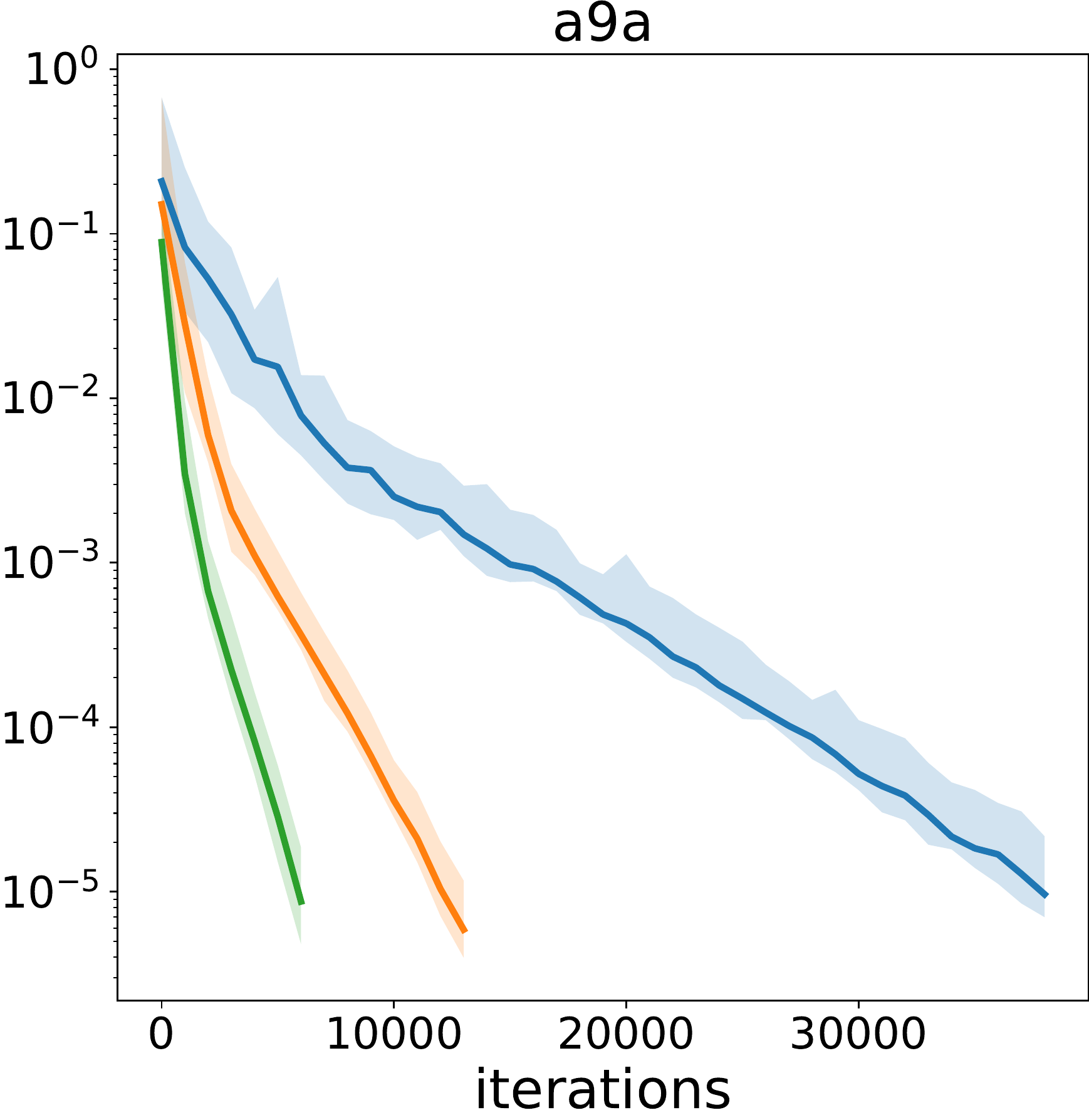}&
\includegraphics[width=.21\linewidth]{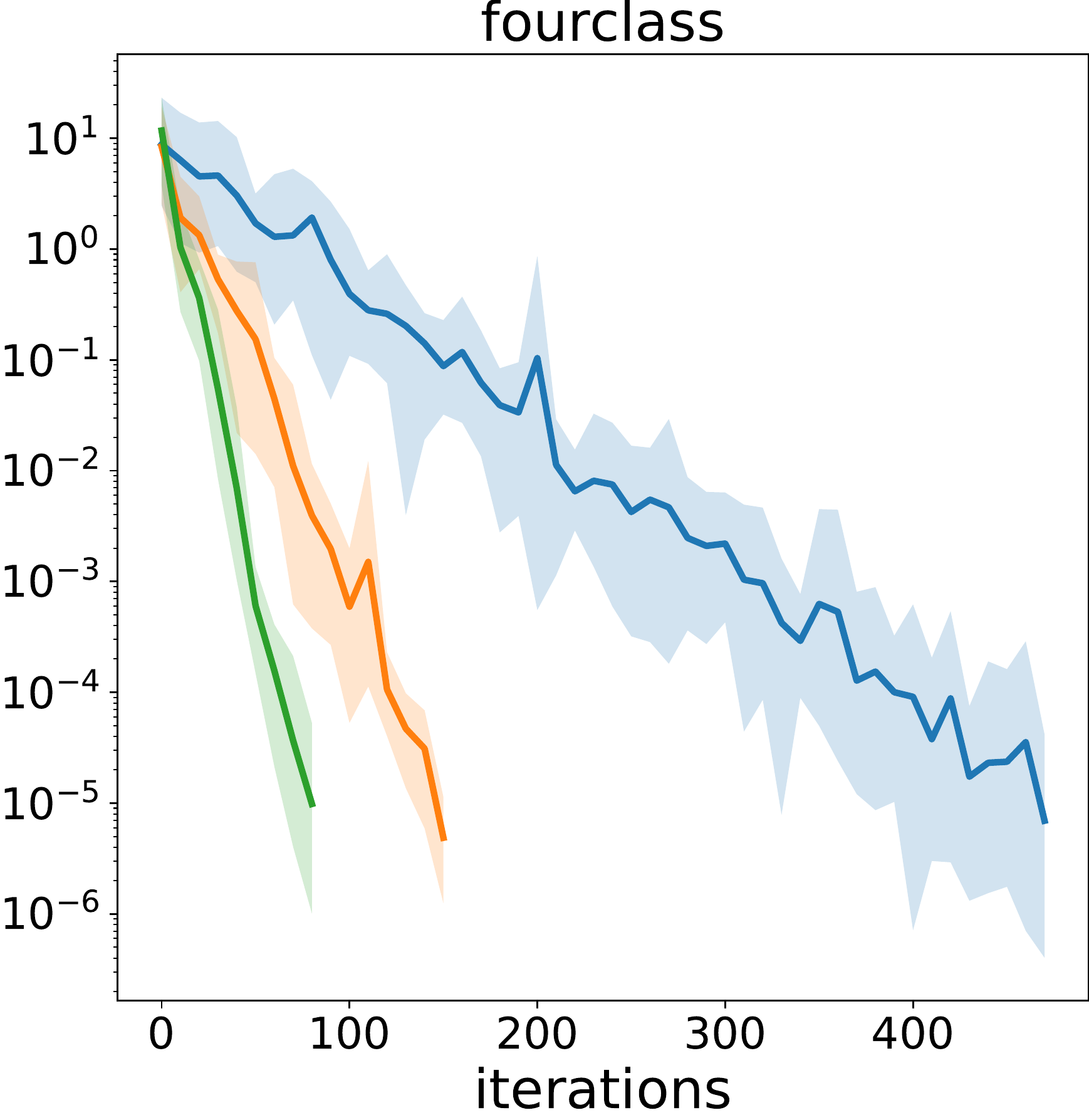}&
\includegraphics[width=.21\linewidth]{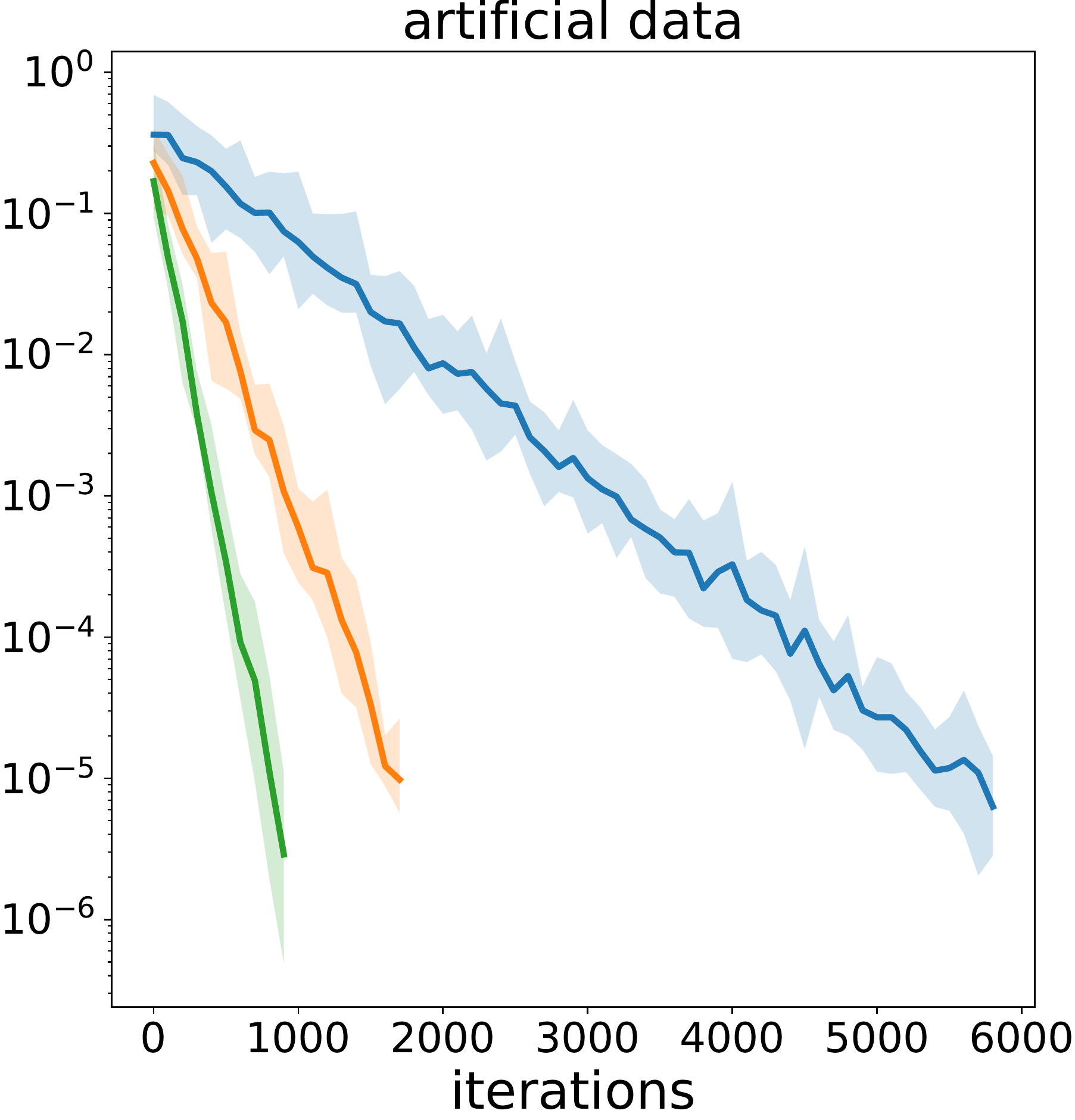} \\
\includegraphics[width=.21\linewidth]{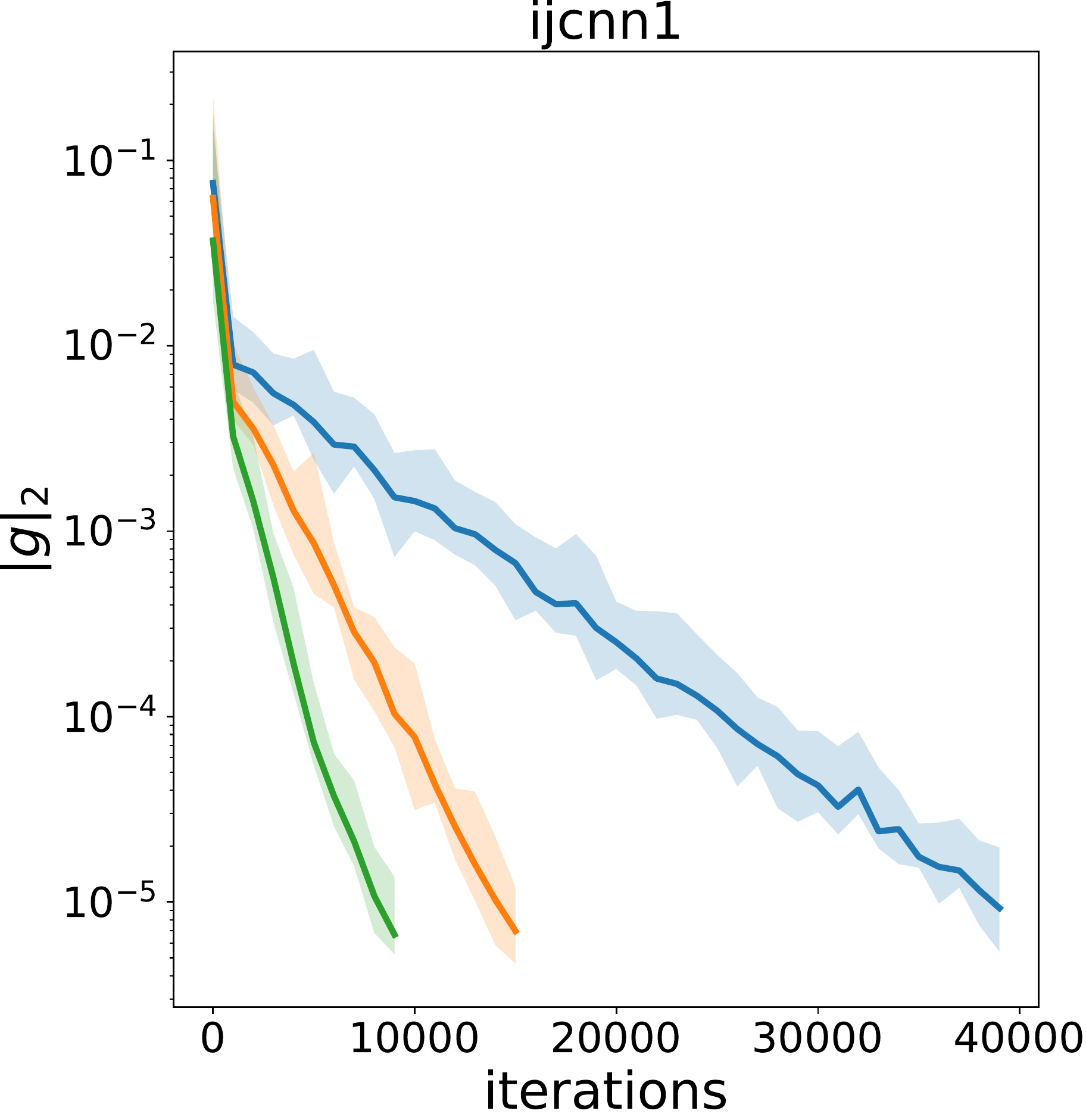}&
\includegraphics[width=.21\linewidth]{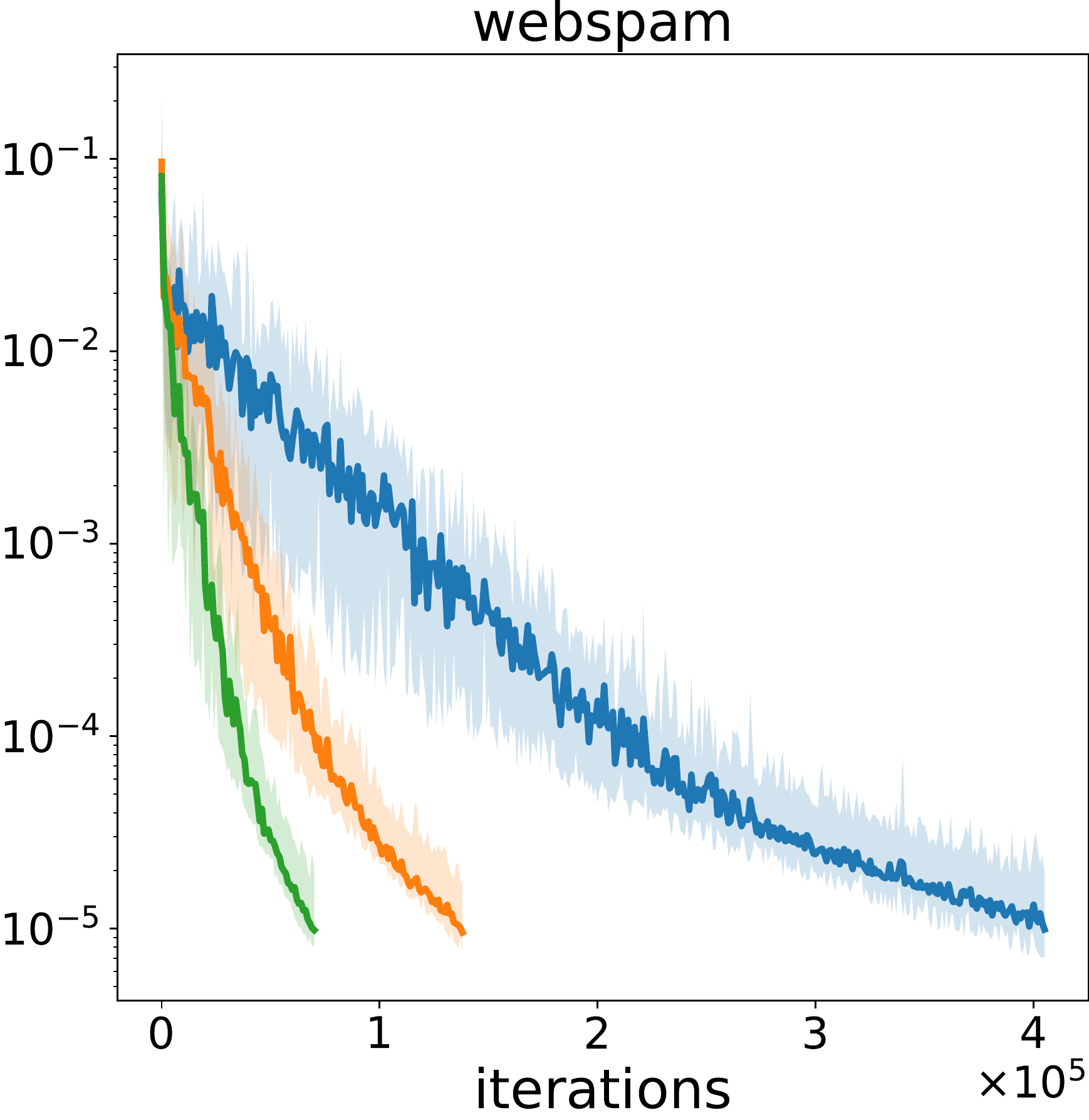}&
\includegraphics[width=.21\linewidth]{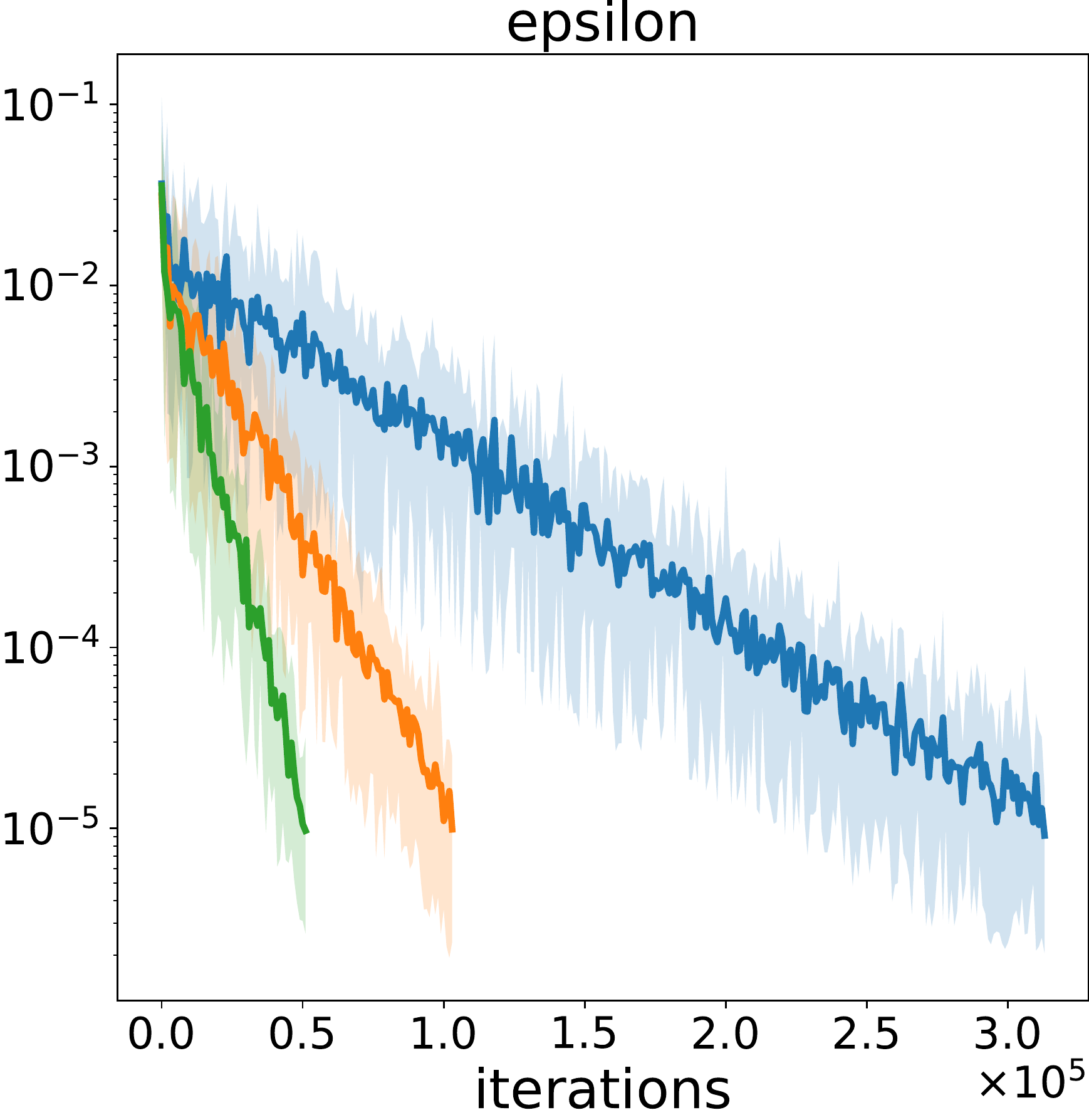}&
\includegraphics[width=.21\linewidth]{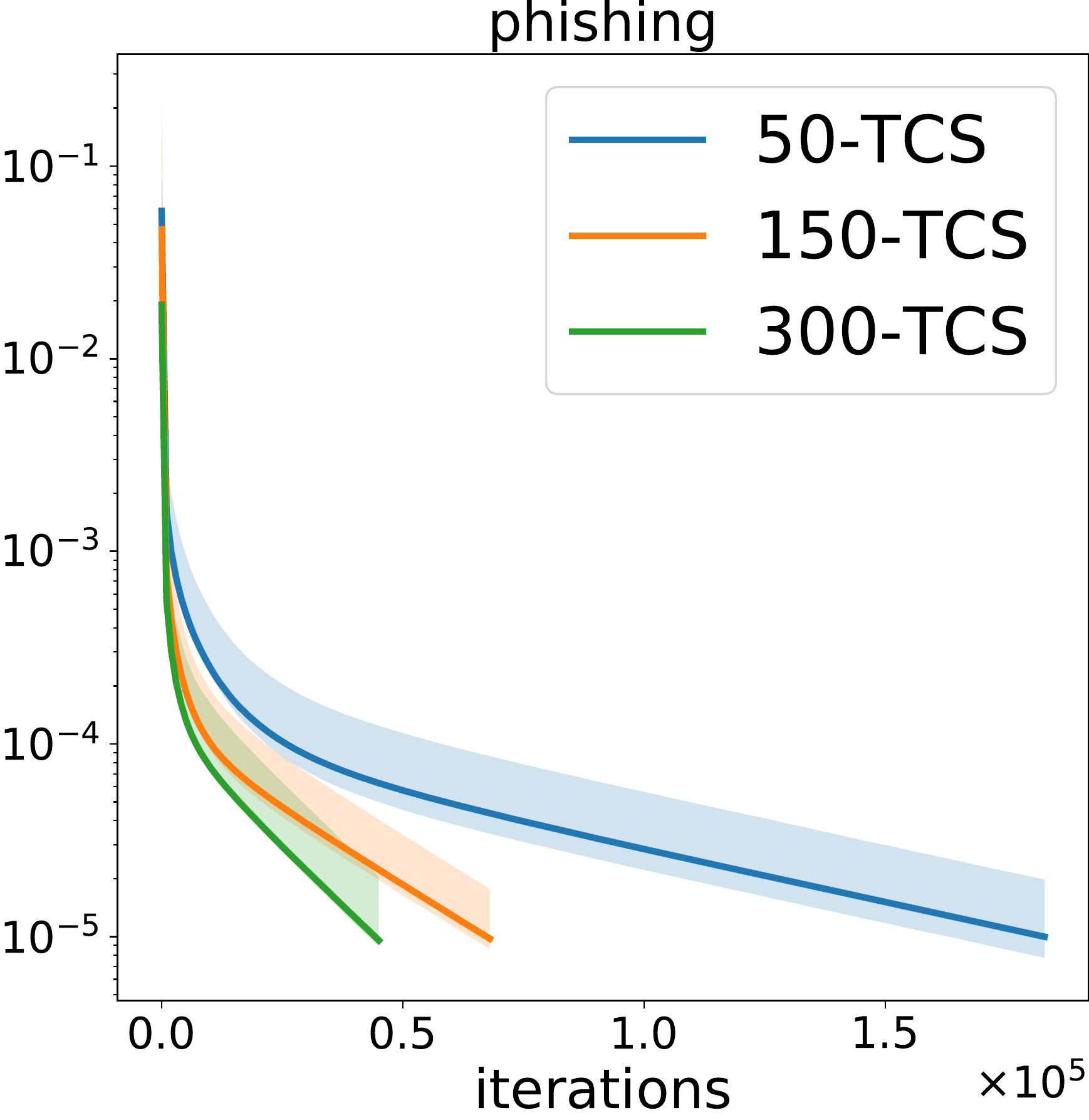}
\end{tabular}
\caption{Comparisons of different sketch sizes for TCS method in terms of the number of iterations.}
\label{fig:iter}
\end{figure}

\paragraph*{Choice of the sketch size $\tau_n$}
As for $\tau_n$, we observe that with bigger sketch size $\tau_n$, the method requires less number of iterations to get convergence. From Figure~\ref{fig:iter}, this is true for all the datasets except for \emph{covtype} dataset. However, choosing bigger sketch size $\tau_n$ will also increase the cost per iteration. Consequently, there exists an optimal sketch size such that the method converges the fastest in time taking account the balance between the number of iterations and the cost per iteration. From the experiments in Figure~\ref{figure}, we show that $\tau_n = 150$ is in general a very good choice for any scale of $n$.



\paragraph*{Choice of the stepsizes}
Different to our global convergence theories, in practice, choosing constant stepsize $\gamma > 1$ may converge faster (see Figure~\ref{fig:a9a}) for certain datasets. Here we need to be careful that the stepsize we mentioned is the stepsize used for the sketch $\mS_k = \begin{bmatrix} 0 & 0 \\ 0 & \mS_n \end{bmatrix}$. As for the sketch $\mS_k = \begin{bmatrix} \mS_d & 0 \\ 0 & 0 \end{bmatrix}$, we always choose stepsize $\gamma = 1$. Because stepsize $\gamma = 1$ solves exactly the linear system. Henceforth, we use $\gamma$ to designate  the stepsize used for the sketch of the last $n$ rows of $F$ \eqref{eq:F}. In our experiments, we find that the choice of the stepsize is related to the condition number (C.N.) of the model. If the dataset is ill-conditioned with a big C.N.\ of the model, $\gamma = 1$ is a good choice (see Figure~\ref{figure} top row and Table~\ref{tab:parameters} first four lines except for a9a); if the dataset is well-conditioned with a small C.N.\ of the model, all $\gamma \in (1, 1.8]$ gets convergence (see Figure~\ref{fig:a9a}). In practice, $\gamma = 1.8$ is a good choice for well-conditioned datasets (see Figure~\ref{figure}). However, from the grid search of stepsizes for a9a (see Figure~\ref{fig:a9a}), we know that the optimal stepsize for a9a is $\gamma = 1.5$.
To avoid tuning the stepsizes, i.e. a grid search procedure, we will apply a stochastic line search process~\cite{Stochastic_Line_Search_nips2019} in the next Section~\ref{sec:sls}.

Furthermore, we observe that the stepsize $\gamma$ is highly related to the smoothness constant $L$. If $L$ is big, then we choose $\gamma$ close to $1$, if $L$ is small, we increase $\gamma$ until $\gamma = 1.8$ (see Table~\ref{tab:parameters}). Such observation remains conjecture.

Finally, to summarize in practice for the TCS method with $d < n$, we choose $\tau_d = d$ and $\tau_n = 150$, we choose $b$ following the guideline introduced above; we always choose stepsize $\gamma = 1$ for the sketch of first $d$ rows \eqref{eq:F}; as for the sketch of last $n$ rows \eqref{eq:F}, we choose stepsize $\gamma = 1$ if the dataset is ill-conditioned and we can choose stepsize $\gamma = 1.8$ if the dataset is well-conditioned.




\section{Stochastic line-search for TCS methods applied in GLM} \label{sec:sls}

In order to avoid tuning the stepsizes, we can modify Algorithm \ref{algo:SkeNeR} by applying a stochastic line-search introduced by \cite{Stochastic_Line_Search_nips2019}. This is because again \SNR can be interpreted as a \SGD method.  It is a \emph{stochastic} line-search because on the $k$th iteration we sample a \emph{stochastic} sketching matrix $\mS_k$, and search for a stepsize $\gamma_k$ satisfying the following condition:
\begin{eqnarray}\label{eq:stoch_line_search}
f_{\mS_k, w^k}\left(w^k - \gamma_k\nabla f_{\mS_k, w^k}(w^k)\right) &\leq& f_{\mS_k, w^k}(w^k) - c \cdot \gamma_k \norm{\nabla f_{\mS_k, w^k}(w^k)}^2 \nonumber \\
&\overset{\eqref{eq:1smooth}}{=}& \left(1-2c\gamma_k\right)f_{\mS_k, w^k}(w^k).
\end{eqnarray}
Here, $c > 0$ is a hyper-parameter, usually a value close to $0$ is chosen in practice.

\subsection{Stochastic line-search for TCS method}

Now we focus on GLMs, which means we develop the stochastic line-search based on \eqref{eq:stoch_line_search} for TCS method. At $k$th iteration, if $\mS_k =
\begin{bmatrix}
\mS_d & 0 \\
0 & 0
\end{bmatrix}$ with $\mS_d = \mI_{B_d}$, we sketch a linear system based on the first $d$ rows of the Jacobian \eqref{eq:DFalphax}. Because of this linearity, the function $f_{\mS_k, (\alpha^k, w^k)}(\alpha^k; w^k)$ is quadratic. Thus \eqref{eq:stoch_line_search} can be re-written as
\begin{eqnarray} \label{eq:stoch_line_search_GLM}
f_{\mS_k, k}\left(\colvec{\alpha^k \\ w^k} - \gamma_k\nabla f_{\mS_k, k}(\alpha^k; w^k)\right) &=& (1-\gamma_k)^2f_{\mS_k, k}(\alpha^k; w^k) \nonumber \\ 
&\leq& (1-2c\gamma_k)f_{\mS_k, k}(\alpha^k; w^k),
\end{eqnarray}
where we use the shorthand $f_{\mS_k, k}(\alpha^k; w^k) \equiv f_{\mS_k, (\alpha^k; w^k)}(\alpha^k; w^k)$. To achieve the Armijo line-search condition \eqref{eq:stoch_line_search_GLM}, it suffices to take $\gamma = 1$ and $0 < c \leq \frac{1}{2}$ which is a common choice. Consequently, we do not need extra function evaluations. It is also well known that stepsize equal to $1$ is optimal as for Newton's method applied in quadratic problems.

If $\mS_k =
\begin{bmatrix}
0 & 0 \\
0 & \mS_n
\end{bmatrix}$ with $\mS_n = \mI_{B_n}$, we have
\begin{eqnarray}
f_{\mS_k, k}(\alpha^k; w^k) &=& \frac{1}{2}\left(\alpha^k_{B_n} + \Phi^k_{B_n}\right)^\top\left([\nabla\Phi^k_{B_n}]^\top\nabla\Phi^k_{B_n} + \mI_{\tau_n}\right)^\dagger\left(\alpha^k_{B_n} + \Phi^k_{B_n}\right). \label{eq:f_alpha_x_S_n}
\end{eqnarray}
and
\begin{align}
& f_{\mS_k, k}\left(\colvec{\alpha^k \\ w^k} - \gamma_k\nabla f_{\mS_k, k}(\alpha^k; w^k)\right) \nonumber \\
&= \frac{1}{2}F(\alpha^k + \gamma_k\Delta\alpha^k; w^k + \gamma_k\Delta w^k)^\top\mH_{\mS_k}(\alpha^k; w^k) F(\alpha^k + \gamma_k\Delta\alpha^k; w^k + \gamma_k\Delta w^k) \nonumber \\
&= \frac{1}{2}F(\alpha^k + \gamma_k\Delta\alpha^k; w^k + \gamma_k\Delta w^k)^\top \colvec{0 \\ \mI_{B_n}} \left([\nabla\Phi^k_{B_n}]^\top\nabla\Phi^k_{B_n} + \mI_{\tau_n}\right)^\dagger \nonumber \\
& \quad \colvec{0 \\ \mI_{B_n}}^\top F(\alpha^k + \gamma_k\Delta\alpha^k; w^k + \gamma_k\Delta w^k) \label{eq:f_alpha_x_S_n_delta}
\end{align}
with
\begin{align}
& \colvec{0 \\ \mI_{B_n}}^\top F(\alpha^k + \gamma_k\Delta\alpha^k; w^k + \gamma_k\Delta w^k) \nonumber \\
&= \alpha^k_{B_n} + \gamma_k\mI_{B_n}^\top\Delta\alpha^k + \phi'_{B_n}\left(\mA_{:, B_n}^\top w^k + \gamma_k\mA_{:, B_n}^\top \Delta w^k\right).
\end{align}
By \eqref{eq:formulaS_n}, we recall that
\begin{eqnarray}
\Delta \alpha^k &=& -\mI_{B_n}\left([\nabla\Phi^k_{B_n}]^\top\nabla\Phi^k_{B_n} + \mI_{\tau_n}\right)^\dagger\left(\alpha^k_{B_n} + \Phi^k_{B_n}\right), \label{eq:delta_alpha_S_n} \\
\Delta w^k &=& -\nabla\Phi^k_{B_n}\left([\nabla\Phi^k_{B_n}]^\top\nabla\Phi^k_{B_n} + \mI_{\tau_n}\right)^\dagger\left(\alpha^k_{B_n} + \Phi^k_{B_n}\right) \label{eq:delta_x_S_n}.
\end{eqnarray}
Note that the cost for evaluating \eqref{eq:f_alpha_x_S_n} and \eqref{eq:f_alpha_x_S_n_delta} are $\cO(\tau_n)$ and $\cO\left(\max\left(\tau_n^3, \tau_nd\right)\right)$ respectively, which are not expensive. Because one part of them are essentially a by-product from the computation of $y_n$, $\Delta \alpha^k$ and $\Delta w^k$ in Algorithm \ref{algo:tau-TCSsimple}. See Algorithm~\ref{algo:tau-TCS+sls} the implementation of TCS combined with the stochastic Armijo line-search. $\beta \in (0,1)$ is a discount factor.

\begin{algorithm}
\caption{$\tau$\texttt{--TCS+Armijo}}\label{algo:tau-TCS+sls}
\begin{algorithmic}[1] \small
\State Choose $(\alpha^0; w^0) \in \R^{n+d}$, $c, \beta, \gamma >0$, $\tau_d,\tau_n\in \N$ and $b \in (0,1).$
\State Let $v \sim B(b)$ be a Bernoulli random variable (the coin toss)
\For{$k = 0, 1, \cdots $}
	\State Sample $v \in \{0,1\}$
	\If {$v=0$}
		\State Sample  $B_d \subset \{1,\ldots, d\}$ with $|B_d | = \tau_d$ uniformly.
		\State Compute $y_d \in \R^{\tau_d}$ the least norm solution to 
		\State $\left(\frac{\mA_{B_d, :}\mA_{B_d, :}^\top}{\lambda^2n^2} + \mI_{\tau_d}\right)y_d  =\frac{\mA_{B_d, :}\alpha^k}{\lambda n} - w_{B_d}^k$
		\State Compute the updates
		\State $\colvec{\Delta \alpha^k \\ \Delta w^k}   = -\colvec{\frac{1}{\lambda n}\mA_{B_d, :}^\top  \\ -\mI_{B_d}} y_d $
		\State $w^{k+1}  = w^k + \Delta w^k $
		\State $\alpha^{k+1}  = \alpha^k + \Delta \alpha^k$
	\Else
		\State Reset $\gamma$ to the initial stepsize.
		\State Sample  $B_n \subset \{1,\ldots, n\}$ with $|B_n | = \tau_n$ uniformly.	
		\State Compute $y_n \in \R^{\tau_n}$ the least norm solution to 
		\State $\left([\nabla \Phi^k_{B_n}]^{\top} \nabla \Phi^k_{B_n} + \mI_{\tau_n}\right)y_n  =\alpha_{B_n}^k + \Phi^k_{B_n}$
		\State Compute the updates
		\State $\colvec{\Delta \alpha^k \\ \Delta w^k}   = -\colvec{\mI_{B_n} \\ \nabla \Phi^k_{B_n}} y_n $
		\While{$f_{\mS_k,k}\left(\colvec{\alpha^k \\ w^k} + \gamma\colvec{\Delta \alpha^k \\ \Delta w^k}\right) > (1-2c\gamma)f_{\mS_k,k}(\alpha^k, w^k)$}
			\State $\gamma \gets \beta \cdot \gamma$
		\EndWhile
		\State $w^{k+1}  = w^k + \gamma \Delta w^k $
		\State $\alpha^{k+1}  = \alpha^k +\gamma  \Delta \alpha^k$
	\EndIf
\EndFor
\State \textbf{return:} last iterate $\alpha^k$, $w^k$
\end{algorithmic}
\end{algorithm}

\subsection{Experimental results for stochastic line search}

For all experiments, we set the initial stepsize $\gamma = 2$ with $\gamma$ the stepsize for the last $n$ rows' sketch and reduce the stepsize by a factor $\beta = 0.9$ when the line-search \eqref{eq:stoch_line_search} is not satisfied. We choose the stepsize $\gamma = 1$ with $\gamma$ the stepsize for the first $d$ rows' sketch and $c = 0.09$.

\begin{figure}
\centering
\begin{tabular}{cccc}
\includegraphics[width=.21\linewidth]{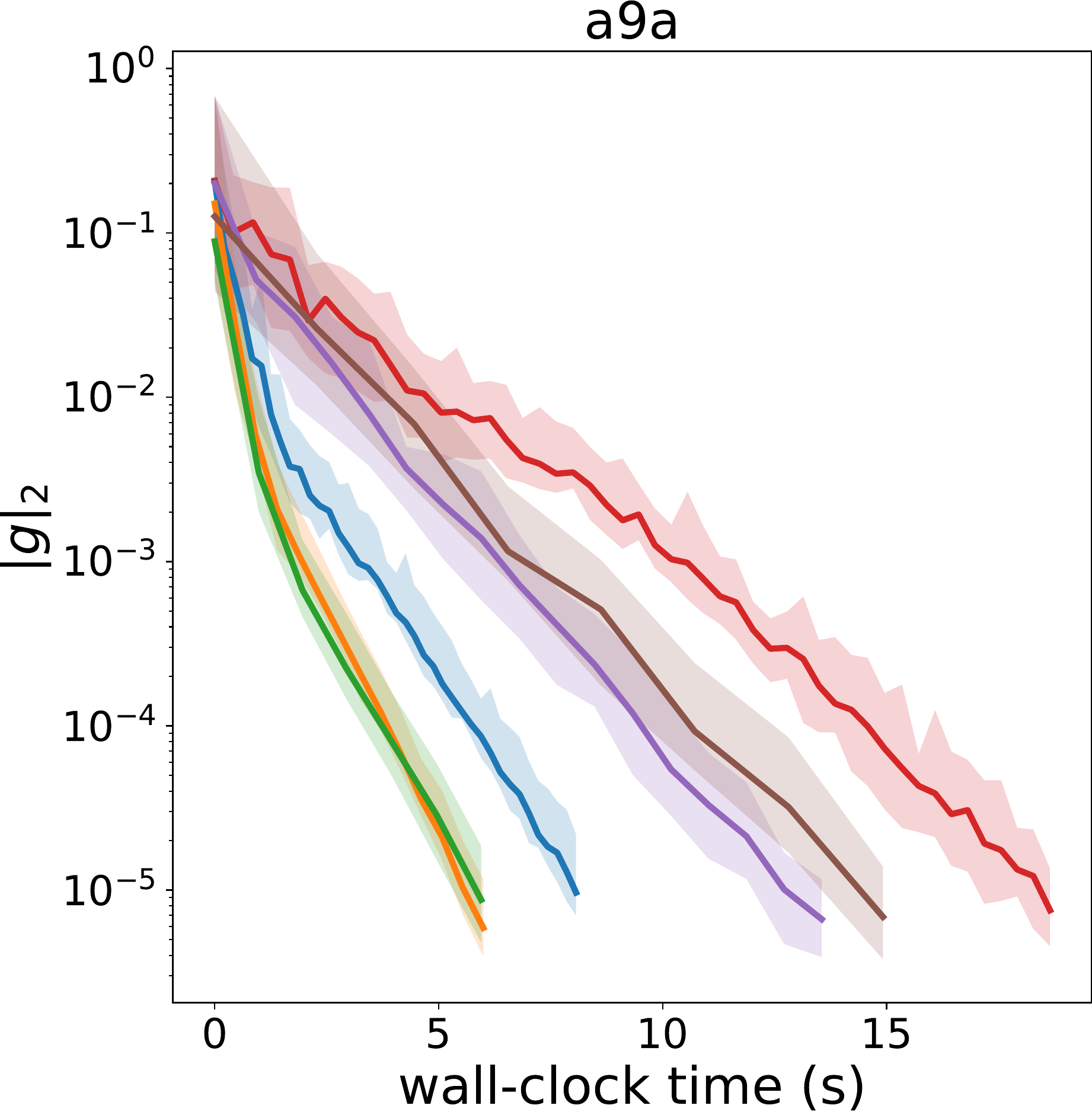}&
\includegraphics[width=.21\linewidth]{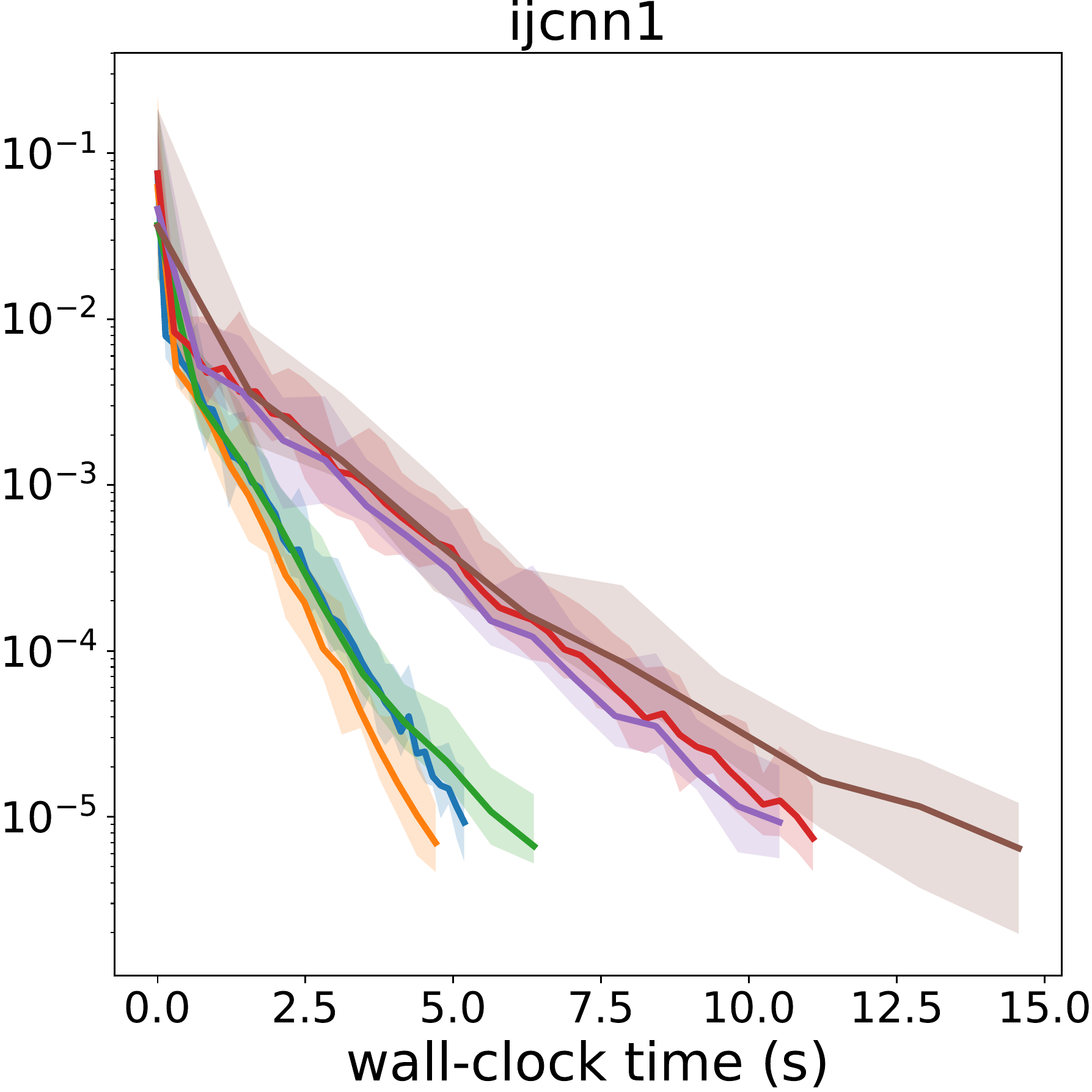}&
\includegraphics[width=.21\linewidth]{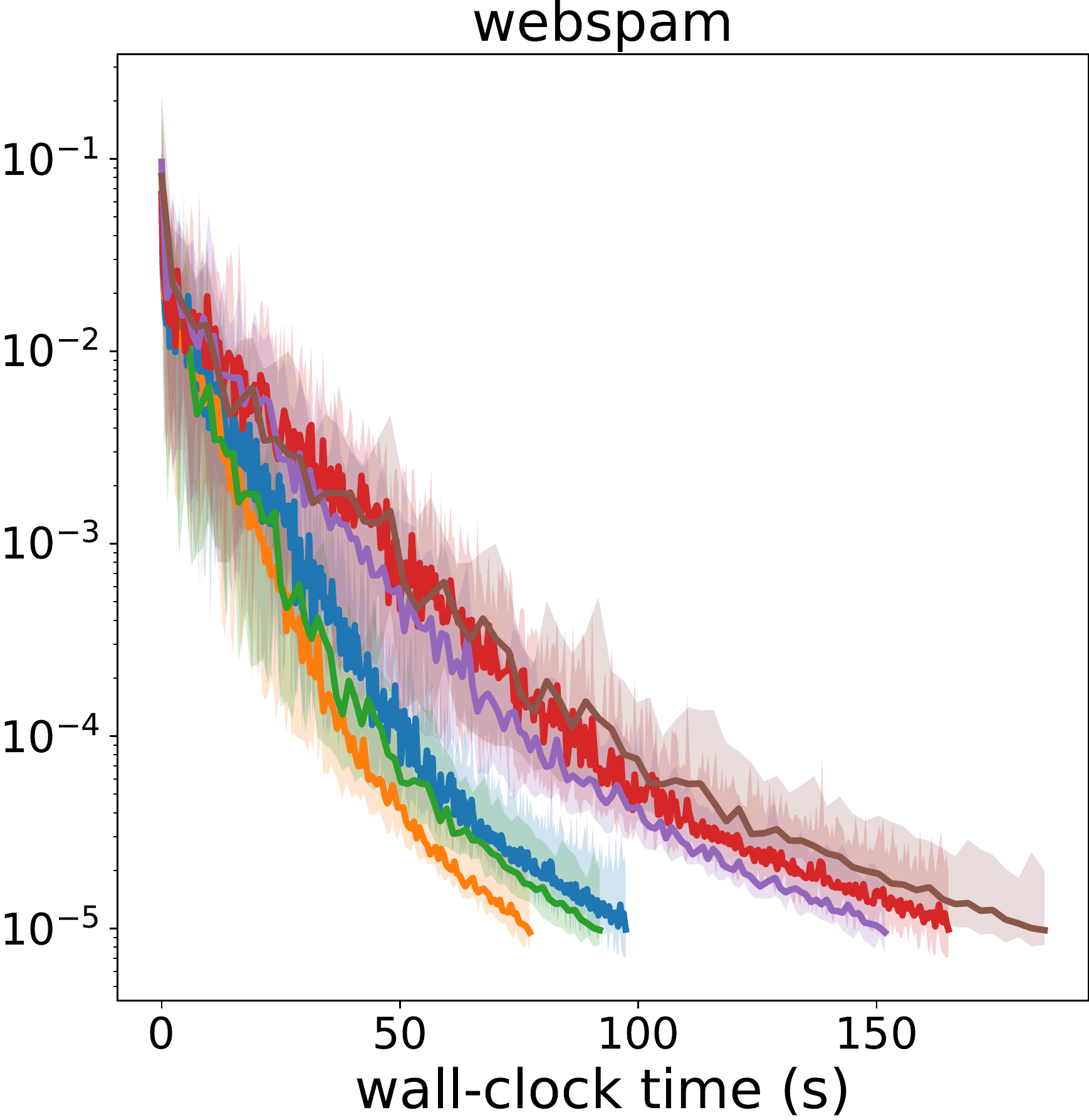}&
\includegraphics[width=.21\linewidth]{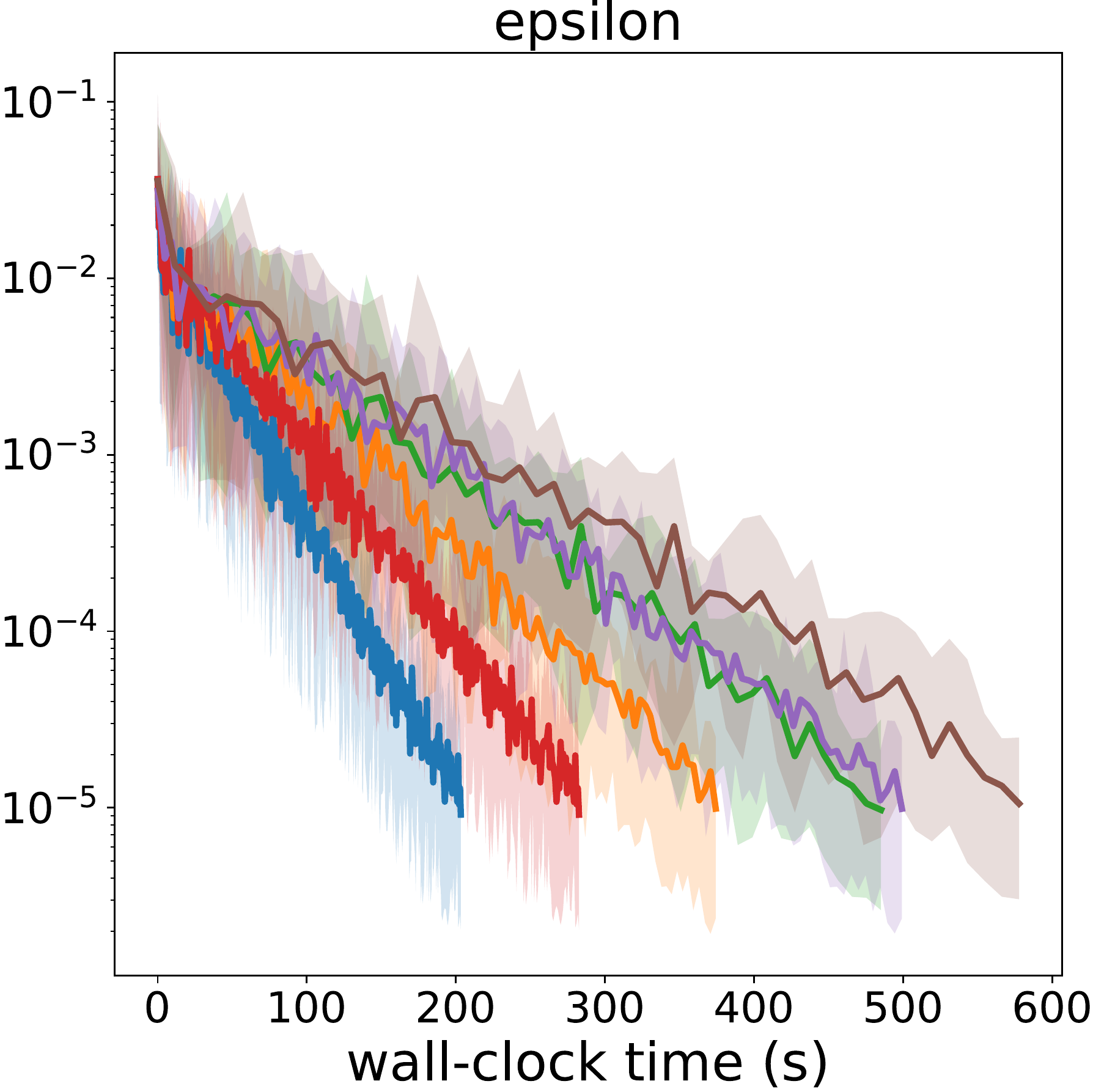} \\
\includegraphics[width=.21\linewidth]{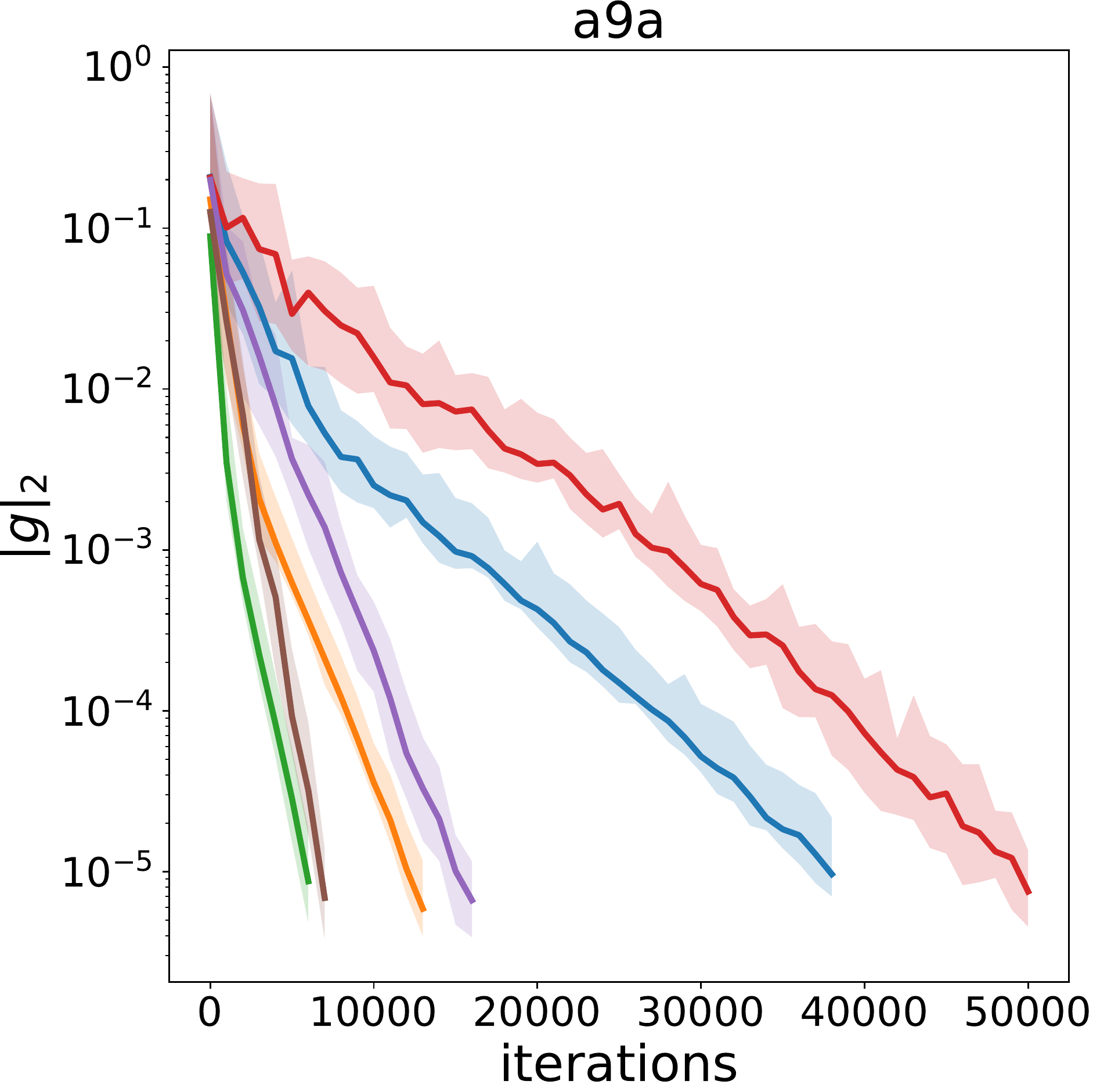}&
\includegraphics[width=.21\linewidth]{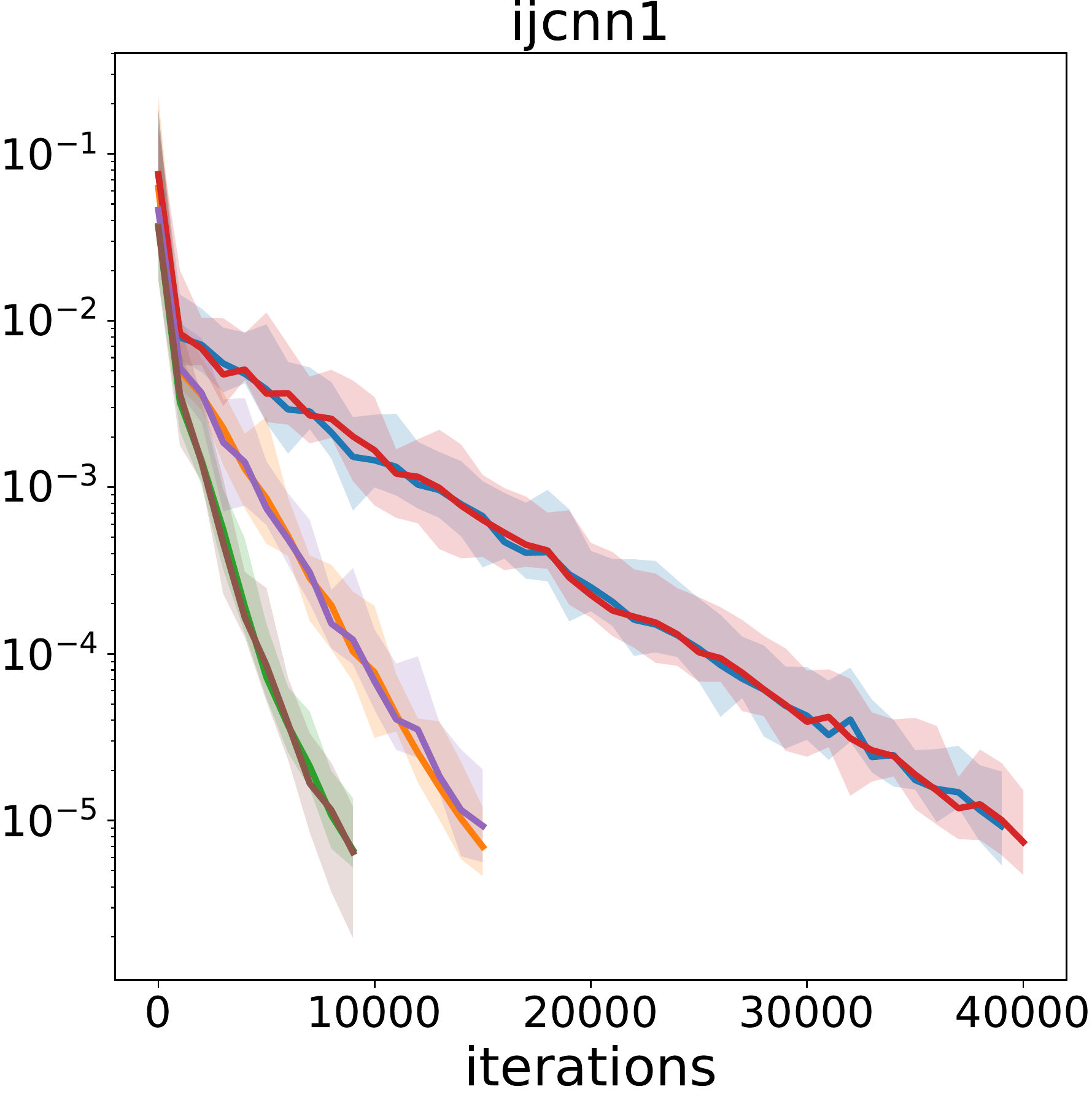}&
\includegraphics[width=.21\linewidth]{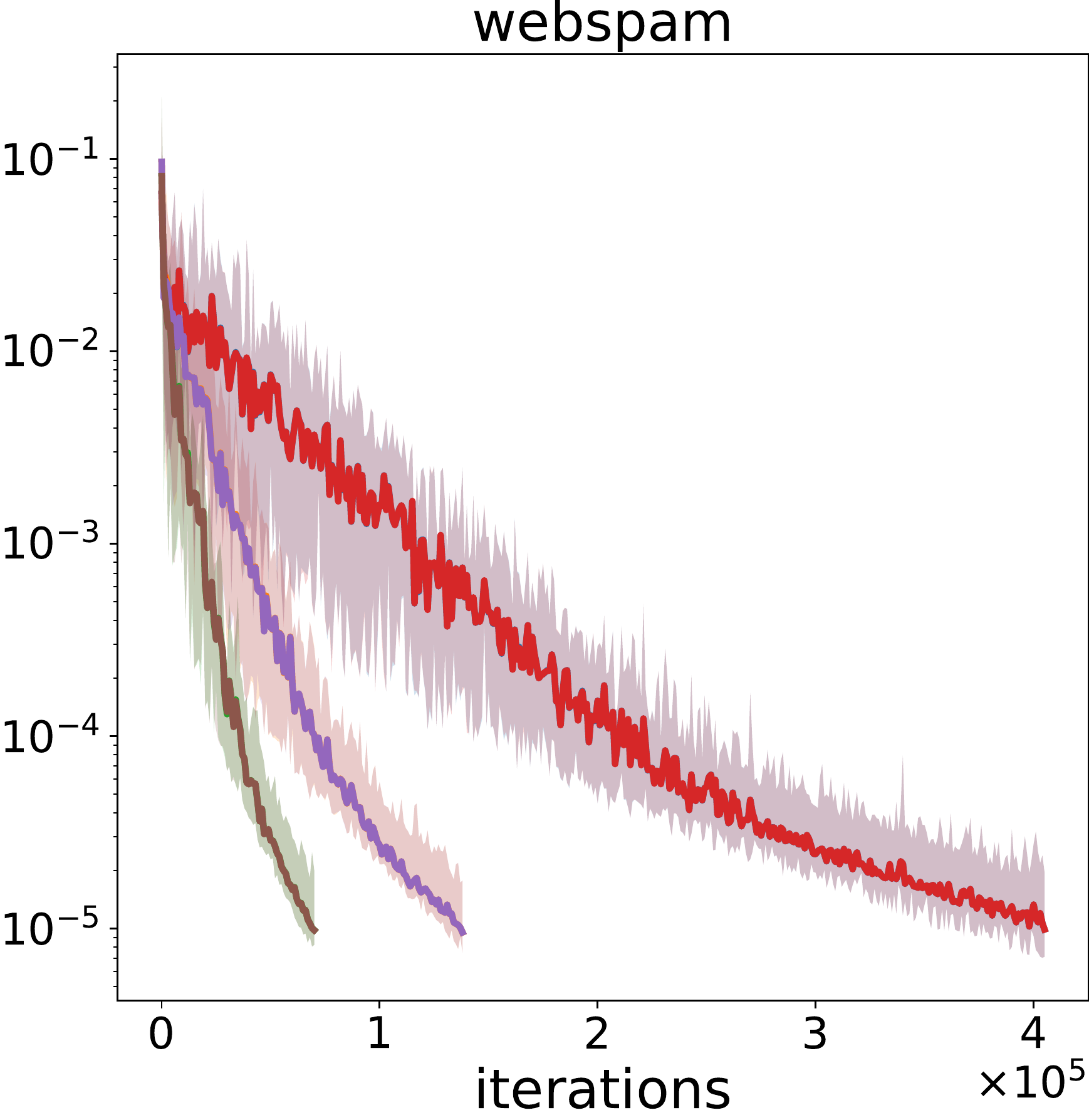}&
\includegraphics[width=.21\linewidth]{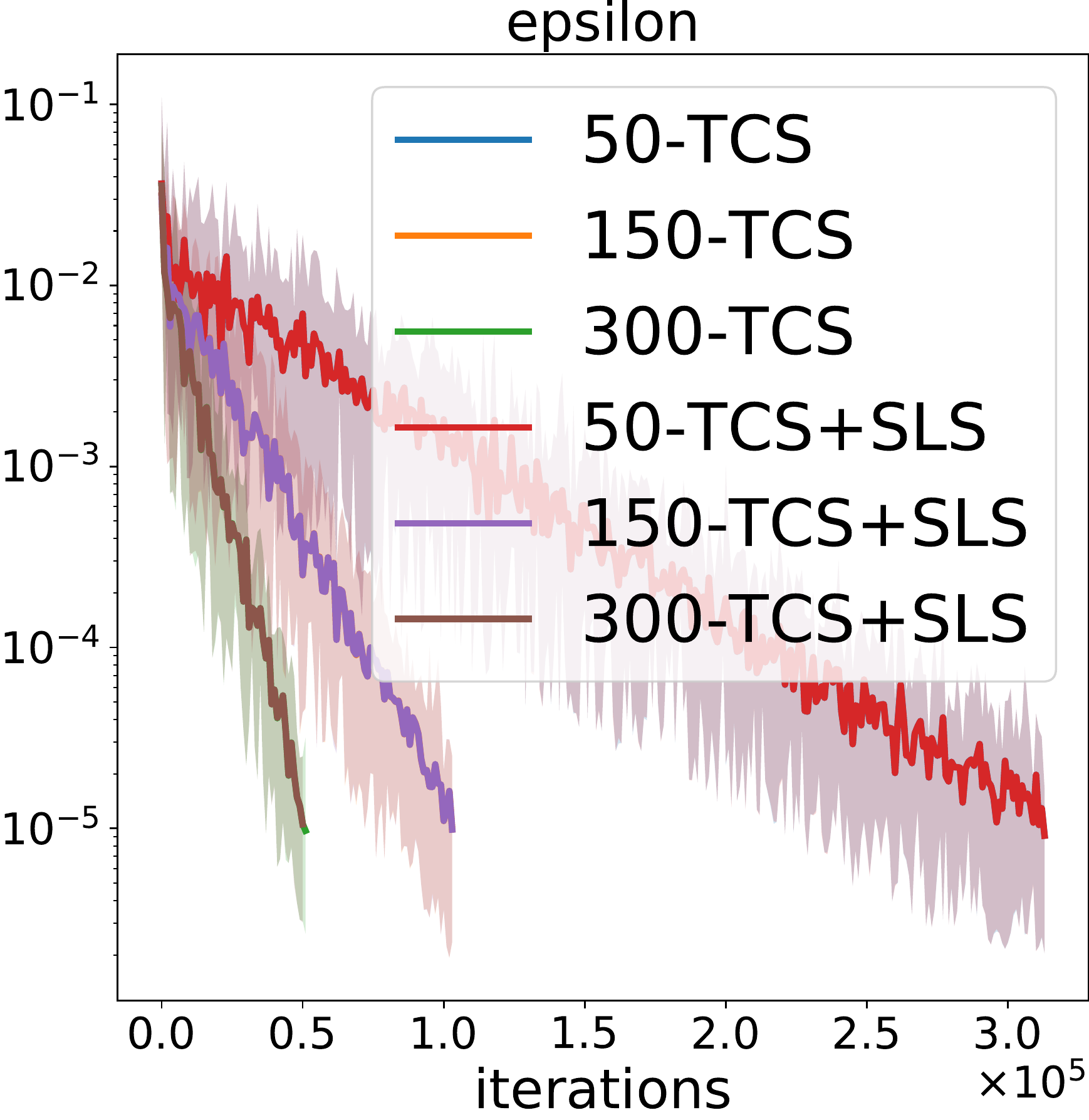}
\end{tabular}
\caption{Experiments for TCS method combined with the stochastic line-search.}
\label{fig:line_search}
\end{figure}

From Figure~\ref{fig:line_search}, we observe that stochastic line search guarantees the convergence of the algorithm and does not tune any parameters. However, it slows down the convergence speed compared to the original algorithm with its rule of thumb parameters' choice. This is expected, as it does extra function evaluations at each step for the stochastic line search procedure.

\bibliographystyle{siamplain}
\bibliography{references}

\end{document}